\documentclass[9pt]{amsart}
\usepackage{amsmath}
\usepackage{amsxtra}
\usepackage{amscd}
\usepackage{amsthm}
\usepackage{amsfonts}
\usepackage{amssymb}
\usepackage[foot]{amsaddr}
\usepackage{cite}
\usepackage{url}
\usepackage{rotating}
\usepackage{eucal}
\usepackage{tikz-cd}
\usepackage[all,2cell]{xy}
\UseAllTwocells
\usepackage{graphicx}
\usepackage{pifont}
\usepackage{comment}
\usepackage{verbatim}
\usepackage{color}
\usepackage{hyperref}
\usepackage{xparse}
\usepackage{upgreek}
\usepackage{MnSymbol}
\newtheorem{cor}[subsubsection]{Corollary}
\newtheorem{lem}[subsubsection]{Lemma}
\newtheorem{goal}[subsubsection]{Goal}
\newtheorem{strategy}[subsubsection]{Strategy}
\newtheorem{prop}[subsubsection]{Proposition}

\newtheorem{warn}[subsubsection]{Warning}

\newtheorem{propconstr}[subsubsection]{Proposition-Construction}

\newtheorem{conj}[subsubsection]{Conjecture}

\newtheorem{thm}[subsubsection]{Theorem}
\newtheorem{defn}[subsubsection]{Definition}
\newtheorem{notn}[subsubsection]{Notation}

\newtheorem{constr}[subsubsection]{Construction}
\newtheorem{rem}[subsubsection]{Remark}
\newtheorem{exam}[subsubsection]{Example}

\newcommand\nc{\newcommand}
\nc\renc{\renewcommand}

\nc\ssec{\subsection}
\nc\sssec{\subsubsection}

\nc\mBA{{\mathbb A}}
\nc\mBB{{\mathbb B}}
\nc\mBC{{\mathbb C}}
\nc\mBD{{\mathbb D}}
\nc\mBE{{\mathbb E}}
\nc\mBF{{\mathbb F}}
\nc\mBG{{\mathbb G}}
\nc\mBH{{\mathbb H}}
\nc\mBI{{\mathbb I}}
\nc\mBJ{{\mathbb J}}
\nc\mBK{{\mathbb K}}
\nc\mBL{{\mathbb L}}
\nc\mBM{{\mathbb M}}
\nc\mBN{{\mathbb N}}
\nc\mBO{{\mathbb O}}
\nc\mBP{{\mathbb P}}
\nc\mBQ{{\mathbb Q}}
\nc\mBR{{\mathbb R}}
\nc\mBS{{\mathbb S}}
\nc\mBT{{\mathbb T}}
\nc\mBU{{\mathbb U}}
\nc\mBV{{\mathbb V}}
\nc\mBW{{\mathbb W}}
\nc\mBX{{\mathbb X}}
\nc\mBY{{\mathbb Y}}
\nc\mBZ{{\mathbb Z}}

\nc\mCA{{\mathcal A}}
\nc\mCB{{\mathcal B}}
\nc\mCC{{\mathcal C}}
\nc\mCD{{\mathcal D}}
\nc\mCE{{\mathcal E}}
\nc\mCF{{\mathcal F}}
\nc\mCG{{\mathcal G}}
\nc\mCH{{\mathcal H}}
\nc\mCI{{\mathcal I}}
\nc\mCJ{{\mathcal J}}
\nc\mCK{{\mathcal K}}
\nc\mCL{{\mathcal L}}
\nc\mCM{{\mathcal M}}
\nc\mCN{{\mathcal N}}
\nc\mCO{{\mathcal O}}
\nc\mCP{{\mathcal P}}
\nc\mCQ{{\mathcal Q}}
\nc\mCR{{\mathcal R}}
\nc\mCS{{\mathcal S}}
\nc\mCT{{\mathcal T}}
\nc\mCU{{\mathcal U}}
\nc\mCV{{\mathcal V}}
\nc\mCW{{\mathcal W}}
\nc\mCX{{\mathcal X}}
\nc\mCY{{\mathcal Y}}
\nc\mCZ{{\mathcal Z}}

\nc\mbA{{\mathbf A}}
\nc\mbB{{\mathbf B}}
\nc\mbC{{\mathbf C}}
\nc\mbD{{\mathbf D}}
\nc\mbE{{\mathbf E}}
\nc\mbF{{\mathbf F}}
\nc\mbG{{\mathbf G}}
\nc\mbH{{\mathbf H}}
\nc\mbI{{\mathbf I}}
\nc\mbJ{{\mathbf J}}
\nc\mbK{{\mathbf K}}
\nc\mbL{{\mathbf L}}
\nc\mbM{{\mathbf M}}
\nc\mbN{{\mathbf N}}
\nc\mbO{{\mathbf O}}
\nc\mbP{{\mathbf P}}
\nc\mbQ{{\mathbf Q}}
\nc\mbR{{\mathbf R}}
\nc\mbS{{\mathbf S}}
\nc\mbT{{\mathbf T}}
\nc\mbU{{\mathbf U}}
\nc\mbV{{\mathbf V}}
\nc\mbW{{\mathbf W}}
\nc\mbX{{\mathbf X}}
\nc\mbY{{\mathbf Y}}
\nc\mbZ{{\mathbf Z}}

\nc\mba{{\mathbf a}}
\nc\mbb{{\mathbf b}}
\nc\mbc{{\mathbf c}}
\nc\mbd{{\mathbf d}}
\nc\mbe{{\mathbf e}}
\nc\mbf{{\mathbf f}}
\nc\mbg{{\mathbf g}}
\nc\mbh{{\mathbf h}}
\nc\mbi{{\mathbf i}}
\nc\mbj{{\mathbf j}}
\nc\mbk{{\mathbf k}}
\nc\mbl{{\mathbf l}}
\nc\mbm{{\mathbf m}}
\nc\mbn{{\mathbf n}}
\nc\mbo{{\mathbf o}}
\nc\mbp{{\mathbf p}}
\nc\mbq{{\mathbf q}}
\nc\mbr{{\mathbf r}}
\nc\mbs{{\mathbf s}}
\nc\mbt{{\mathbf t}}
\nc\mbu{{\mathbf u}}
\nc\mbv{{\mathbf v}}
\nc\mbw{{\mathbf w}}
\nc\mbx{{\mathbf x}}
\nc\mby{{\mathbf y}}
\nc\mbz{{\mathbf z}}

\nc\mfa{{\mathfrak a}}
\nc\mfb{{\mathfrak b}}
\nc\mfc{{\mathfrak c}}
\nc\mfd{{\mathfrak d}}
\nc\mfe{{\mathfrak e}}
\nc\mff{{\mathfrak f}}
\nc\mfg{{\mathfrak g}}
\nc\mfh{{\mathfrak h}}
\nc\mfi{{\mathfrak i}}
\nc\mfj{{\mathfrak j}}
\nc\mfk{{\mathfrak k}}
\nc\mfl{{\mathfrak l}}
\nc\mfm{{\mathfrak m}}
\nc\mfn{{\mathfrak n}}
\nc\mfo{{\mathfrak o}}
\nc\mfp{{\mathfrak p}}
\nc\mfq{{\mathfrak q}}
\nc\mfr{{\mathfrak r}}
\nc\mfs{{\mathfrak s}}
\nc\mft{{\mathfrak t}}
\nc\mfu{{\mathfrak u}}
\nc\mfv{{\mathfrak v}}
\nc\mfw{{\mathfrak w}}
\nc\mfx{{\mathfrak x}}
\nc\mfy{{\mathfrak y}}
\nc\mfz{{\mathfrak z}}

\NewDocumentCommand{\ot}{e{_^}}{
  \mathbin{\mathop{\otimes}\displaylimits
    \IfValueT{#1}{_{#1}}
    \IfValueT{#2}{^{#2}}
  }
}
\NewDocumentCommand{\boxt}{e{_^}}{
  \mathbin{\mathop{\boxtimes}\displaylimits
    \IfValueT{#1}{_{#1}}
    \IfValueT{#2}{^{#2}}
  }
}
\NewDocumentCommand{\mt}{e{_^}}{
  \mathbin{\mathop{\times}\displaylimits
    \IfValueT{#1}{_{#1}}
    \IfValueT{#2}{^{#2}}
  }
}
\NewDocumentCommand{\convolve}{e{_^}}{
  \mathbin{\mathop{\star}\displaylimits
    \IfValueT{#1}{_{#1}}
    \IfValueT{#2}{^{#2}}
  }
}
\NewDocumentCommand{\colim}{e{_^}}{
  \mathbin{\mathop{\operatorname{colim}}\displaylimits
    \IfValueT{#1}{_{#1}\,}
    \IfValueT{#2}{^{#2}\,}
  }
}
\NewDocumentCommand{\laxlim}{e{_^}}{
  \mathbin{\mathop{\operatorname{laxlim}}\displaylimits
    \IfValueT{#1}{_{#1}\,}
    \IfValueT{#2}{^{#2}\,}
  }
}
\NewDocumentCommand{\oplaxlim}{e{_^}}{
  \mathbin{\mathop\operatorname{oplax-lim}\displaylimits
    \IfValueT{#1}{_{#1}\,}
    \IfValueT{#2}{^{#2}\,}
  }
}

\newcommand{\xrightleftarrows}[1]{\mathrel{\substack{\xrightarrow{#1} \\[-.9ex] \xleftarrow{#1}}}} 
 
\renewcommand{\subset}{\subseteq}
\renewcommand{\supset}{\supseteq}

\DeclareMathOperator{\Id}{Id} 
\DeclareMathOperator{\coFib}{coFib} 
\DeclareMathOperator{\Rank}{Rank} 
\DeclareMathOperator{\Maps}{Maps} 
\DeclareMathOperator{\LFun}{LFun} 
\DeclareMathOperator{\Ind}{Ind} 
\DeclareMathOperator{\Arr}{Arr} 
\DeclareMathOperator{\Tw}{Tw} 
\DeclareMathOperator{\Corr}{Corr} 
\DeclareMathOperator{\adjoint}{\xrightleftarrows{\rule{0.5cm}{0cm}} } 
\DeclareMathOperator{\Ad}{Ad} 
\DeclareMathOperator{\oblv}{oblv} 
\DeclareMathOperator{\ins}{ins} 
\DeclareMathOperator{\ev}{ev} 
\DeclareMathOperator{\Av}{Av} 
\DeclareMathOperator{\DGCat}{DGCat} 
\DeclareMathOperator{\hmod}{-mod} 

\DeclareMathOperator{\PreStk}{PreStk} 
\DeclareMathOperator{\Spec}{Spec} 
\DeclareMathOperator{\pt}{pt} 
\DeclareMathOperator{\Div}{Div}  
\DeclareMathOperator{\Sing}{Sing} 
\newcommand{\twisttimes}{\widetilde{\times}} 

\DeclareMathOperator{\Lie}{Lie} 
\DeclareMathOperator{\Rep}{Rep} 
\DeclareMathOperator{\Par}{Par} 
\DeclareMathOperator{\Nilp}{Nilp} 
\DeclareMathOperator{\Sat}{Sat} 
\DeclareMathOperator{\Whit}{Whit} 
\DeclareMathOperator{\Eis}{Eis} 
\DeclareMathOperator{\CT}{CT} 
\DeclareMathOperator{\Bun}{Bun} 
\DeclareMathOperator{\Gr}{Gr} 
\DeclareMathOperator{\Glue}{Glue} 
\DeclareMathOperator{\Sph}{Sph} 
\DeclareMathOperator{\SI}{SI} 
\DeclareMathOperator{\Hecke}{Hecke} 
\DeclareMathOperator{\Zas}{Zas} 
\DeclareMathOperator{\PsId}{Ps-Id} 
\DeclareMathOperator{\HL}{HL} 

\DeclareMathOperator{\Coh}{Coh} 
\DeclareMathOperator{\QCoh}{QCoh} 
\DeclareMathOperator{\IndCoh}{IndCoh} 
\DeclareMathOperator{\DMod}{DMod} 

\newcommand{\op}{\mathrm{op}} 
\newcommand{\st}{\mathrm{st}} 
\newcommand{\rev}{\mathrm{rev}} 
\newcommand{\rel}{\mathrm{rel}} 
\newcommand{\ext}{\mathrm{ext}} 
\newcommand{\enh}{\mathrm{enh}} 
\newcommand{\ren}{\mathrm{ren}} 
\newcommand{\pos}{\mathrm{pos}} 
\newcommand{\level}{\mathrm{level}} 
\newcommand{\temp}{\mathrm{temp}} 
\newcommand{\atemp}{\mathrm{atemp}} 
\newcommand{\htemp}{\mathrm{-temp}} 
\newcommand{\hatemp}{\mathrm{-atemp}} 
\newcommand{\hgen}{{\mathrm{-gen}}} 
\newcommand{\disj}{\mathrm{disj}} 
\newcommand{\hyphen}{\mathrm{-}}
\newcommand{\hpull}{\mathrm{-pull}} 
\newcommand{\hpush}{\mathrm{-push}} 

\newcommand{\Dmod}{\DMod} 
\newcommand{\ICoh}{\IndCoh}
\newcommand{\virg}[1]{``#1"}
\DeclareMathOperator{\LS}{LS}

\nc{\LSGch}{{\LS_\Gch}}
\nc{\LSPch}{{\LS_\Pch}}
\nc{\LSMch}{{\LS_\Mch}}

\nc{\Gch}{{\check{G}}}
\nc{\Pch}{{\check{P}}}
\nc{\Mch}{{\check{M}}}
\nc{\Qch}{{\check{Q}}}

\nc{\Nch}{\check{\mathcal{N}}}

\nc{\LL}{\mathbb{L}}

\nc{\kk}{\mathbbm{k}}

\nc{\longto}{\longrightarrow}

\DeclareMathOperator{\Poinc}{Poinc}
\DeclareMathOperator{\ch}{ch}
\DeclareMathOperator{\Fun}{Fun}
\DeclareMathOperator{\coeff}{coeff}

\nc{\PQ}{P \supset Q}

\begin{document}

\newpage
\title[Automorphic Gluing]{Automorphic Gluing}
\author{Dario Beraldo and Lin Chen}
\begin{abstract}
We prove a gluing theorem on the automorphic side of the geometric Langlands correspondence: roughly speaking, we show that the difference between $\Dmod(\Bun_G)$ and its full subcategory $\Dmod(\Bun_G)^\temp$ of tempered objects is compensated by the categories $\Dmod(\Bun_M)^\temp$ for all standard Levi subgroups $M \subsetneq G$. 
This theorem is designed to match exactly with the spectral gluing theorem, an analogous result occurring on the other side of the geometric Langlands conjecture.
Along the way, we state and prove several facts that might be of independent interest. For instance, for any parabolic $P \subseteq G$, we show that the functors $\CT_{P,*}:\Dmod(\Bun_G) \to \Dmod(\Bun_M)$ and $\Eis_{P,*}:\Dmod(\Bun_M) \to \Dmod(\Bun_G)$ preserve tempered objects, whereas the standard Eisenstein functor $\Eis_{P,!}$ preserves anti-tempered objects.
\end{abstract}

\maketitle
\tableofcontents

\section{Introduction}

This paper is a contribution to the geometric Langlands program: we formulate and prove a gluing theorem occurring on the automorphic side of the geometric Langlands conjecture. We believe this theorem is important for at least three different reasons: 
\begin{itemize}
\item
it is a direct geometric analogue of Langlands' classification theorem for representations of $G(\mathbb{R})$;
\item
it bypasses the usage of the \emph{extended Whittaker category}: the latter is a main ingredient in Gaitsgory's outline of the proof of geometric Langlands, but it has been proven to be hard to deal with beyond the cases of $G=GL_n$ and $G=PGL_n$;
\item
it is designed to match precisely with the \emph{spectral gluing theorem}, the gluing theorem occurring on the other side of the geometric Langlands duality.
As far as we are aware, the only way to construct a functor that connects the two sides of the duality employs these two gluing theorems in an essential way.
\end{itemize}

Moreover, some of our techniques are quite general (second adjunctions, Weyl combinatorics, the semi-infinite category, miraculous parabolic duality, $\ldots$) and we expect them to become useful in other contexts.

\ssec{Rough statement of the main result}


Let $G$ be a connected reductive group and $X$ a smooth projective curve, both defined over an algebraically closed ground field $k$ of characteristic zero. The automorphic side of the geometric Langlands conjecture is $\Dmod(\Bun_G)$, the differential graded (DG) category of D-modules on the stack $\Bun_G$ of $G$-bundles on $X$. 
It turns out that $\Dmod(\Bun_G)$ contains a remarkable full subcategory: the subcategory $\Dmod(\Bun_G)^\temp$ of tempered objects. Tempered objects were introduced in \cite[Section 12]{arinkin2015singular} and studied in \cite{beraldo2021tempered}, \cite{beraldo2021geometric}, \cite{faergeman2021arinkin}.

\begin{exam}
The general definition of $\Dmod(\Bun_G)^\temp$ will be given later: there are several different ways to present it, but they all involve the Hecke action as a common ingredient. The case of $G$ of semisimple rank one is an exception: for instance, an object $\mCF \in \Dmod(\Bun_{SL_2})$ is tempered iff its cohomology with compact supports vanishes, see \cite{beraldo2021geometric}. The \virg{only if} direction is true for general $G$ by \cite{beraldo2021tempered}. 
\end{exam}


Our \emph{automorphic gluing theorem} explains how $\Dmod(\Bun_G)$ is built out of tempered objects for $G$ and for its Levi subgroups. 
Very roughly and heuristically, we may write
\begin{equation} \label{eqn:rough-aut-gluing}
\Dmod(\Bun_G)
\simeq
\underset{P \in \Par}\Glue \,
\Dmod(\Bun_M)^\temp,
\end{equation}
where $\Par$ denotes the poset of standard parabolics subgroups of $G$ (relative to a Borel subgroup $B$ fixed once and for all) and $M = \operatorname{Levi}(P)$. Of course, all the technicalities and subtleties of the theorem are hidden in the symbol \virg{$\Glue$}. We will come back to them later, in Section \ref{ssec:intro-construction-gluing}.

\begin{rem}
Here we see the resemblance between the above statement and the description of $G(\mathbb{R})$-representations due to Langlands. 
\end{rem}


\ssec{Geometric Langlands}

Now, let us recall the statement of the geometric Langlands conjecture and explain the point where the automorphic theorem will be exploited.


 Consider the Langlands dual group $\Gch$ and the stack $\LSGch$ of $\Gch$-local systems on $X$. This stack is derived in general and mildly singular (quasi-smooth).
The spectral side of the geometric Langlands conjecture is the DG category $\ICoh_{\Nch}(\LSGch)$ of \emph{ind-coherent sheaf with nilpotent singular support} on $\LSGch$.

In their ground-breaking work \cite{arinkin2015singular}, Arinkin and Gaitsgory defined $\ICoh_{\Nch}(\LSGch)$ and argued that
\begin{equation} \label{eqn:GLC}
\ICoh_{\Nch}(\LSGch)
\stackrel ? \simeq
\Dmod(\Bun_G)
\end{equation}
is a plausible formulation of the geometric Langlands conjecture. This formulation of the conjecture prompts the definition of tempered D-modules: by construction, $\ICoh_{\Nch}(\LSGch)$ contains the more familiar $\QCoh(\LSGch)$ as a full subcategory, and such a subcategory ought to correspond to $\Dmod(\Bun_G)^{\temp}$ under geometric Langlands.


In a second break-through (\cite{arinkin2018category}), Arinkin and Gaitsgory showed that $\ICoh_{\Nch}(\LSGch)$  can be embedded fully faithfully into a DG category glued out of $\QCoh(\LSMch)$, where $\Mch$ runs thorugh all standard Levi subgroups of $\Gch$ (including $\Gch$ itself):
\begin{equation} \label{eqn:AG-spectral-gluing}
\ICoh_{\Nch}(\LSGch)
\hookrightarrow
\underset{P \in \Par}{\Glue'} \,
\QCoh(\LSMch).
\end{equation}
Again, the technicalities are all hidden in the meaning of the symbol \virg{$\Glue'$}.
However, this type of gluing is not optimal, in that it does not yield an equivalence but only a fully faithful functor.


Building up on \cite{arinkin2018category}, in \cite{beraldo2020spectral}, the first author reglued the DG categories on the RHS of \eqref{eqn:AG-spectral-gluing} to obtain an equivalence. Roughly, the main theorem of \cite{beraldo2020spectral} yields an equivalence
\begin{equation} \label{eqn:rough-spectral-gluing}
\ICoh_{\Nch}(\LSGch)
\simeq
\underset{P \in \Par}{\Glue} \,
\QCoh(\LSMch).
\end{equation}

\sssec{}

We can now explain the importance of \eqref{eqn:rough-aut-gluing} vis-a-vis \eqref{eqn:rough-spectral-gluing}. 
Namely, we are not aware of a direct way to write a plausible functor between $\ICoh_{\Nch}(\LSGch)$ and $\Dmod(\Bun_G)$. On the other hand, thanks to the \emph{vanishing theorem} of \cite{gaitsgory2010vanishing}, there is an evident functor
$$
\LL_G^{\temp}:
\QCoh(\LSGch) \longto
\Dmod(\Bun_G)^\temp
$$
which ought to be an equivalence.

\begin{rem}

In more detail, $\LL_G^\temp$ is defined as follows. The vanishing theorem states that $\QCoh(\LSGch)$ acts on $\Dmod(\Bun_G)$ preserving $\Dmod(\Bun_G)^\temp$; we set $\LL_G^{\temp}$ to be the result of this action on $\Poinc_!$, a canonical tempered object of $\Dmod(\Bun_G)$.

\end{rem}


The next step of our program (to be performed in another publication) is to prove that the functors $\LL_M^{\temp}$, for all Levi subgroups $M \subset G$, assemble together to yield a functor
$$
\underset{P \in \Par}{\Glue} \,
\QCoh(\LSMch)
\longto
\underset{P \in \Par}\Glue \,
\Dmod(\Bun_M)^\temp,
$$
that is, a functor $\LL_G$ as desired. It will follow from this construction that $\LL_G$ is an equivalence if and only if so is $\LL_M^\temp$, for any Levi subgroup $M$. 

\medskip

In other words, our automorphic gluing theorem allows to reduce the geometric Langlands conjecture to the more concrete (but still very hard) \emph{tempered Langlands conjecture} $\QCoh(\LSGch) \stackrel{?}\simeq
\Dmod(\Bun_G)^\temp$.

\ssec{Bypassing the Whittaker categories}

In \cite{gaitsgory2015outline}, Gaitsgory outlined a strategy to prove the geometric Langlands conjecture. Let us recall this strategy and highlight the role played by our theorem.


In a nutshell, Gaitsgory's strategy amounts to embed both sides of \eqref{eqn:GLC} into a larger DG category and compare the essential images. This larger DG category is the extended Whittaker category $\Whit(G,\ext)$, whose definition we now recall.


It is easier to start with a function-theoretic version. We assume here that $G$ has connected center.
Let $\ch_G$ be the vector space of (additive) characters on $N(\mbA)/N(F)$. We first replace $\Dmod(\Bun_G)$ with the vector space of functions $\Fun(G(F)  \backslash G(\mbA) / G(\mbO))$. Then we consider the vector space of \emph{extended Whittaker functions}
$$
\Fun( G(\mbA) / G(\mbO) \times \ch_G)^{T(F) \ltimes N(\mbA), ev}
$$
consisting of functions on $ G(\mbA) / G(\mbO) \times \ch_G$ satisfying:
\begin{itemize}
\item
$f(n \cdot g,\chi )=\chi(n)\cdot f(g,\chi)$ for all $n \in N(\mbA)$;
\item

$f(t \cdot g,\Ad_t(\chi))= f(g,\chi) $ for all $t\in T(F)$.

\end{itemize}


There is an evident operator 
$$
\coeff_{G,\ext}:
\Fun(G(F) \backslash G(\mbA) / G(\mbO) )
\longto
\Fun(G(\mbA) / G(\mbO) \times \ch_G)^{T(F) \ltimes N(\mbA), \ev}.
$$
Moreover, $\Fun(G(\mbA) / G(\mbO) \times \ch_G)^{T(F) \ltimes N(\mbA), \ev}$ splits into components parametrized by the elements of the coordinate stratification of $\ch_G$.


The standard trick with the mirabolic subgroup shows that $\coeff_{G,\ext}$ is injective when $G=GL_n$ and $G=PGL_n$. In these two cases, we have therefore embedded $\Fun(G(F) \backslash G(\mbA) / G(\mbO) )$ into a simpler vector space. In fact, in spite of the long-winded definition, $\Fun(G(\mbA) / G(\mbO) \times \ch_G)^{T(F) \ltimes N(\mbA), \ev}$ is simpler than $\Fun(G(F) \backslash G(\mbA) / G(\mbO) )$: its definition only involves \emph{abelian} global objects, namely $\ch_G$ and $T(F)$.


The above constructions render to the setting of geometric Langlands as follows (see \cite{beraldo2019whittaker} for details):
 
\begin{itemize}

\item
there is a DG category $\Whit(G,\ext)$, equipped with a functor
$$
\coeff_{G,\ext}: 
\Dmod(\Bun_G) \longto \Whit(G,\ext);
$$

\item
$\Whit(G,\ext)$ naturally decomposes as a \emph{lax limit} of the \emph{partial Whittaker categories}
$$ \Whit(G,\ext) \simeq \laxlim_{P \in \Par} {\Whit(G,P)}; $$

\item
for $G=GL_n$ and $G= PGL_n$, the functor $\coeff_{G,\ext}$ is fully faithful (the proof uses again the mirabolic subgroup, this time combined with a blow-up construction).

\end{itemize}


Let us now summarize Gaitsgory's strategy as presented in \cite{gaitsgory2015outline}: assuming Langlands duality for the Levi subgroups of $G$, one can construct a diagonal fully faithful arrow in the diagram below:
\begin{equation}  \label{diag:slanted}
\xymatrix{
    \ICoh_{\Nch}(\LSGch) \ar[d]^-{\textrm{[spectral gluing]}}_-\simeq 
    & & \Dmod(\Bun_G) \ar[dd]^-{ \coeff_{G,\ext} }\\
    \Glue_P \QCoh(\LSMch) \ar[rrd]^-{\subset} \\
    & & \Whit(G,\ext)
}
\end{equation}

Now, for $G = GL_n$ and $G=PGL_n$, the right vertical arrow is fully faithful too, so one can construct a geometric Langlands equivalence by checking that the essential images of the two functors into $\Whit(G,\ext)$ agree.
Gaitsgory then proposed the following statement (if true, then the same method as above would prove the geometric Langlands conjecture):

\begin{conj}[Gaitsgory]
The functor $\coeff_{G,\ext}$ is fully faithful for any reductive group $G$. 
\end{conj}


This conjecture remains open at the moment. Besides this fundamental issue, working with $\Whit(G,\ext)$ is not ideal for other (less serious) reasons: 
\begin{itemize}

\item
when the center of $G$ is disconnected, the definition must be modified and it becomes more complicated;

\item

as mentioned, $\Whit(G,\ext)$ is only a lax-limit of the partial Whittaker categories, while having an actual limit would be desirable;

\item

we are not aware of a DG category that precisely matches $\Whit(G,\ext)$ on the spectral side (for this reason, the bottom arrow in our previous diagram is slanted).

\end{itemize}


For all these reasons, we propose to bypass $\Whit(G,\ext)$ altogether. It can be shown that $\coeff_{G,\ext}$ naturally factors through the automorphic gluing DG category, thereby improving the symmetry of \eqref{diag:slanted}:

\begin{equation}  
\xymatrix{
    \ICoh_{\Nch}(\LSGch) \ar[d]^-{\textrm{[spectral gluing]}}_-\simeq 
    & & \Dmod(\Bun_G) \ar[d]^-{\textrm{[automorphic gluing]}}_-\simeq \\
    \Glue_P \QCoh(\LSMch) \ar@{.>}[rrd]^-{\subset} \ar[rr] 
    & & \Glue_{P \in \Par} \Dmod(\Bun_M)^\temp \ar@{.>}[d]
    \\
    & & \Whit(G,\ext)
}
\end{equation}

\ssec{Construction of the glued category} \label{ssec:intro-construction-gluing}

Let us now be more precise on how the automorphic gluing theorem is formulated. We will see that our rough statement \eqref{eqn:rough-aut-gluing} must be modified in two ways: first, the components of the gluing category must be made $\Sph_G$-linear; second, we must express their functoriality via twisted arrows of parabolics.


As mentioned, in general the tempered condition defining the full subcategory $\Dmod(\Bun_G)^\temp \subset \Dmod(\Bun_G)$ is expressed in terms of the Hecke action, that is, the action of the spherical monoidal DG category $\Sph_G$ on $\Dmod(\Bun_G)$. 
\footnote{This action requires the choice of $x \in X$ and $\mCO \simeq \mCO_x$: as proven in \cite{faergeman2021arinkin}, any such choice yields the same full subcategory.}
Accordingly, we postulate that all the categories appearing in the automorphic gluing theorem must be equipped with an action of $\Sph_G$. Likewise, all connecting functors appearing in the theorem must be $\Sph_G$-linear.

\begin{rem}
In general, we define $\mCC^{\temp}$ for any DG category $\mCC$ equipped with an action of $\Sph_G$. By definition, $\mCC^{\temp} \subset \mCC$ is a full subcategory. Its right orthogonal is called the full subcategory of \emph{anti-tempered objects} and denoted by $\mCC^{\atemp}$. See Section \ref{ssec-temper} for more details.
\end{rem}


Here we encounter our first modification: for $P \subset G$ a proper parabolic subgroup with Levi quotient $M$, the DG categories $\Dmod(\Bun_M)$ and $\Dmod(\Bun_M)^\temp$ carry action of $\Sph_M$, but not actions of $\Sph_G$. There is a canonical way to obtain $\Sph_G$-linear DG categories from these two, and it leads us to $\mCI(G,P)$ and $\mCW(G,P)$ respectively.


The definition goes as follows (full details will be provided in the section on preliminaries). Following \cite{barlev2012d, chen2020deligne}, consider the prestack $\Bun_G^{P\hgen}$ parametrized $G$-bundles on $X$ with a generic reduction to $P$. 
The DG category $\Dmod(\Bun_G^{P\hgen})$ is well-defined and equipped with a natural action of $\Sph_G$.

\begin{rem}
Similarly to our description of the Whittaker category, $\Bun_G^{P\hgen}$ is a geometric analogue of $P(F) \backslash G(\mCA)/G(\mCO)$ and $\Dmod(\Bun_G^{P\hgen})$ is analogue to the vector space $\Fun(P(F) \backslash G(\mCA)/G(\mCO))$. 
\end{rem}


Considering the obvious maps
\begin{equation} \label{eqn:fancy-CT-span}
\Bun_M \xleftarrow{q} \Bun_P \xrightarrow{i_P} \Bun_G^{P\hgen},
\end{equation}
we let $\mCI(G,P)$ and $\mCW(G,P)$ be the DG full subcategories of $\Dmod(\Bun_G^{P\hgen})$ defined by the following fiber diagrams:
\[
\xymatrix{
	\mCW(G,P) \ar[r] \ar[d] &
	\mCI(G,P) \ar[r] \ar[d] &
	\DMod(\Bun_G^{P\hgen}) \ar[d]^-{\iota_P^!} \\
	\DMod(\Bun_M)^{\temp} \ar[r] &
	\DMod(\Bun_M) \ar[r]^-{q^!} &
	\DMod(\Bun_P)
}
\]


It is easy to see that $\Dmod(\Bun_G^{P\hgen})$ is equipped with an action of $\Sph_G$.
However, it is not at all clear from the above definition that the DG categories $\mCI(G,P)$ and $\mCW(G,P)$ admit an action of $\Sph_G$. This is one of our first technical result, proven in Proposition \ref{prop-sph-acts-on-IGP}:

\begin{prop} \label{prop:Sph-action on I and W}
The natural $\Sph_G$-action on $\Dmod(\Bun_G^{P\hgen})$ preserves $\mCI(G,P)$ and $\mCW(G,P)$.
\end{prop}

\begin{rem}
We can compare with the situation on the spectral side: there we see that $\QCoh(\LSMch)$ and $\ICoh_{\Nch_M}(\LSMch)$ carry actions of $\mBH(\LSMch)$ but not actions of $\mBH(\LSGch)$. Using $\Pch$, there is a canonical way to produce two $\mBH(\LSGch)$-linear DG categories; for details, see the introduction of \cite{beraldo2020spectral}. 
\end{rem}

\begin{rem}
This comparison brings up a recurring theme: while the spectral gluing theorem involves the \emph{formal geometry} of stacks of local systems, the automorphic theorem involves the \emph{rational geometry} of the stack of bundles.
\end{rem}


The DG categories $\mCW(G,P)$ are the $\Sph_G$-linear categories that replace the $\Dmod(\Bun_M)^\temp$ in our very preliminary statement automorphic gluing theorem. Thus, our rough statement \eqref{eqn:rough-aut-gluing} now reads
\begin{equation} \label{eqn:less-rough-aut-gluing}
\Dmod(\Bun_G)
\simeq
\underset{P \in \Par}\Glue \,
\mCW(G,P),
\end{equation}
but once again the symbol $\Glue$ carried a vague meaning, which we are going to clarify next. This will be our second and last modification, alluded to at the beginning of this section.



To arrive at a precise statement, we need to investigate how the components $\mCW(G,P)$ and $\mCW(G,Q)$ are related for $P \supset Q$. 
The key idea is to realize that the assignment $P \mapsto \mCW(G,P)$ can be extended to a functor out of twisted arrows of parabolics. Namely, for any inclusion $P \supset Q$, we define a DG category $\mCW(G,\PQ)$ as follows. Let $M$ denote the Levi quotient(/subgroup) of $P$ and consider the following natural diagram:\footnote{That this diagram generalizes \eqref{eqn:fancy-CT-span}: the latter is a special case of the former for $P=Q$.}
\begin{equation} \label{eqn:very-fancy-CT-span}
\Bun_M^{Q \cap M \hgen} \xleftarrow{ } \Bun_P^{Q\hgen} \xrightarrow{ } \Bun_G^{Q\hgen}.
\end{equation}
Pull-push along this diagram yields a functor
$$
\Dmod(\Bun_G^{Q\hgen}) \longto \Dmod(\Bun_M^{Q \cap M \hgen})
$$
that, as shown in the main body of the text, restricts to a functor
$$
\mCI(G,Q) \longto \mCI(M,Q \cap M).
$$
Since $ \mCI(M,Q \cap M)$ is equipped with an action of $\Sph_M$, it makes sense to consider its $M$-tempered subcategory: we then define $\mCW(G, \PQ)$ by the fiber product
\[
\xymatrix{
	\mCW(G,\PQ) \ar[r]\ar[d] & \mCI(G,Q)\ar[d] \\
	\mCI(M,Q\cap M)^{\temp} \ar[r] & \mCI(M,Q\cap M).
}
\]

\begin{rem}
When making definitions like the one above, we need to appeal to tempered objects relative to different reductive groups ($G$ or $M$). Hence, to avoid confusions we may use the notations '$G\mathrm{-}\temp$' or '$M\mathrm{-}\temp$'. For instance, in the above diagram, we might have used $\mCI(M,Q \cap M)^{M\mathrm{-}\temp}$ to signal that tempered objects of $I(M,Q \cap M)$ are taken with respect to the $\Sph_M$-action on $\mCI(M,Q \cap M)$ given by Proposition \ref{prop:Sph-action on I and W}.
\end{rem}


In Proposition \ref{prop-W-preserve-by-Sph}, we will prove that the $\Sph_G$-action on $\mCI(G,Q)$ preserves the full subcategory $\mCW(G,\PQ)$.
Now observe that, for $P=Q$, the DG category $\mCW(G,\PQ)$ coincides with our old $\mCW(G,P)$.
In general, for any inclusion $\PQ$, we have a diagram
$$
\mCW(G,P) 
\leftarrow
\mCW(G,\PQ) 
\rightarrow
\mCW(G,Q)
$$ 
connecting $\mCW(G,P) $ with $\mCW(G,Q)$. Explicitly:
\begin{itemize}
\item
The arrow $\mCW(G,\PQ) 
\rightarrow
\mCW(G,Q)$ is fully faithful: it is defined by proving that the structure inclusion $\mCW(G,\PQ) \subset
\mCI(G,Q)$ lands in $\mCW(G,Q)$.

\item

The arrow $\mCW(G,\PQ) \to \mCW(G,P)$ is defined by proving that the composition
$$
\mCW(G,\PQ) \hookrightarrow \Dmod(\Bun_G^{Q\hgen}) \xrightarrow{!\mathit{-push}}
\Dmod(\Bun_G^{P\hgen})
$$
lands in $\mCW(G,P) \subset \Dmod(\Bun_G^{P\hgen})$.

\end{itemize}
By construction, these functors are all $\Sph_G$-linear.


Thus, we see that the formation of $\mCW(G,\PQ)$ is covariant in $P$ and contravariant in $Q$. We express this rigorously by constructing a functor
$$
\mCW(G,- \supset -):
\Tw(\Par)^\op
\longrightarrow
\Sph_G \hmod.
$$
This functoriality is not just a technical exercise with abstract nonsense. One of the main inputs necessary for the construction is the following theorem, which is of interest in its own right.

\begin{thm} \label{thm:CT-temp}
For any $P \in \Par$, the constant term functor $\CT_{P,*}: \Dmod(\Bun_G) \to \Dmod(\Bun_M)$ sends $G$-tempered objects to $M$-tempered objects.
\end{thm}

\begin{rem}
The statement is obvious for $P=B$, but not so for other parabolics. The corresponding statement on the spectral side says that the functor $\CT_{\Pch}^{spec}: \ICoh_{\Nch}(\LSGch) \to \ICoh_{\Nch}(\LSMch)$ sends $\QCoh(\LSGch)$ to $\QCoh(\LSMch)$. This is easy to prove using the rules of propagation of singular support and the fact that the map $\LSPch \to \LSMch$ is quasi-smooth.
\end{rem}


In fact, what we really need is a statement stronger than the above theorem: the fact that the structure functor $\mCI(G,P) \to \Dmod(\Bun_M)$ sends $G$-tempered objects to $M$-tempered objects. This is the first part of Theorem \ref{thm-iota!*-temperedness}.

\begin{rem}
Along the same lines, we prove that the 'non-standard' Eisenstein series functor $\Eis_{P,*} : \DMod(\Bun_M) \to \DMod(\Bun_G)$ sends $M$-tempered objects to $G$-tempered objects, while the standard Eisenstein series functor $\Eis_{P,!} : \DMod(\Bun_M) \to \DMod(\Bun_G)$ sends $M$-anti-tempered objects to $G$-anti-tempered objects. See Theorem \ref{thm-EisCT-temperedness}.
Contrarily to the statement for $\CT_{P,*}$, the corresponding facts on the spectral side are not obvious (for instance, see \cite[\S 1.7]{beraldo2021DL} for a special case).

\end{rem}


It is easy to connect $\Dmod(\Bun_G)$ with $\mCW(G,\PQ)$: in one direction, we simply have a composition of obvious functors
$$
\mCW(G,\PQ)
\hookrightarrow 
\mCI(G,Q)
\hookrightarrow
\Dmod(\Bun_G^{Q\hgen})
\xrightarrow{(p_Q^{\enh})_!}
\Dmod(\Bun_G).
$$
By definition, the composition of the second and third arrow is the \emph{enhanced Eisenstein functor}, see our Section \ref{ssec-EisCT}.
In the main body of the text, we will show these arrows right adjoints: we obtain a functor
$$
\gamma_{G,\PQ}:
\Dmod(\Bun_G)
\longto
\mCW(G,\PQ).
$$
We can now state the automorphic gluing theorem:

\begin{thm} \label{mainthm}
The above arrows yield a $\Sph_G$-linear equivalence
$$
\gamma_G:
\Dmod(\Bun_G)
\xrightarrow{\;\; \simeq \;\;}
\lim_{[\PQ] \in \Tw(\Par)^\op}
\mCW(G,\PQ).
$$
\end{thm}

\begin{rem}

We have already noticed that $\mCW(G,P ) = \mCW(G, P \supset P)$ for any $P \in \Par$. We regard these DG categories as the main components of the gluing; the intermediate ones (that is, the DG categories $\mCW(G, \PQ)$ with $P \neq Q$) are necessary to assemble the main components together.
If one wants to avoid twisted arrows, one could consider only the main components and obtain a gluing statement that is weaker in two respects:  first, these main components assemble only into a lax-limit, not into a limit; second, the natural functor $\Dmod(\Bun_G) \to \laxlim_{P \in \Par} \mCW(G,P)$  that one can write down will be fully faithful, but not an equivalence.
\end{rem}


\ssec{Proving that \texorpdfstring{$\gamma_G$}{automorphic gluing} is an equivalence}

Let us sketch our method to prove that the above automorphic gluing functor $\gamma_G$ is an equivalence. More detailed outlines on the partial steps appear at the beginning of sections \ref{sect-aut-gluing-thm} and \ref{sect-gluing-atemp-proof}.


Since $\gamma_G$ is by construction $\Sph_G$-linear, it induces a functor on the tempered subcategories:
$$
(\gamma_G)^\temp:
\Dmod(\Bun_G)^\temp
\longrightarrow
\lim_{[\PQ] \in \Tw(\Par)^\op}
\mCW(G,\PQ)^\temp.
$$
We first check that this functor is an equivalence: for this, we observe that the structure inclusion $\mCW(G,\PQ)^{\temp} \subset \mCI(G,Q)^{\temp}$ is an equivalence. Since the latter DG category does not depend on $P$, the above limit on the right-hand-side greatly simplifies (twisted arrows disappear) and we are able to conclude.


Next, we look at the subcategories of anti-tempered objects and formulate a different-looking gluing theorem:

\begin{thm} [Anti-tempered gluing theorem]
\label{thm:antitemp-gluing}
There is a natural $\Sph_G$-linear equivalence
$$
\beta_G:
\Dmod(\Bun_G)^\atemp
\longrightarrow
\lim_{R \in (\Par')^\op}
\mCI(G,R)^{\atemp}.
$$
Here $\Par'$ denotes the poset of proper standard parabolics (so, the parabolic $G$ is excluded).
\end{thm}

This theorem looks more tractable than Theorem \ref{mainthm}: the limit is not taken with respect to twisted arrows, and the DG categories $\mCI(G,R)$ are simpler than their relatives $\mCW(G,R)$. We will use induction on the semisimple rank of $G$ to show that Theorem \ref{thm:antitemp-gluing} implies Theorem \ref{mainthm}.


It remains to prove Theorem \ref{thm:antitemp-gluing}. By contruction, $\beta_G$ is equipped with a left adjoint $\beta_G^L$, so we need to show that the unit and counit of the adjunction are both isomorphisms.

The proof for the counit $\beta_G^L \circ \beta_G \to id$ follows by combining the Deligne-Lusztig duality of the second author (see \cite{chen2020deligne}) with the Ramanujan conjecture of the first author (see \cite{beraldo2021geometric}).

Finally, the proof for the unit $id \to \beta_G \circ \beta_G^L$ is one of the core results of this paper and it brings together several techniques: Weyl combinatorics, miraculous duality, Drinfeld compactifications, etc. We refer to the beginning of Section \ref{sect-gluing-atemp-proof} for a detailed outline of the strategy.


\ssec{Structure of the paper}


In Section \ref{sect-preliminaries}, we recall some essential preliminary notions: tempered and anti-tempered objects, Eisenstein series and constant term functors, the prestack $\Bun_G^{P\hgen}$ and the parabolic DG category $\mCI(G,P)$. Then we show that $\Sph_G$ naturally acts on $\mCI(G,P)$; as explained in the discussions above, this action is necessary for our constructions and proofs.


In Section \ref{sect-CT}, we prove another important result, Theorem \ref{thm:CT-temp}, on preservation of tempered objects. Actually, we prove more results of similar kind, see Theorem \ref{thm-EisCT-temperedness} and Theorem \ref{thm-iota!*-temperedness} for the full list. These statements seem to be of interest in their own right. They are proven by passing to a local situation: we consider the \emph{semi-infinite} parabolic category $\SI_P$ and analyze its relation with $\Sph_M$ and $\Sph_G$.


In Section \ref{sect-WGPQ}, we use the results of the previous sections to construct the DG categories $\mCW(G,\PQ)$ and study their functoriality via twisted arrows. We can then define the automorphic gluing functor 
$$
\gamma_G: \Dmod(\Bun_G) \to \lim_{[\PQ] \in \Tw(\Par)^\op} \mCW(G,\PQ)
$$ 
and its left adjoint.


In Section \ref{sect-aut-gluing-thm}, we first state the anti-tempered gluing theorem, that is, the fact that a natural functor
$$
\beta_G: \Dmod(\Bun_G)^{\atemp} \to \lim_{R \in \Par'} \mCI(G,R)^{\atemp}
$$ 
is an equivalence. We use induction on the semisimple rank of $G$ to prove that this theorem implies our main theorem.


Finally, in Section \ref{sect-gluing-atemp-proof}, we prove the anti-tempered gluing theorem, that is, the fact that $\beta_G$ is an equivalence. As mentioned, the most difficult part is to prove that $\beta_G$ is essentially surjective.

\ssec{Acknowledgements}

We are grateful to Dennis Gaitsgory for teaching us most of the mathematics used in this paper. The first author thanks Dima Arinkin, Ian Grojnowski and Sam Raskin for several conversations and suggestions. The second author thanks Kevin Lin for a discussion on Weyl combinatorics.
The work of the second author was supported by a grant from the Simons Foundation (816048, LC) during his stay at the Institute for Advanced Study.

\section{Preliminaries} \label{sect-preliminaries}

In this section, we recall some of the main players of this paper. We fix a parabolic subgroup $P \subseteq G$ throughout. Let $U$ be its unipotent radical and $M=P/U$ be its Levi quotient group.

In \S \ref{ssec-temper}, we recall the notion of \emph{$G$-tempered objects} in a category acted on by $\Sph_G$.
In \S \ref{ssec-EisCT}, we recall the \emph{(geometric) Eisenstein series functors} $\Eis_{P,!}, \Eis_{P,*}$ and the \emph{constant term functors} $\CT_{P,!}, \CT_{P,*}$.
In \S \ref{ssec-IGP}, we recall the \emph{parabolic category} $\mCI(G,P)$ and its properties.
Finally, \S \ref{ssec-IGP-proof-1} is dedicated to the proof of Proposition \ref{prop-sph-acts-on-IGP}: this results equips $\mCI(G,P)$ with a natural action of $\Sph_G$, thereby letting us consider $G$-tempered objects in $\mCI(G,P)$.

\subsection{\texorpdfstring{$G$}{G}-temperedness}
\label{ssec-temper}
In this subsection, we recall the notion of temperedness. Details and proofs can be found in \cite{beraldo2021geometric} and \cite{beraldo2021tempered}.

Consider the monoidal category 
\[\Sph_G:= \DMod( G(\mCO)\backslash G(\mCK)/G(\mCO) )\]
with monoidal structure given by convolution. Recall the \emph{derived Satake equivalence} (see \cite[\S 12]{arinkin2015singular}):
\[\Sat_G: \Sph_G \simeq \IndCoh_{\Nilp(\check \mfg^*)/\check G}( \pt\mt_{\check \mfg}\pt/\check G ) ,\]
where the RHS is a certain full subcategory of $\IndCoh( \pt\mt_{\check \mfg}\pt/\check G)$ that contains $\QCoh(\pt\mt_{\check \mfg}\pt/\check G)$.

\begin{defn}
 Define
\[ \Sph_G^{\temp}:= \Sat_G^{-1}( \QCoh(\pt\mt_{\check \mfg}\pt/\check G)) \subset \Sph_G,\]
which is a two-sided monoidal ideal of $\Sph_G$. Objects in this full subcategory are called \emph{$G$-tempered}, or simply \emph{tempered}.

The category $\Sph_G^{\temp}$ is compactly generated and the inclusion functor $\Sph_G^\temp \to \Sph_G$ preserves compact objects. It has a continuous right adjoint
\[ \temp_G: \Sph_G\to \Sph_G^\temp \]
called the \emph{($G$-)temperization functor}.
\end{defn}

\begin{defn}
Define 
\[ \Sph_G^{\atemp}:= \ker( \temp_G) \subset \Sph_G \]
to be the kernel of $ \temp_G$, which is also a two-sided monoidal ideal of $\Sph_G$. Objects in this full subcategory are called \emph{($G$-)anti-tempered}.

The category $\Sph_G^{\atemp}$ is also compactly generated, and the inclusion functor $\Sph_G^{\atemp}\to \Sph_G$ has a left adjoint
\[ \atemp_G: \Sph_G\to \Sph_G^\atemp .\]
Then we also have 
\[ \Sph_G^{\temp}\simeq \ker( \atemp_G) \subset \Sph_G \]
\end{defn}

Consider the perverse t-structure on $\Sph_G$, i.e., that induced from the perverse t-structure on $\DMod(\Gr_G)$. We have

\begin{lem}[\!\!{\cite[Corollary 2.3.7]{beraldo2021geometric}}] \label{lem-atemp-t-structure} The full subcategory $ \Sph_G^{\atemp}\subset \Sph_G$ is equivalent to $\Sph_G^{\le -\infty}$.
\end{lem}

\begin{defn} For any left $\Sph_G$-module category $\mCC$, define
\[ \mCC^{G\mathrm{-}\temp}:= \Sph_G^{\temp} \ot_{\Sph_G} \mCC,\;\mCC^{G\mathrm{-}\atemp}:= \Sph_G^{\atemp} \ot_{\Sph_G} \mCC .\]
The $\Sph_G$-linear adjoint functors
\[ \Sph_G^{\temp} \adjoint \Sph_G \adjoint \Sph_G^{\atemp} \]
induce adjoint functors
\[ \mCC^{G\mathrm{-}\temp} \adjoint \mCC \adjoint \mCC^{G\mathrm{-}\atemp}.\]
We view $\mCC^{G\mathrm{-}\temp}$ and $\mCC^{G\mathrm{-}\atemp}$ as full subcategories of $\mCC$: these are the full subcategories of \emph{($G$-)tempered} and  \emph{($G$-)anti-tempered} objects of $\mCC$. We still have
\[ \mCC^{G\hatemp}\simeq \ker( \temp_G) \subset \mCC ,\;\;  \mCC^{G\htemp}\simeq \ker( \atemp_G) \subset \mCC\]

We define similar notions for right $\Sph_G$-module categories. Equivalently, they can be defined by viewing any right $\Sph_G$-module category as a left one via the anti-monoidal involution on $\Sph_G$ induced by taking the inverse in $G(\mCO)\backslash G(\mCK)/G(\mCO)$.
\end{defn}

\begin{defn} Define
\[ \mathbf{1}_{\Sph_G}^{\temp}:= \temp_G(\mathbf{1}_{\Sph_G}) \]
to be the temperization of the monoidal unit of $\Sph_G$. We call this object the \emph{tempered monoidal unit} of $\Sph_G$.
\end{defn}

\begin{lem} \label{lem-temp-unit} For any left $\Sph_G$-module category $\mCC$ and object $\mCF\in \mCC$, we have 
	\begin{itemize}
		\item $\mCF$ is tempered iff there is an isomorphism $\mathbf{1}_{\Sph_G}^{\temp}\convolve \mCF \simeq \mCF$;
		\item $\mCF$ is anti-tempered iff  $\mathbf{1}_{\Sph_G}^{\temp}\convolve \mCF \simeq 0$.
	\end{itemize}
\end{lem}

\proof Follows from the definitions.

\qed[Lemma \ref{lem-temp-unit}]

\begin{exam} \label{exam-T-temper}
When $G=T$ is a torus, we have $\Sph_T^{\temp}=\Sph_T$. Hence $\mCC^{T\mathrm{-}\temp}=\mCC$.
\end{exam}

\begin{exam} \label{exam-temper-object-BunG}
For any fixed closed point $x\in X$ and fixed identification $\mCO_x\simeq \mCO$, consider the $\Sph_G$-action on $\DMod(\Bun_G)$ given by Hecke modifications at $\mCO_x\simeq \mCO$. By \cite{faergeman2021arinkin}, the resulting subcategory $\DMod(\Bun_G)^{G\mathrm{-}\temp}$ does \emph{not} depend on the choice of $x$ and $\mCO_x\simeq \mCO$.
\end{exam}

We also need the following technical lemma:

\begin{lem} \label{lem-cpt-gen-temp-general}
For any left $\Sph_G$-module category $\mCC$, if $\mCC$ is compactly generated, then so is $\mCC^{G\mathrm{-}\temp}$, and the functor $\mCC^{G\mathrm{-}\temp}\to \mCC$ preserves compact objects.
\end{lem}

\proof Note that the second claim follows from the first one because the functor $\mCC^{G\mathrm{-}\temp}\to \mCC$ has a continuous right adjoint $\temp_G$.

By \cite[Corollary 1.9.4]{drinfeld2015compact}, we only need to show the acting functors $\Sph_G\ot \mCC \to \mCC$ and $\Sph_G^{\temp}\ot \Sph_G \to \Sph_G^{\temp}$ have continuous right adjoints. This follows from Lemma \ref{lem-sphG-quasi-rigid} and Lemma \ref{lem-action-functor-quasi-rigid}.

\qed[Lemma \ref{lem-cpt-gen-temp-general}]

\subsection{The Eisenstein series and the constant term functors}
\label{ssec-EisCT}
In this subsection, we recall the Eisenstein series and the constant term functors. Details and proof can be found in \cite{braverman2002geometric}, \cite{drinfeld2016geometric}.

Consider the following diagram
\[ \Bun_M\xleftarrow{q} \Bun_P \xrightarrow{p} \Bun_G.\]
The map $q$ is quasi-compact, smooth and universally homologically contractible. The latter property means that, for any base-change $q'$ of $q$, the functor $(q')^!$ is fully faithful. The map $p$ is schematic and quasi-compact when restricted to each connected component of the source.

\begin{defn} Define the following four functors:
\[ 
\begin{aligned}
\Eis_{P,!}=p_!\circ q^*,\, \Eis_{P,*}= p_*\circ q^! : \DMod(\Bun_M) \to \DMod(\Bun_G);\\
\CT_{P,!}=q_!\circ p^*,\, \CT_{P,*}= q_*\circ p^! : \DMod(\Bun_G) \to \DMod(\Bun_M).
\end{aligned}
\]
\end{defn}

\begin{rem} Even though $p_!, p^*$ and $q_!, q^*$ are only partially defined functors, the compositions $p_!\circ q^*$ and $q_!\circ p^*$ are well-defined.
\end{rem}

\begin{thm}[{2nd adjointness for $\Bun_G$, \cite[Theorem 1.2.3]{drinfeld2016geometric}}] \label{thm-2nd-adjointness-DG}
 For a pair of opposite parabolic subgroups $(P,P^-)$ of $G$, there is a canonical equivalence $\CT_{P^-,!}\simeq \CT_{P,*}$.
\end{thm}

\subsection{The parabolic category \texorpdfstring{$\mCI(G,P)$}{I(G,P)}}
\label{ssec-IGP}
In this subsection, we review the definition of the parabolic category $\mCI(G,P)$.

Consider the prestack $\Bun_G^{P\hgen}:=\mathbf{Maps}_{\mathrm{gen}}(X, \mBB G\gets \mBB P)$ classifying $G$-torsors on $X$ equipped with generic reductions to $P$ (see \cite{barlev2012d}, and also \cite[\S 0.2]{chen2020deligne}). We have a diagram
\[ \Bun_M \xleftarrow{q} \Bun_P \xrightarrow{\iota_P} \Bun_G^{P\hgen}.  \]

\begin{rem} \label{rem-BunG-Pgen-to-BunG-pseudo-proper}
The prestack $\Bun_G^{P\hgen}$ is \emph{not} an algebraic stack because it does not have an atlas. Nevertheless, it is a reasonable geometric object. Namely, by \cite[Remark 4.1.9]{barlev2012d}, $\Bun_G^{P\hgen}$ is isomorphic to the fppf sheafification of a simplicial colimit $\colim_{\mathbf{\Delta}^\op} Y_\bullet$ such that $Y_n$ are algebraic stacks and each connected component of $Y_n$ is proper over $\Bun_G$. In fact, one can choose $Y_0:=\widetilde{\Bun}_P$ or $\overline{\Bun}_P$ to be the Drinfeld compactifications (see \cite{braverman2002geometric}) and $Y_\bullet$ be the Cech nerve of the map $Y_0\to \Bun_G^{P\hgen}$. It follows that we have
\begin{equation} \label{eqn-Dmod-bunG-P-gen}
 \DMod(\Bun_G^{P\hgen})\simeq \lim_{[n]\in \mathbf{\Delta}, !\hpull} \DMod(Y_n)\simeq \colim_{[n]\in \mathbf{\Delta}^\op, !\hpush} \DMod(Y_n). 
 \end{equation}
\end{rem}

The following result is implicit in \cite[Appendix B]{chen2020deligne}:

\begin{lem} \label{lem-stratification-BunG-Pgen}
For any affine finite type test scheme $S\to \Bun_G^{P\hgen}$, the base-change $\Bun_P\mt_{ \Bun_G^{P\hgen}} S$ is a finite disjoint union of a locally closed subschemes of $S$, and the map $\Bun_P\mt_{ \Bun_G^{P\hgen}} S\to S$ is both injective and surjective. In particular, $\Bun_P\to  \Bun_G^{P\hgen}$ is qcqs schematic.
\end{lem}

\proof By Remark \ref{rem-BunG-Pgen-to-BunG-pseudo-proper} and fppf descent, we can assume that $S\to \Bun_G^{P\hgen}$ factors through $\widetilde{\Bun}_P$. By \cite[\S B.1.3]{chen2020deligne}, there is a Cartesian square
\begin{equation} \label{eqn-stratum-BunG-Pgen}
\begin{tikzcd}
   _\mathrm{dfstr}\widetilde{\Bun}_P
    \arrow[r] \arrow[d]
    \arrow[dr, phantom, "\lrcorner", very near start]
    &\Bun_P \arrow[d,"\iota_P"']  \\
      \widetilde{\Bun}_P \arrow[r,"\widetilde{\iota}_P"']
    & \Bun_G^{P\hgen},
  \end{tikzcd}
\end{equation}
where $_\mathrm{dfstr}\widetilde{\Bun}_P$ is the disjoint union of the defect strata of $\widetilde{\Bun}_P$. This implies the claim of the lemma.

\qed[Lemma \ref{lem-stratification-BunG-Pgen}]

The above lemma suggests $\Bun_G^{P\hgen}$ should be viewed as glued from the connected components of $\Bun_P$. Hence $\DMod(\Bun_G^{P\hgen})$ should be viewed as glued from the direct summands of $\DMod(\Bun_P)$. Recall the functor $q^*:\DMod(\Bun_M)\to \DMod(\Bun_P)$ is fully faithful. Hence the above gluing induces a full subcategory of $\DMod(\Bun_G^{P\hgen})$ glued from the direct summands of $\DMod(\Bun_M)$. More precisely, we have the following definition:

\begin{defn}[\!\!{\cite[Definition 0.2.4]{chen2020deligne}}] 

Define $\mCI(G,P)\subset \DMod(\Bun_G^{P\hgen})$ to be the full subcategory such that we have the following Cartesian diagram
\[
\begin{tikzcd}
   \mCI(G,P)
    \arrow[r,"\subset"] \arrow[d,"\iota_M^!"]
    \arrow[dr, phantom, "\lrcorner", very near start]
    &\DMod(\Bun_G^{P\hgen}) \arrow[d,"\iota_P^!"']  \\
      \DMod(\Bun_M) \arrow[r,"q^*"',"\subset"]
    & \DMod(\Bun_P).
  \end{tikzcd}
\]
\end{defn}

\begin{rem} Since the map $q$ is smooth, in the above definition, we can use $q^!$ instead of $q^*$ and obtain the same full subcategory $\mCI(G,P)$.
\end{rem}

\begin{warn} The notation $\iota_M^!$ is only heuristic: it is \emph{not} induced by $!$-pullback along any (well-behaved) map $\iota_M$.

In fact, $\mCI(G,P)$ is equivalent to the category of D-modules on the pushout $\Bun_M \sqcup_{\Bun_P} \Bun_G^{P\hgen}$. However, this prestack is ill-behaved. Then $\iota_M$ can be viewed as the map 
\[\Bun_M \to \Bun_M \sqcup_{\Bun_P} \Bun_G^{P\hgen}.\]
\end{warn}

\begin{rem} \label{rem-adelic-picture}
Let us also provide the analogue picture in number theory. The stack $\Bun_M$ (resp. $\Bun_P$, $\Bun_G^{P\hgen}$) is an algebro-geoemtric incarnation of the double quotient $M(F)\backslash M(\mBA)/M(\mBO)$ (resp. $P(F)\backslash P(\mBA)/P(\mBO)$, $P(F)\backslash G(\mBA)/G(\mBO)$), where $F$ is the field of rational functions on $X$ and $\mBA$ (resp. $\mBO$) is the ring of adeles (resp. integral adeles) of $F$. Therefore the above prestack should be viewed as an algebro-geometric incarnation of $U(\mBA)M(F)\backslash G(\mBA)/G(\mBO)$, which is indeed isomorphic to $M(F)\backslash M(\mBA)/M(\mBO)$ as sets.
\end{rem}

\begin{lem}[\!\!{\cite[Corollary 0.2.8]{chen2020deligne})}]  \label{lem-generator-IGP}
We have adjoint functors:
\[  \iota_{M,!}: \DMod(\Bun_M) \adjoint \mCI(G,P): \iota_M^!,\]
where $\iota_M^!$ is induced by $\iota_P^!$ and is conservative, while $\iota_{M,!}$ is induced by $\iota_{P,!}\circ q^*$. The category $\mCI(G,P)$ is compactly generated by the image of compact objects of $\DMod(\Bun_M)$ under the functor $\iota_{M,!}$.
\end{lem}

\begin{defn} \label{defn-av-adelic}
Write $\oblv^{U(\mBA)}: \mCI(G,P)\to \DMod(\Bun_G^{P\hgen})$ for the forgetful functor. The above lemma implies it preserves compact objects, hence it admits a continuous right adjiont 
\[\Av_*^{U(\mBA)}: \DMod(\Bun_G^{P\hgen})\to \mCI(G,P).\]

In above functor, the notation $U(\mBA)$ is motivated by Remark \ref{rem-adelic-picture}.
\end{defn}

\begin{defn} Recall that the connected components $\Bun_{M,\lambda}$ of $\Bun_M$ (resp. $\Bun_P$) are labelled by $\lambda\in \Lambda_{G,P}$, which is the quotient of the coweight lattice $\Lambda_G$ by the $\mBZ$-span of simple coroots contained in $M$. Hence we obtain a $\Lambda_{G,P}$-grading on $\DMod(\Bun_{M})$. Let $\Bun_{P,\lambda}$ and $\widetilde{\Bun}_{P,\lambda}$ be the inverse images of $\Bun_{M,\lambda}$ along the maps $\Bun_P\to \widetilde{\Bun}_P \to \Bun_M$. 

For each $\lambda\in \Lambda_{G,P}$, we have adjoint functors
\[  (\iota_{M,\lambda})_!: \DMod(\Bun_{M,\lambda}) \adjoint \mCI(G,P): \iota_{M,\lambda}^!.\]
The partial order $\ge$ on $\Lambda_G$ induces a partial order $\ge$ on $\Lambda_{G,P}$. Let $\DMod(\Bun_{M})^{\ge \lambda}:= \bigoplus_{\mu \ge \lambda} \DMod(\Bun_{M,\mu}) $ and let
$\mCI(G,P)^{\ge \lambda}\subset \mCI(G,P)$ be the full subcategory generated under colimits by the image of $\iota_{M,!}(\DMod(\Bun_{M})^{\ge \lambda})$. Equivalently, it is the full subcategory of objects $\mCF$ such that $\iota_{M,\mu}^!(\mCF)\simeq 0$ unless $\mu\ge \lambda$.

We obtain a (decreasing) filtration on $\mCI(G,P)$ indexed by the poset $(\Lambda_{G,P},\ge)$.
\end{defn}

\begin{rem} The category $\mCI(G,P)^{\ge \lambda}\subset \mCI(G,P)$ should be viewed as the category of D-modules supported on a closed sub-prestack of $\Bun_M \sqcup_{\Bun_P} \Bun_G^{P\hgen}$.

\end{rem}

\begin{lem} \label{lem-filtration-on-IGP}
The functor $\iota_{M,!}$ and $\iota_{M}^!$ are compatible with the above filtrations on $\DMod(\Bun_{M})$ and $\mCI(G,P)$. Also, the endo-functor $\mathrm{coFib}(\Id\to \iota_{M}^! \circ \iota_{M,!})$ on $\DMod(\Bun_M)$ strictly increases the grading.
\end{lem}

\proof The functor $\iota_{M,!}$ is obviously compatible with the filtrations. To prove this for $\iota_{M}^!$, we only need to show that the endo-functor $\iota_{M}^!\circ \iota_{M,!}$ preserves $\DMod(\Bun_{M})^{\ge \lambda}$. In other words, we need to show $\iota_{M,\mu}^!\circ (\iota_{M,\lambda})_!\simeq 0$ unless $\mu\ge \lambda$. By definition, we only need to prove the corresponding claim for $\iota_{P,\mu}^!\circ (\iota_{P,\lambda})_!$, where $\iota_{P,\lambda}$ is the map $\Bun_{P,\lambda}\to \Bun_G^{P\hgen}$. Consider the diagram (\ref{eqn-stratum-BunG-Pgen}). It suffices to prove that $\iota_{P,\mu}^!\circ (\widetilde{\iota}_{P,\lambda})_!\simeq 0$ unless $\mu\ge \lambda$. 
By proper base-change, it suffices to show that the fiber product
\[ \widetilde{\Bun}_{P,\lambda}\mt_{\Bun_G^{P\hgen}} \Bun_{P,\mu} \]
is empty unless $\mu\ge \lambda$. In other words, we need to show the following: if $\widetilde{\Bun}_{P,\lambda}$ has a defect stratum over $\Bun_{P,\mu}$, then $\mu\ge \lambda$. But this was proved in \cite{braverman2002geometric}. This concludes the proof that $\iota_{M}^!$ is compatible with the filtrations.

For the second claim, we just need to show that $\Id \to \iota_{M,\lambda}^!\circ (\iota_{M,\lambda})_!$ is invertible. This follows from diagram (\ref{eqn-stratum-BunG-Pgen}) and the fact that $\Bun_{P,\lambda}$ is an open stratum of $\widetilde{\Bun}_{P,\lambda}$.

\qed[Lemma \ref{lem-filtration-on-IGP}]

\begin{cor} \label{cor-IGP-*-functors} The functor $\iota_M^!$ has a unique \emph{continuous} right inverse $\iota_{M,*}$, and the functor $\iota_{M,!}$ has a unique \emph{continuous} left inverse $\iota_M^*$. These functors are compatible with the filtrations on $\DMod(\Bun_{M})$ and $\mCI(G,P)$. Also, we have adjoint functors:
\[  \iota_{M}^*: \mCI(G,P) \adjoint \DMod(\Bun_{M,\lambda}): (\iota_{M})_*.\]
\end{cor}

\proof The proof below is standard, i.e., essentially the same as that for similar statements about D-modules on stratified spaces.

Consider the functor $\iota_{P,*}$. Note that $\Bun_P\mt_{\Bun_G^{P\hgen}} \Bun_P \simeq \Bun_P$. Hence we have $\iota_P^!\circ\iota_{P,*}\simeq \Id$. It follows that $\iota_{P,*}$ sends the image of $q^*$ into $\mCI(G,P)$. Let $ \iota_{M,*}:  \DMod(\Bun_M)\to \mCI(G,P)$ to the composition $\iota_{P,*}  \circ q^*$. By construction, it is a right inverse of $\iota_M^!$. The left adjoint $ \iota_{M,\lambda}^*$ of each $(\iota_{M,\lambda})_*$ is well-defined because the composition $\iota_{M,\lambda}^*\circ \iota_{M,!}$, which is the restriction functor $\DMod(\Bun_M) \to \DMod(\Bun_{M,\lambda})$, is well-defined.

By construction, the functor $\iota_{M,*}$ and $\iota_M^* :=\oplus_{\lambda} \iota_{M,\lambda}^*$ are compatible with filtrations. 

To prove the adjoint pair $(\iota_{M}^*,\iota_{M,*})$, we only need to show $\bigoplus (\iota_{M,\lambda})_*\to \prod  (\iota_{M,\lambda})_*$ is an equivalence. The conservative functor $\iota_{M}^!$ obviously commutes with colimits, and commutes with limits because it has a left adjoint $\iota_{M,!}$. Hence we only need to show $\bigoplus \iota_{M}^!\circ (\iota_{M,\lambda})_*\to \prod \iota_{M}^!\circ  (\iota_{M,\lambda})_*$ is an equivalence but this is obvious.

It remains to prove the uniqueness. By the Barr--Beck--Lurie theorem, $\mCI(G,P)\simeq (\iota_{M}^! \circ \iota_{M,!})\mathrm{-mod}( \DMod(\Bun_M) )$. By Lemma \ref{lem-filtration-on-IGP}, for each object in  $\DMod(\Bun_{M,\lambda})\subset \DMod(\Bun_{M})$, there is a unique $(\iota_{M}^! \circ \iota_{M,!})$-module structure on it. Unwinding the definitions, this means there is a unique functor $F_\lambda:\DMod(\Bun_{M,\lambda}) \to \mCI(G,P)$ such that $\iota_{M,\mu}^!\circ F_\lambda \simeq 0$ for $\mu\neq \lambda$ and $\iota_{M,\lambda}^!\circ F_\lambda\simeq \Id$. Therefore any \emph{continuous} right inverse of $ \iota_M^!$ must be $\oplus_\lambda F_\lambda$. This proves the uniqueness of $F$.

On the other hand, for any continuous left inverse $E$ of $\iota_{M,!}$, its restriction $E_\lambda:\mCI(G,P)\to \DMod(\Bun_{M,\lambda})$ is continuous and satisfies $E_\lambda\circ (\iota_{M,\mu})_! \simeq 0$ for $\mu\neq \lambda$ and $E_\lambda\circ (\iota_{M,\lambda})_! \simeq \Id$. Passing to right adjoints, we see $F_\lambda$ is the right adjoint of $E_\lambda$. Hence $E$ is also unique. 

\qed[Corollary \ref{cor-IGP-*-functors}]

\begin{defn} \label{defn-iota-*-gen}
Let $\mCI(G,P)^{\iota_*\hgen}\subset \mCI(G,P)$ be the full subcategory generated under colimits by the image of $\iota_{M,*}: \DMod(\Bun_{M}) \to \mCI(G,P)$. 
\end{defn}

\begin{cor} \label{cor-iota-*-gen-description} For $\lambda\in\Lambda_{G,P}$, suppose $\mCF\in \mCI(G,P)$ is such that $\iota_{M,\mu}^!(\mCF)\simeq 0$ unless $\mu \le \lambda$, then $\mCF\in \mCI(G,P)^{\iota_*\hgen}$.
\end{cor}

\proof Follows from Lemma \ref{lem-filtration-on-IGP} and Corollary \ref{cor-IGP-*-functors} by a standard argument.

\qed[Corollary \ref{cor-iota-*-gen-description}]

The following result seems to be well-known but we can not find a proof in the literature. We provide a proof in \S \ref{ssec-IGP-proof-1}.

\begin{prop} \label{prop-sph-acts-on-IGP}
For any parabolic subgroup $P$ of $G$ and fixed closed point $x\in X$, consider the $\Sph_{G,x}$-action on $\DMod(\Bun_G^{P\hgen})$ given by Hecke modifications at $x$. Then this action preserves $\mCI(G,P)\subset \DMod(\Bun_G^{P\hgen})$.
\end{prop}

Now we study the functorial properties of the parabolic categories $\mCI(G,P)$.

\begin{defn} \label{defn-Eis-CT-enh}
For any parabolic subgroups $Q\subset P$ of $G$, define adjoint functors
\[  \Eis_{Q\to P}^\enh: \mCI(G,Q)\adjoint \mCI(G,P): \CT_{P\gets Q}^\enh\]
as follows. Consider the map 
\[ p_{Q\to P}^\enh: \Bun_G^{Q\hgen} \to \Bun_G^{P\hgen}\]
sending a generic $Q$-reduction to its induced $P$-reduction. The functor $(p_{Q\to P}^\enh)_!$ is well defined (see Remark \ref{rem-!-push-gen-exist} below) and sends $\mCI(G,Q)$ into $\mCI(G,P)$ (see \cite[Proposition 0.2.10]{chen2020deligne}). The left adjoint $\Eis_{Q\to P}^\enh$ is defined as its restriction. 

The functor $\CT_{P\gets Q}^\enh$ is defined as its continuous right adjoint. By definition, we have 
\[\CT_{P\gets Q}^\enh\simeq \Av_*^{U_Q(\mBA)}\circ (p_{Q\to P}^\enh)^!.\]
\end{defn}

\begin{warn} The functor $(p_{Q\to P}^\enh)^!$ does not sends $\mCI(G,P)$ to $\mCI(G,Q)$. 
\end{warn}

\begin{rem} \label{rem-!-push-gen-exist}
Recall we have a map $\overline{\Bun}_Q\to \overline{\Bun}_P$ compatible with the map $p_{Q\to P}^\enh: \Bun_G^{Q\hgen} \to \Bun_G^{P\hgen}$. Hence by the equivalence (\ref{eqn-Dmod-bunG-P-gen}) (with the choice $Y_0:= \overline{\Bun}_P$), well-definedness of $(p_{Q\to P}^\enh)_!$ follows from well-definedness of proper $!$-pushforward functors.
\end{rem}

\begin{rem} \label{rem-enh-eisct-vs-eisct}
By definition, we have 
\[\Eis_{P,!}\simeq \Eis_{P\to G}^\enh\circ \iota_{M,!},\;\; \CT_{P,*}\simeq \iota_{M}^!\circ \CT_{G\gets P}^\enh. \]

On the other hand, by Remark \ref{rem-BunG-Pgen-to-BunG-pseudo-proper}, the functor $(p^\enh_{P\to G})_!$ can be calculated as the simplicial colimit of 
\[ \DMod(\Bun_G^{P\hgen}) \xrightarrow{!\hpull} \DMod(Y_n) \xrightarrow{!\hpush} \DMod(\Bun_G), \]
where $Y_n$ is the Cech nerve of the map $\widetilde{\Bun}_P\to \Bun_G$. Note that the above $!$-pushforward is also $*$-pushforward. It follows that the composition 
\[ \DMod(\Bun_P)\xrightarrow{*\hpush} \DMod(\Bun_G^{P\hgen})\xrightarrow{!\hpush} \DMod(\Bun_G)  \]
is equivalent to the $*$-pushforward functor. In particular,
\[\Eis_{P,*}\simeq \Eis_{P\to G}^\enh\circ \iota_{M,*}\circ [2\mathrm{rel.dim.}(\Bun_P,\Bun_M)],\]
where $\mathrm{rel.dim.}(\Bun_P,\Bun_M)$ is the locally constant function on $\Bun_M$ whose values are the dimensions of the fibers of $\Bun_P\to \Bun_M$.
\end{rem}

The following result will be used repeatedly in this paper without explicit reference:

\begin{lem} \label{lem-EisCT-enh-Sph-linear}
For any inclusion $Q \subset P$ of parabolic subgroups of $G$, consider the $\Sph_{G,x}$-actions on $\mCI(G,Q)$ and $\mCI(G,P)$. The adjoint functors 
\[ \Eis: \mCI(G,Q)\adjoint \mCI(G,P):\CT \]
are canonically $\Sph_{G,x}$-linear.
\end{lem}

\proof The $!$-pullback functor along $\Bun_G^{Q\hgen}\to \Bun_G^{P\hgen}$ is canonically $\Sph_{G,x}$-linear by the base-change isomorphisms. By definition, $\CT$ is the composition of this functor with the localization functor $\DMod( \Bun_G^{Q\hgen})\to \mCI(G,Q)$, i.e., with the right adjoint to the embedding $\mCI(G,Q)\to \DMod( \Bun_G^{Q\hgen})$. By definition, this embedding is canonically $\Sph_{G,x}$-linear, hence so is its right adjoint (Appendix \ref{sect-lax-linear}). It follows that $\CT$ is canonically $\Sph_{G,x}$-linear. Using Appendix \ref{sect-lax-linear} again, its left adjoint $\Eis$ is also canonically $\Sph_{G,x}$-linear.

\qed[Lemma \ref{lem-EisCT-enh-Sph-linear}]

\subsection{Proof of Proposition \ref{prop-sph-acts-on-IGP}}
\label{ssec-IGP-proof-1}

In this subsection, we prove Proposition \ref{prop-sph-acts-on-IGP}.

For $\mCM\in \mCI(G,P)$ and $\mCF\in \Sph_{G,x}$, we need to show that $\mCF \convolve_{\Sph_{G,x}}\mCM \in \mCI(G,P)$. Consider the prestack $\Hecke_{G,x}^{P\hgen}$ classifying quadruples $(\mCP_G^l,\mCP_G^r,\theta,\mCP_P^{\mathrm{gen}})$, where:
\begin{itemize}
\item $\mCP_G^l$ and $\mCP_G^r$ are $G$-torsors on $X$;
\item $\theta:\mCP_G^l|_{X\setminus x}\to \mCP_G^r|_{X\setminus x}$ is an isomorphism, which in particular allows to identify the generic $G$-torsors underlying $\mCP_G^l$ and $\mCP_G^r$;
\item $\mCP_P^{\mathrm{gen}}$ is a $P$-reduction of this generic $G$-torsor.
\end{itemize}
Let $h_l,h_r: \Gr_{G,x}\widetilde{\mt} \Bun_G^{P\hgen}\to \Bun_G^{P\hgen}$ be the maps sending the above data to $(\mCP_G^l, \mCP_P^{\mathrm{gen}})$ and $(\mCP_G^r, \mCP_P^{\mathrm{gen}})$ respectively. Note that the map $h_r$ witnesses $\Hecke_{G,x}^{P\hgen}$ as a twisted product
\[ \Hecke_{G,x}^{P\hgen} \simeq \Gr_{G,x}\widetilde{\mt} \Bun_G^{P\hgen}. \]
Then we have
\[ \mCF \convolve_{\Sph_{G,x}}\mCM \simeq h_{l,*}(\mCF\widetilde{\boxtimes}\mCM) .\]
Write $\iota_P: \Bun_P \to  \Bun_G^{P\hgen}$. We only need to show that $\iota_P^!\circ h_{l,*}(\mCF\widetilde{\boxtimes}\mCM)$ is in the essential image of the pullback functor $\DMod(\Bun_{M})\to \DMod(\Bun_P)$, where $M$ is the Levi quotient group of $P$.

The fiber product 
\[\Bun_P\mt_{\Bun_G^{P\hgen},h_l} \Hecke_{G,x}^{P\hgen} \] classifies $(\mCP_P^l,\mCP_G^r,\delta)$, where $\mCP_P^l$ (resp. $\mCP_G^r$) is a $P$-torsor ($G$-torsor) on $X$, and $\delta: G\mt^P \mCP_P^l|_{X\setminus x}\to \mCP_G^r|_{X\setminus x}$ is an isomorphism. In other words, we have
\[\Bun_P\mt_{\Bun_G^{P\hgen},h_l} \Hecke_{G,x}^{P\hgen} \simeq \Bun_P\mt_{\Bun_G,h_l} \Hecke_{G,x} \] 
Hence we have the following Cartesian square
\[
  \begin{tikzcd}
     \Bun_P\mt_{\Bun_G,h_l} \Hecke_{G,x}
    \arrow[r,"v"] \arrow[d,"\mathrm{pr}_1"]
    \arrow[dr, phantom, "\lrcorner", very near start]
    &\Hecke_{G,x}^{P\hgen} \arrow[d,"h_l"]  \\
      \Bun_P \arrow[r,"\iota_P"]
    & \Bun_G^{P\hgen}.
  \end{tikzcd}
\]
Using the base-change isomorphism, we obtain
\[  \iota_P^!\circ h_{l,*}(\mCF\widetilde{\boxtimes}\mCM) \simeq \mathrm{pr}_{1,*}\circ v^!( \mCF\widetilde{\boxtimes}\mCM ). \]
Note that we also have
\[ \Bun_P\mt_{\Bun_G,h_l} \Hecke_{G,x}\simeq \Bun_P \widetilde{\mt} \Gr_{G,x},\; \Hecke_{P,x}\simeq \Bun_P \widetilde{\mt} \Gr_{P,x} .\]
It is well-known that $\Gr_{P,x}\to \Gr_{G,x}$ is a stratification (see e.g. \cite[Appendix C.3]{chen2020nearby}), hence so is the map 
\[u:\Hecke_{P,x} \to \Bun_P\mt_{\Bun_G,h_l} \Hecke_{G,x}.\]
Hence it remains to show that $ \mathrm{pr}_{1,*}\circ u_*\circ u^!  \circ v^!( \mCF\widetilde{\boxtimes}\mCM )$ is contained in the essential image of the pullback functor $\DMod(\Bun_{M})\to \DMod(\Bun_P)$. Note that we can write $\Hecke_{P,x}$ as
\[\Hecke_{P,x} \simeq \Gr_{P,x} \widetilde{\mt}\Bun_P.\]
Thanks to this, we obtain
\[\mathrm{pr}_{1,*}\circ u_*\circ u^!  \circ v^!( \mCF\widetilde{\boxtimes}\mCM ) \simeq h_{l,*} (w^!(\mCF) \widetilde{\boxtimes} \iota_P^!(\mCM)),
\]
where $w^!(\mCF) \in \Sph_{P,x}$ is the pullback of $\mCF$ along $w:P(\mCO_x)\backslash P(\mCK_x)/P(\mCO_x)\to G(\mCO_x)\backslash G(\mCK_x)/G(\mCO_x)$ and $h_l: \Hecke_{P,x}\to \Bun_P$ is the left projection map. In other words, the RHS is just $w^!(\mCF)\convolve_{\Sph_{P,x}} \iota_P^!(\mCM) $. By assumption, $\iota_P^!(\mCM) \simeq q^!(\mCN)$ for some $\mCN\in \DMod(\Bun_{M})$, where $q:\Bun_P\to \Bun_{M}$ is the projection. Now Lemma \ref{lem-sphP-preserve-BunM} below implies $w^!(\mCF)\convolve_{\Sph_{P,x}} \iota_P^!(\mCM)$ is contained in the essential image of $q^!$ as desired.

\qed[Proposition \ref{prop-sph-acts-on-IGP}]

The following lemma was used above:
\begin{lem} \label{lem-sphP-preserve-BunM}
For any parabolic subgroup $P$ of $G$ and its Levi quotient group $M$, the $\Sph_{P,x}$-action on $\DMod(\Bun_P)$ preserves the essential image of the functor $q^!:\DMod(\Bun_{M})\to \DMod(\Bun_P)$.
\end{lem}

\proof Consider the maps $\mathrm{pr}_1,\mathrm{pr}_2:\Bun_P\mt_{\Bun_{M}} \Bun_P \to \Bun_P$. Since $q$ is universally homologically contractible, the functors $\mathrm{pr}_1^!,\mathrm{pr}_2^!$ are fully faithful. By smooth descent for D-modules, an object $\mCK\in \DMod(\Bun_P)$ is in the essential image of $q^!$ iff there exists an equivalence $\mathrm{pr}_1^!(\mCK) \simeq \mathrm{pr}_2^!(\mCK)$.

Let $y\neq x$ be another closed point of $X$. Consider the fiber product
\[H:=\Hecke_{P,y}\mt_{\Hecke_{M,y}} \Bun_M\]
classifying Hecke modifications of $P$-torsors at $y$ whose induced modifications of $M$-torsors are trivial. By Lemma \ref{lem-UHC-unipotent-hecke} below, the natural map
\[ \theta: H \to \Bun_P\mt_{\Bun_{M}} \Bun_P \]
is universally homologically contractible.

Combining the above two paragraphs, we see an object $\mCK\in \DMod(\Bun_P)$ is in the essential image of $q^!$ iff there exists an equivalence $(\mathrm{pr}_1\circ \theta)^!(\mCK) \simeq (\mathrm{pr}_2\circ \theta)^!(\mCK)$.

Now let $\mCF\in \Sph_{P,x}$ and $\mCM = q^!(\mCN)\in \DMod(\Bun_P)$. Using the assumption $y\neq x$, it is easy to see 
\[ H\mt_{\mathrm{pr}_1\circ \theta,\Bun_P, h_l} \Hecke_{P,x} \simeq  H\mt_{\mathrm{pr}_2\circ \theta,\Bun_P, h_l} \Hecke_{P,x} \simeq \Hecke_{P,x\cup y}\mt_{\Hecke_{M,y}} \Bun_M.\]
This space classifies Hecke modifications of $P$-torsors at $x\cup y$ whose induced modifications of $M$-torsors at $y$ are trivial. Thus, for $i=1$ or $2$, we have maps
\[v_i: H':=\Hecke_{P,x\cup y}\mt_{\Hecke_{M,y}} \Bun_M\to \Hecke_{P,x}  \]
and Cartesian diagrams
\[
  \begin{tikzcd}
    H'
    \arrow[r,"v_i"] \arrow[d,"\mathrm{pr}_1"]
    \arrow[dr, phantom, "\lrcorner", very near start]
    &\Hecke_{P,x} \arrow[d,"h_l"]  \\
      H \arrow[r,"\mathrm{pr}_i\circ \theta"]
    & \Bun_P.
  \end{tikzcd}
\]
Using the base-change isomorphisms, it remains to show that
\[v_1^!( \mCF \widetilde{\boxtimes} q^!{\mCN} ) \simeq v_2^!( \mCF \widetilde{\boxtimes} q^!{\mCN} ). \]
But this follows from the fact that the maps
\[  H' \xrightarrow{v_i} \Hecke_{P,x} \to P(\mCO_x)\backslash P(\mCK_x)/P(\mCO_x),\; H' \xrightarrow{v_i} \Hecke_{P,x} \xrightarrow{h_r} \Bun_P \xrightarrow{q} \Bun_M \]
are independent of $i$.

\qed[Lemma \ref{lem-sphP-preserve-BunM}]

The following lemma was used above:
\begin{lem} \label{lem-UHC-unipotent-hecke}
For any closed point $y\in X$, the natural map 
\[ \theta:\Hecke_{P,y}\mt_{\Hecke_{M,y}} \Bun_M\to  \Bun_P\mt_{\Bun_{M}} \Bun_P \]
is universally homologically contractible.
\end{lem}

\proof For any finite type affine test scheme $s: S\to \Bun_P\mt_{\Bun_{M}} \Bun_P$, we need to show that base-change $\theta': S'\to S$ of $\theta$ along $s$ is homologically contractible, i.e., $(\theta')^!$ is fully faithful. By definition, $s$ corresponds to $(\mCF_P^l,\mCF_P^r,\alpha)$ where $\mCF_P^l$ and $\mCF_P^r$ are $P$-torsors on $X\times S$ and $\alpha$ is an isomorphism between their induced $M$-torsors. Then, for any finite type affine test $S$-scheme $T$, the groupoid $\Maps_{S}(T,S')$, which is actually a set, classifies isomorphisms 
\[\beta: \mCF_P^l|_{(X\setminus y) \times T} \to \mCF_P^r|_{(X\setminus y)\times T} \]
 between $P$-torsors on $(X\setminus x)\times T$ such that the induced isomorphism
 \[ M\mt^P \beta: M\mt^P \mCF_P^l|_{(X\setminus y) \times T} \to  M\mt^P \mCF_P^r|_{(X\setminus y)\times T} \]
is the base-change of $\alpha$ along $(X\setminus y) \times T\to X\mt S$. 

We first show that $\Maps_{S}(T,S')$ is non-empty. To this end, we can assume $T=S$. Since $(X\setminus y) \times S$ is affine, any $P$-torsor on it has an $M$-reduction. Choose $M$-reductions $\mCF_M^l$ and $\mCF_M^r$ of $\mCF_P^l|_{(X\setminus y) \times S}$ and $\mCF_P^r|_{(X\setminus y)\times S}$. Then we have
\[ M\mt^P \mCF_P^l|_{(X\setminus y) \times S} \simeq M\mt^P P\mt^M  \mCF_M^l \simeq \mCF_M^l \]
and similar isomorphisms for $\mCF_M^r$. Hence the base-change of $\alpha$ induces an isomorphism between $\mCF_M^l$ and $\mCF_M^r$, which in turn yields an element in $\Maps_{S}(S,S')$.

Let $H_T$ be the (abstract) group of those automorphisms of the $P$-torsor $\mCF_P^l|_{(X\setminus y) \times T} $ such that the induced automorphism of $M\mt^P \mCF_P^l|_{(X\setminus y) \times T}$ is the identity map. It follows that $\Maps_{S}(T,S')$ is a non-empty $H_T$-torsor. On the other hand, using the $M$-reductions $\mCF_M^l$ of $\mCF_P^l|_{(X\setminus y) \times S} $, we see that $H_T$ can be identified with the group of $M$-equivariant maps $\mCF_M^l|_{ (X\setminus y) \times T } \to U$. Choose a $\mBG_M$-action on $G$ that contracts $U$ into the unit element and fixes $M$. Using this $\mBG_m$-action, we see that $T\mapsto H_T$ is represented by a contractible group prestack over $S$. It follows that $T\mapsto \Maps_{S}(T,S')$ is represented by a contractible prestack over $S$. By definition, this prestack is $S'$.

\qed[Lemma \ref{lem-UHC-unipotent-hecke}]

\section{Tempered objects in the parabolic category \texorpdfstring{$\mCI(G,P)$}{I(G,P)}} \label{sect-CT}

Throughout this section, we fix a parabolic subgroup $P$ of $G$. Let $U$ be its unipotent radical and $M=P/U$ be its Levi quotient group. We fix a splitting $M\to P$ and view $M$ also as a Levi subgroup. 

We fix a closed point $x\in X$ and an identification $\mCO_x \simeq \mCO$, which provide an equivalence $\Sph_{G}\simeq \Sph_{G,x}$. We use this equivalence to define $G$-tempered objects in $\DMod(\Bun_G)$ and $\mCI(G,P)$ (see Example \ref{exam-temper-object-BunG} and Proposition \ref{prop-sph-acts-on-IGP}).

The goal of this section is to prove the following theorems:

\begin{thm} \label{thm-EisCT-temperedness} We have:
\begin{itemize}
	\item[(1)] The functor $\CT_{P,*}:  \DMod(\Bun_G) \to \DMod(\Bun_{M})$ sends $G$-tempered objects to $M$-tempered objects.

	\item[(2)] The functor $\Eis_{P,!}: \DMod(\Bun_M) \to \DMod(\Bun_{G})$ sends $M$-anti-tempered objects to $G$-anti-tempered objects.

	\item[(3)] The functor $\Eis_{P,*}: \DMod(\Bun_M) \to \DMod(\Bun_{G})$ sends $M$-tempered objects to $G$-tempered objects.
\end{itemize}
\end{thm}

\begin{thm} \label{thm-iota!*-temperedness} \label{thm-iota-temperedness} We have:
\begin{itemize}
	\item[(1)] The functor $\iota_M^!:  \mCI(G,P)\to \DMod(\Bun_M)$ sends $G$-tempered objects to $M$-tempered objects.

	\item[(2)] The functor $\iota_{M,!}: \DMod(\Bun_M) \to  \mCI(G,P)$ sends $M$-anti-tempered objects to $G$-anti-tempered objects.

	\item[(3)] The functor $\iota_{M,*}: \DMod(\Bun_M) \to  \mCI(G,P)$ sends $M$-tempered objects to $G$-tempered objects.

	\item[(4)] The functor $\iota_{M}^*:  \mCI(G,P) \to  \DMod(\Bun_M)$ sends $G$-anti-tempered objects to $M$-anti-tempered objects.
\end{itemize}
\end{thm}

Theorem \ref{thm-EisCT-temperedness} can be deduced from Theorem \ref{thm-iota!*-temperedness} as follows.

\proof[Proof of Theorem \ref{thm-EisCT-temperedness}.] Recall the functors 
\[ \Eis_{P\to G}^\enh:\mCI(G,P)\adjoint \DMod(\Bun_G): \CT_{G\gets P}^\enh \]
are $\Sph_{G,x}$-linear. Hence by Remark \ref{rem-enh-eisct-vs-eisct}, we see Theorem \ref{thm-EisCT-temperedness}(1)-(3) follow respectively from Theorem \ref{thm-iota!*-temperedness}(1)-(3).

\qed[Theorem \ref{thm-EisCT-temperedness}]

Our strategy for Theorem \ref{thm-iota!*-temperedness} is to deduce it from its local analogue. This is done as follows.

In \S \ref{ssec-SI}, we recall the local analogue of $\mCI(G,P)$, known as the \emph{semi-infinite category} 
\[\SI_P := \DMod(\Gr_G)^{U(\mCK)M(\mCO)}.\]

In \S \ref{ssec-SI-temper}-\ref{ssec-HL-temper-2}, we state and prove the local analogue of Theorem \ref{thm-iota!*-temperedness} (see Theorem \ref{thm-iota!*-temperedness-local}).

In \S \ref{ssec-SI-temper}, we reduce this local result to proving that a certain hyperbolic localization functor $\HL_P^{!}$ (resp. $\HL_P^*$): $\Sph_G\to \Sph_M$ preserves temperedness (resp. detects anti-temperedness), see Proposition \ref{prop-HL!-temper} and Proposition \ref{prop-HL*-temper}. The main ingredient for this reduction is the duality between $\SI_P$ and $\SI_{P^-}$ worked out in \cite{chen2020nearby}, \cite{chen2021thesis}.

In \S \ref{ssec-HL-temper-1}, we prove the claim about $\HL_P^{!}$. The main ingredients in its proof include: (i) V. Lafforgue's equivalence $\Sph_G \simeq \DMod(\Bun_G(\mBP^1))$, see \cite{lafforgue2009quelques} and also \cite[\S 3.1]{beraldo2021tempered}; (ii) the compatibility between $\HL_P^{!}$ and $\CT_{P,*}$ via Lafforgue's equivalence (Proposition \ref{prop-local-global-CT-P1}); (ii) the explicit description of tempered objects in $\DMod(\Bun_G(\mBP^1))$ in \cite[\S 3.5]{beraldo2021geometric}.

In \S \ref{ssec-HL-temper-2}, we prove the claim about $\HL_P^{*}$. It can be easily reduced to the case $P=B$, and in that case it follows from the well-known fact (see \cite{mirkovic2007geometric}) that the hyperbolic localization functor $\HL_B^{*}$ sends the heart of $\Sph_G$ to the heart of $\Sph_T$ up to a cohomological shift.

In \S \ref{ssec-SI-to-IGP-temper}, we deduce Theorem \ref{thm-iota!*-temperedness} from its local analogue. The main ingredient is a factorization of the functor $\iota_{M,!}$ as
\[\DMod(\Bun_M)\to \SI_{P,x}\ot_{\Sph_{M,x}} \DMod(\Bun_M) \to \mCI(G,P) \]
which relates the local and global picture.

In \S \ref{ssec-CTtemmperlocalspectal}, we provide a (conjectural) explanation of the results about $\SI_P$ on the spectral side of the Langlands duality.

\subsection{Recollection: the semi-infinite category \texorpdfstring{$\SI_P$}{ }}
\label{ssec-SI}

In this section, we recall the semi-infinite category
\[
\SI_{P}:= \DMod(\Gr_G)^{ U(\mCK)M(\mCO) },
\]
which is the $U(\mCK)M(\mCO)$-invariant category of $\DMod(\Gr_G)$ (See \cite[\S 2]{raskin2016chiral} or \cite[\S 1.4]{chen2020nearby}). It carries commuting actions of the spherical categories $\Sph_G$ and $\Sph_M$. Namely, if we realize $\Gr_G$ as $G(\mCO)\backslash G(\mCK)$, then the $\Sph_G$-action is induced by convolution from the left, while the $\Sph_M$-action is induced by the right $M(\mCK)$-action on $G(\mCO)\backslash G(\mCK)$.

As mentioned at the beginning of this section, $\SI_P$ should be viewed as a local analogue of $\mCI(G,P)$. Let us describe the local analogue of the adjoint functors
\begin{eqnarray*}  
\iota_{M,!}:& \DMod(\Bun_M) \adjoint \mCI(G,P) &: \iota_M^!, \\
 \iota_{M}^*:& \mCI(G,P) \adjoint \DMod(\Bun_M) &: \iota_{M,*}.
\end{eqnarray*}

\begin{constr}
We have adjoint functors (See \cite[Lemma 2.3.4]{chen2020nearby})
\begin{eqnarray*}
\iota_{P,!}: \DMod(\Gr_P)^{U(\mCK)} \adjoint \DMod(\Gr_G)^{U(\mCK)}: \iota_P^!,\\
\iota_{P}^*: \DMod(\Gr_G)^{U(\mCK)} \adjoint \DMod(\Gr_P)^{U(\mCK)}  : \iota_{P,*}
\end{eqnarray*}
induced by pushforwards and pullbacks along\footnote{In \S \ref{ssec-IGP}, we wrote the map $\Bun_P\to \Bun_G^{P\hgen}$ using the same notation $\iota_P$. We will articulate the source and target of the maps if there is danger of ambiguity. In any case, the map $\Gr_P\to \Gr_G$ should be viewed as a local analogue of $\Bun_P\to \Bun_G^{P\hgen}$.} $\iota_P:\Gr_P\to \Gr_G$. 

The functors $\iota_P^!$ and $\iota_{P,*}$ are obviously $M(\mCK)$-linear; by \cite[Lemma D.4.4]{gaitsgory2020local}, so are their left adjoints. Hence, the above adjoint functors induce $\Sph_M$-linear adjoint functors on $M(\mCO)$-invariants, i.e., between $\SI_P$ and $\DMod(\Gr_P)^{U(\mCK)M(\mCO)}$.

On the other hand, the $!$-pullback functor along $\Gr_P\to \Gr_M$ induces equivalences (see \cite[Lemma 2.3.2]{chen2020nearby}):
\[\DMod(\Gr_M) \simeq \DMod(\Gr_P)^{U(\mCK)} ,\; \Sph_M \simeq \DMod(\Gr_P)^{U(\mCK)M(\mCO)}.  \]
Combining with the above paragraph, we obtain $\Sph_M$-linear adjoint functors:
\begin{eqnarray*}
\iota_{M,!}:\Sph_M \adjoint \SI_P: \iota_M^!,\;\;\; \iota_{M}^*: \SI_P \adjoint \Sph_M  : \iota_{M,*}
\end{eqnarray*}
\end{constr}

It is well-known that the connected components of $\Gr_P$ provide a stratification on $\Gr_G$ (see \cite[\S C.3]{chen2020nearby}). It follows that $\iota_M^!$ is conservative and that $\SI_P$ is compactly generated by $\iota_{M,!}(\mCF)$ for compact objects $\mCF\in \Sph_M$. Also, $\iota_{M,*}$ (resp. $\iota^*$) is the unique continuous right inverse of $\iota_M^!$ (resp. left inverse of $\iota_{M,!}$).

\begin{thm}[\!\!{\cite[Theorem A, Theorem C]{chen2021thesis}}] \label{thm-inv-inv-duality}
Let $P^-$ be any parabolic subgroup of $G$ opposite to $P$ and $U^-$ be its unipotent radical. Then the categories $\DMod(\Gr_G)^{U(\mCK)}$ and $\DMod(\Gr_G)^{U^-(\mCK)}$ are dual to each other, with the pairing functor given by
\[ \DMod(\Gr_G)^{U(\mCK)}\ot \DMod(\Gr_G)^{U^-(\mCK)} \xrightarrow{\oblv} \DMod(\Gr_G)\ot\DMod(\Gr_G) \to \mathrm{Vect}, \]
where the rightmost arrow is the Verdier self-pairing on $\Gr_G$. Moreover, via this duality and the Verdier self-duality on $\Gr_M$, the functors 
\begin{eqnarray*}
\DMod(\Gr_M) &\simeq \DMod(\Gr_P)^{U(\mCK)}&\xrightarrow{\iota_{P,!}} \DMod(\Gr_G)^{U(\mCK)},\\
\DMod(\Gr_M) &\simeq \DMod(\Gr_{P^-})^{U^-(\mCK)}&\xrightarrow{\iota_{P^-,!}} \DMod(\Gr_G)^{U^-(\mCK)}
\end{eqnarray*}
are conjugate to each other.
\end{thm}

One can easily check that the Verdier self-pairing $\DMod(\Gr_G)\ot\DMod(\Gr_G) \to \mathrm{Vect}$ is compatible with the $\Sph_G$-action on $\DMod(\Gr_G)$. It follows that the duality in Theorem \ref{thm-inv-inv-duality} is compatible with the $\Sph_G$-actions on $\DMod(\Gr_G)^{U(\mCK)}$ and $\DMod(\Gr_G)^{U^-(\mCK)}$. Similarly, the Verdier self-pairing on $\Gr_G$ is compatible with the $\mCM(\mCK)$-action. It follows we can take $\mCM(\mCO)$-invariants and obtain the following result.

\begin{cor}[C.f. {\cite[Corollary 1.4.5]{chen2020nearby}}] \label{cor-inv-inv-duality-L^+M}
Let $P^-$ be any parabolic subgroup of $G$ opposite to $P$. Then the categories $\SI_P$ and $\SI_{P^-}$ are dual to each other. This duality is compatible with the $\Sph_M$-actions and $\Sph_G$-actions. Moreover, via this duality and the Verdier self-duality of $\Sph_M$, the functors 
\begin{eqnarray*}
\iota_{M,!}: \Sph_M \to \SI_P,\;\; \iota_{M,!}^-: \Sph_M \to \SI_{P^-}
\end{eqnarray*}
are conjugate to each other.
\end{cor}

\begin{cor} \label{cor-inv-inv-duality-functors}
Via the above duality between $\SI_P$ and $\SI_{P^-}$, and the Verdier self-duality of $\Sph_M$, we have:
\begin{itemize}
	\item[(1)] The functors 
	\begin{eqnarray*}
		\iota_{M}^!: \SI_P \to \Sph_M ,\;\; \iota_{M,!}^-: \Sph_M \to \SI_{P^-}
	\end{eqnarray*}
	are dual to each other.
	\item[(2)] The functors 
	\begin{eqnarray*}
		\iota_{M}^*: \SI_P \to \Sph_M ,\;\; \iota_{M,*}^-: \Sph_M \to \SI_{P^-}
	\end{eqnarray*}
	are dual to each other. 
\end{itemize}
\end{cor}

\proof (1) follows from the formal fact that the conjugate functor is left adjoint to the dual functor (see e.g. \cite[\S 1.5]{gaitsgory2016functors}). (2) follows from the fact that $\iota_{M,*}$ (resp. $\iota^*$) is the unique right inverse of $\iota_M^!$ (resp. left inverse of $\iota_{M,!}$).

\qed[Corollary \ref{cor-inv-inv-duality-functors}]

The last topic of this subsection is to describe the local analogue of the functors 
\[\Eis^\enh_{P\to G}:\mCI(G,P) \adjoint  \DMod(\Bun_G) :\CT_{G\gets P}^\enh. \]

\begin{constr} We have adjoint functors
\[ \Av_!^{U(\mCK)/U(\mCO)}: \Sph_G \adjoint \SI_P: \Av_*^{G(\mCO)/P(\mCO)}. \]
More precisely, the right adjoint is the composition
\[\SI_P:= \DMod(\Gr_G)^{U(\mCK)M(\mCO)} \xrightarrow{\oblv} \DMod(\Gr_G)^{P(\mCO)} \xrightarrow{\Av_*} \DMod(\Gr_G)^{G(\mCO)} =: \Sph_G .\]
Its left adjoint is well-defined and equivalent to
\[  - \convolve_{\Sph_G} \Delta_0: \Sph_G\to \SI_P, \]
where $\Delta_0:= \iota_{M,!}(\mathbf{1}_{\Sph_M})\in \SI_P$ is the $0$-standard object. Note that both functors are $\Sph_G$-linear by Appendix \ref{sect-lax-linear}.
\end{constr}

The readers might object the above analogy because the functor $\Eis^\enh_{P\to G}:\mCI(G,P)\to \DMod(\Bun_G)$ is a left adjoint, while the functor $\Av_*^{G(\mCO)/P(\mCO)}: \SI_P\to  \Sph_G $ is a right adjoint. But by Corollary \ref{cor-SIP-SphG-2nd-adj} below, $\Av_*^{G(\mCO)/P(\mCO)}$ is also left adjoint to $\Av_!^{U(\mCK)/U(\mCO)}$ up to a cohomological shift. To prove this, we need the following result:

\begin{lem}[C.f. {\cite[Corollary 6.2.3]{raskin2016chiral}}] \label{lem-semi-inf-vs-parahoric} Consider the parahoric subgroup $Ph_P:= G(\mCO)\mt_{G} P$ of $G(\mCK)$. For any $\mCC$ equipped with a strong $G(\mCK)$-action, the composition
\[ \mCC^{ U(\mCK)M(\mCO) }\xrightarrow{\oblv} \mCC^{P(\mCO)} \xrightarrow{\Av_*} \mCC^{Ph_P} \]
is an equivalence.
\end{lem}

\proof If $P=B$, then $Ph_P$ is the Iwahori subgroup of $G$ and the lemma is just \cite[Corollary 6.2.3]{raskin2016chiral}. The general case can be proved in exactly the same way.

\qed[Lemma \ref{lem-semi-inf-vs-parahoric}]

In particular, we have a pair of inverse functors
\[ \Av_!^{U(\mCK)/U(\mCO)}: \DMod(\Gr_G)^{Ph_P} \adjoint \SI_P:  \Av_*^{Ph_P/P(\mCO)}.  \]
Again, both functors are $\Sph_G$-linear by Appendix \ref{sect-lax-linear}.

\begin{lem} \label{lem-!-avg-vs-oblv-via-Raskin}
The composition
\[ \Sph_G \xrightarrow{ \Av_!^{U(\mCK)/U(\mCO)} } \SI_P \xrightarrow{  \Av_*^{Ph_P/P(\mCO)}}  \DMod(\Gr_G)^{Ph_P} \]
is equivalent to the forgetful functor $\oblv^{G(\mCO)\to Ph_P}$.
\end{lem}

\proof Since both functors are $\Sph_G$-linear, we only need to show they send the unit $\mathbf{1}_{\Sph_G}$ to isomorphic objects, i.e.,
\[ \Av_*^{Ph_P/P(\mCO)}(\Delta_0) \simeq \oblv^{G(\mCO)\to Ph_P}(\mathbf{1}_{\Sph_G}). \]
Applying $\Av_!^{U(\mCK)/U(\mCO)}$, we only need to show
\[ \Delta_0 \simeq \Av_!^{U(\mCK)/U(\mCO)}\circ \oblv^{G(\mCO)\to Ph_P}(\mathbf{1}_{\Sph_G}). \]
But this is obvious.

\qed[Lemma \ref{lem-!-avg-vs-oblv-via-Raskin}]

\begin{cor} \label{cor-SIP-SphG-2nd-adj}
We have adjoint functors
\[ \Av_*^{G(\mCO)/P(\mCO)}[2\dim(G/P)] : \SI_P \adjoint \Sph_G: \Av_!^{U(\mCK)/U(\mCO)} . \]
In other words, the left adjoint and right adjoint of $\Av_!^{U(\mCK)/U(\mCO)}$ are equivalent up to a cohomological shift.
\end{cor}

\proof By Lemma \ref{lem-!-avg-vs-oblv-via-Raskin}, it suffices to prove that the left adjoint and right adjoint of $\oblv^{G(\mCO)\to Ph_P}:\Sph_G\to \DMod(\Gr_G)^{Ph_P}$ are equivalent up to the cohomological shift $[2\dim(G/P)]$. This follows from the fact that $G(\mCO)/Ph_P \simeq G/P$ is smooth and proper.

\qed[Corollary \ref{cor-SIP-SphG-2nd-adj}]

\subsection{Tempered objects in the semi-infinite category \texorpdfstring{$\SI_P$}{ }}
\label{ssec-SI-temper}

In this subsection, we prove the following local analogue of Theorem \ref{thm-iota!*-temperedness}. See also \S \ref{ssec-CTtemmperlocalspectal} for an explanation of its counterpart on the spectral side.

\begin{thm} \label{thm-iota!*-temperedness-local} We have
\begin{itemize}
	\item[(1)] The functors $\iota_M^!:  \SI_P \to \Sph_M$ sends $G$-tempered objects to $M$-tempered objects.

	\item[(2)] The functor $\iota_{M,!}: \Sph_M \to \SI_P$ sends $M$-anti-tempered objects to $G$-anti-tempered objects.

	\item[(3)] The functor $\iota_{M,*}: \Sph_M \to \SI_P$ sends $M$-tempered objects to $G$-tempered objects.

	\item[(4)] The functor $\iota_{M}^*: \SI_P  \to  \Sph_M$ sends $G$-anti-tempered objects to $M$-anti-tempered objects.
\end{itemize}
\end{thm}

Let us first mention the following corollary of the theorem:

\begin{cor} \label{cor-G-temp-imply-M-temp-local}
For any parabolic subgroup $P$ of $G$, any $G$-tempered object in $\SI_{P}$ is $M$-tempered.
\end{cor}

\proof Follows from Theorem \ref{thm-iota!*-temperedness-local}(1) and the fact that $\iota_M^!$ is conservative and $\Sph_M$-linear.

\qed[Corollary \ref{cor-G-temp-imply-M-temp-local}]

\begin{rem} We will also prove the global analogue of the corollary (see Proposition \ref{prop-G-temp-imply-M-temp-global} and remark below it) once we define $M$-tempered objects (or in fact \emph{$P$-tempered} objects) in $\mCI(G,P)$. Such definition is \emph{not} obvious because there is no $\Sph_{M}$-action on $\mCI(G,P)$.
\end{rem}

\proof[Proof of Theorem \ref{thm-iota!*-temperedness-local}] The following lemma can be easily proven by unwinding definitions:

\begin{lem} \label{lem-temper-and-duality} Let $\mCC$ and $\mCD$ be $\Sph_G$-module categories. Suppose $\mCC$ and $\mCD$ are dual to each other and the duality is compatible with the $\Sph_G$-actions. Then 
\begin{itemize}
	\item $\mCC^{G\htemp}$ and $\mCD^{G\hyphen \temp}$ are dual to each other, and the functors
	\[ \oblv: \mCC^{G\htemp}\to \mCC ,\; \; \temp_G: \mCD\to \mCD^{G\htemp} \]
	are dual to each other;
	\item $\mCC^{G\hatemp}$ and $\mCD^{G\hyphen \atemp}$ are dual to each other, and the functors
	\[ \oblv: \mCC^{G\hatemp}\to \mCC ,\; \; \atemp_G: \mCD\to \mCD^{G\hatemp} \]
	are dual to each other.
\end{itemize}
\end{lem}

Let $P^-$ be any parabolic subgroup of $G$ opposite to $P$. By the above lemma and Corollary \ref{cor-inv-inv-duality-functors}, the following statements are equivalent:
\begin{itemize}
	\item The functor $\iota_M^!:  \SI_P \to \Sph_M$ sends $G$-tempered objects to $M$-tempered objects.
	\item The following composition is zero:
	\[ \SI_P^{G\htemp} \to \SI_P\xrightarrow{\iota_M^!} \Sph_M \to \Sph_M^{M\hatemp}.  \]
	\item The following composition is zero:
	\[   \Sph_M^{M\hatemp}\to \Sph_M  \xrightarrow{\iota_{M,!}^-} \SI_{P^-} \to \SI_{P^-}^{G\htemp}.  \]
	\item The functor $\iota_{M,!}^-: \Sph_M \to \SI_{P^-}$ sends $M$-anti-tempered objects to $G$-anti-tempered objects.
\end{itemize}
In other words, (2) for $P^-$ is equivalent to (1) for $P$. Similarly, (3) for $P^-$ is equivalent to (4) for $P$.

Hence we just need to prove (1) and (4). Our strategy is to replace the source category $\SI_P$ by $\Sph_G$. Namely, we deduce (1) and (4) from the following two results, which are proved respectively in \S \ref{ssec-HL-temper-1} and \S \ref{ssec-HL-temper-2}:

\begin{prop} \label{prop-HL!-temper} The composition
	\[ \HL^!_P: \Sph_G \xrightarrow{\Av_!^{U(\mCK)/U(\mCO)}} \SI_P \xrightarrow{\iota_{M}^!} \Sph_M \]
	 sends $G$-tempered objects to $M$-tempered objects.
\end{prop}

\begin{prop} \label{prop-HL*-temper} Consider the composition
	\[ \HL^*_P: \Sph_G \xrightarrow{\Av_!^{U(\mCK)/U(\mCO)}} \SI_P \xrightarrow{\iota_{M}^*} \Sph_M. \]
	An object $\mCF\in \Sph_G$ is $G$-anti-tempered iff $\HL^*_P(\mCF)$ is $M$-anti-tempered.
\end{prop}

\begin{rem} \label{rem-hyperbolic-localization}
The notation $\HL$ stands for \emph{hyperbolic localization} because the underlying object of $\HL^*_P(\mCF)$ on $\Gr_M$ is given by the hyperbolic localization of $\mCF$ in the literature, i.e., the object obtained by first $*$-pullback along $\Gr_P\to \Gr_G$ then $!$-pushforward along $\Gr_P\to \Gr_M$. This can be seen by passing to right adjoints. Namely, when restricted to each connected component $\Gr_{M,\lambda}$ of $\Gr_M$, both right adjoints are
 \[ \DMod(\Gr_{M,\lambda})^{M(\mCO)} \xrightarrow{!\hpull} \DMod(\Gr_{P,\lambda})^{M(\mCO)} \xrightarrow{*\hpush} \DMod(\Gr_G)^{M(\mCO)}\xrightarrow{\Av_*} \DMod(\Gr_G)^{G(\mCO)}.   \]
\end{rem}

\begin{rem} Theorem \ref{thm-EisCT-temperedness}(1) is the global analogue of Proposition \ref{prop-HL!-temper}. We do not know any global analogue of Proposition \ref{prop-HL*-temper}.
\end{rem}

\begin{rem} For the purpose of proving Theorem \ref{thm-iota!*-temperedness-local}, we only need the ``only if'' part of Proposition \ref{prop-HL*-temper}. However, we think the stated stronger claim is of independent interest.
\end{rem}

Let us continue to deduce Theorem \ref{thm-iota!*-temperedness-local}(1),(4) from the above propositions. We explain this deduction for (1), the argument for (4) is completely similar.

Since $\iota_{M,!}: \Sph_M\to \SI_P$ is $\Sph_M$-linear, we have
\[ \iota_{M,!}(\mCG) \simeq \Delta_0 \convolve_{\Sph_M} \mCG  \]
for any $\mCG\in \Sph_M$, and the category $\SI_P$ is generated under colimits by such objects. Let $\mathbf{1}^\temp_{\Sph_G}$ be the tempered unit of $\Sph_G$. Then $\SI_P^{G\mathrm{-}\temp}\subset \SI_P$ is generated under colimits by 
\[\mathbf{1}^\temp_{\Sph_G}\convolve_{\Sph_G}\Delta_0 \convolve_{\Sph_M} \mCG .\] 
Since $\iota_M^!$ is also $\Sph_M$-linear, we have
\[\iota_M^!(\mathbf{1}^\temp_{\Sph_G}\convolve_{\Sph_G}\Delta_0 \convolve_{\Sph_M} \mCG) \simeq \iota_M^!(\mathbf{1}^\temp_{\Sph_G}\convolve_{\Sph_G}\Delta_0) \convolve_{\Sph_M} \mCG, \]
and we only need to show that the RHS is $M$-tempered. Since $\Sph_M^{\temp}\subset \Sph_M$ is a two-sided ideal, we only need to show
\[\iota_M^!(\mathbf{1}^\temp_{\Sph_G}\convolve_{\Sph_G}\Delta_0)\]
is $M$-tempered. But this object is $\HL_P^!( \mathbf{1}^\temp_{\Sph_G} )$, which is $M$-tempered by Proposition \ref{prop-HL!-temper}.

\qed[Theorem \ref{thm-iota!*-temperedness-local}]

\subsection{Hyperbolic localizations and temperedness: I}
\label{ssec-HL-temper-1}

In this subsection, we prove Proposition \ref{prop-HL!-temper}, which states that the functor $\HL^!_P: \Sph_G\to \Sph_M$ preserves temperedness.

As mentioned in the beginning of this section, we need to use the following result of Lafforgue:

\begin{thm}[\!\!{\cite{lafforgue2009quelques}}, see also {\cite[\S 3.1]{beraldo2021tempered}}]
\label{thm-lafforgue}
 Fix an embedding $\mCO  \to \mBP^1$ and identify $\Bun_G(\mBP^1)$ with $G(\mCO)\backslash G(\mCK)/G(\mCO^-)$. Then we have a pair of inverse functors:
\[ \Av_!^{G(\mCO^-)/G}: \Sph_G \adjoint \DMod(\Bun_G(\mBP^1)): \Av_*^{G(\mCO)/G}. \]
More precisely, $\Av_*^{G(\mCO)/G}$ is the composition
\[ \DMod(\Bun_G(\mBP^1)) \simeq \DMod(\Gr_G)^{G(\mCO^-)} \xrightarrow{\oblv} \DMod(\Gr_G)^G\xrightarrow{\Av_*} \DMod(\Gr_G)^{G(\mCO)} =: \Sph_G\]
and $\Av_!^{G(\mCO^-)/G}$ is defined similarly.
\end{thm}

Note that both functors in the theorem are $\Sph_G$-linear by Appendix \ref{sect-lax-linear}.

We first deduce Proposition \ref{prop-HL!-temper} from the following result\footnote{This result seems to be well-known to experts. We provide a proof at the end of this subsection because we cannot find one in the literature.}.

\begin{prop}\label{prop-local-global-CT-P1}
We have a canonical commutative diagram
\[
\xymatrix{
	\Sph_G \ar[d]^-{\simeq}_-{  \Av_!^{G(\mCO^-)/G} } \ar[r]^-{\HL_P^!}  &
	\Sph_M \ar[d]^-{\simeq}_-{  \Av_!^{M(\mCO^-)/M} } \\
	\DMod(\Bun_G(\mBP^1)) \ar[r]^-{\CT_{P,*}} & 
	\DMod(\Bun_M(\mBP^1)).
}
\]
\end{prop}

\begin{warn} \label{warn-D(BunP(P1))-not-SIP}
The functor $\HL_P^!$ factors through $\SI_P$ and $\CT_{P,*}$ factors through $\DMod(\Bun_P(\mBP^1))$, but the category $\SI_P$ is \emph{not} equivalent to $\DMod(\Bun_P(\mBP^1))$.
\end{warn}

\proof[Proof of Proposition \ref{prop-HL!-temper}.] By Proposition \ref{prop-local-global-CT-P1}, we only need to show when $X=\mBP^1$, the functor $\CT_{P,*}$ sends $G$-tempered objects to $M$-tempered objects.

The full subcategory $\DMod(\Bun_G(\mBP^1))^{G\htemp}$ can be explicitly described as follows.

For any curve $X$, let $\DMod(\Bun_G)^{{*\hgen}}\subset \DMod(\Bun_G)$ to be the full subcategory generated under colimits by objects of the form $j_*(\mCF)$, where $j:U\to \Bun_G$ is a \emph{quasi-compact} open substack of $\Bun_G$ and $\mCF\in \DMod(U)$. Then \cite[Theorem 3.5.1]{beraldo2021geometric} says
\[ \DMod(\Bun_G(\mBP^1))^{\temp}\simeq \DMod(\Bun_G(\mBP^1))^{*\hgen}\]
as full subcategories of $\DMod(\Bun_G(\mBP^1)) $.

On the other hand, using the fact that $p:\Bun_P\to \Bun_G$ is quasi-compact when restricted to each connected component of the source, we see that $\CT_{P,*}$ sends $\DMod(\Bun_G)^{{*\hgen}}$ into $\DMod(\Bun_M)^{{*\hgen}}$. Hence, when $X=\mBP^1$, the functor $\CT_{P,*}$ preserves tempered objects by the last paragraph.

\qed[Proposition \ref{prop-HL!-temper}]

Before proving Proposition \ref{prop-local-global-CT-P1}, let us mention the following corollary:

\begin{cor} \label{cor-local-global-Eis-P1}
We have a canonical commutative diagram
\[
\xymatrix{
	\Sph_G \ar[d]^-{\simeq}_-{  \Av_!^{G(\mCO^-)/G} } & &
	\SI_P  \ar[ll]_-{\Av_*^{G(\mCO)/P(\mCO)}[2\dim(G/P)]} &
	\Sph_M \ar[d]^-{\simeq}_-{  \Av_!^{M(\mCO^-)/M} } \ar[l]_-{\iota_{M,!}} \\
	\DMod(\Bun_G(\mBP^1)) & & &
	\DMod(\Bun_M(\mBP^1)) \ar[lll]^-{\Eis_{P,!}}.
}
\]
\end{cor}

\proof Follows from Proposition \ref{prop-local-global-CT-P1} and Corollary \ref{cor-SIP-SphG-2nd-adj}.

\qed[Corollary \ref{cor-local-global-Eis-P1}]

We prove Proposition \ref{prop-local-global-CT-P1} in the rest of this subsection. In spite of Warning \ref{warn-D(BunP(P1))-not-SIP}, we want to relate $\CT_*$ to $\SI_P$. This is achieved in the following construction, which works for any curve $X$. 

\begin{propconstr} \label{propconstr-extended-EisCT} There are canonical $\Sph_{G,x}$-linear adjoint functors
\[ \Eis_{P,\ext,x}: \SI_{P,x}\ot_{\Sph_{M,x}} \DMod(\Bun_M) \adjoint \DMod(\Bun_G): \CT_{P,\ext,x}
\]
such that the composition
\[\DMod(\Bun_M)\simeq \Sph_{M,x}\ot_{\Sph_{M,x}} \DMod(\Bun_M) \adjoint  \SI_{P,x}\ot_{\Sph_{M,x}} \DMod(\Bun_M) \adjoint \DMod(\Bun_G)
\]
is equivalent to $(\Eis_{P,!},\CT_{P,*})$. 
\end{propconstr}

\begin{rem} The functors $\Eis_{P,\ext,x}$ and $\CT_{P,\ext,x}$ were introduced by D. Gaitsgory using the Drinfeld compactification. Unfortunately, his construction is only documented in the case $P=B$ in unpublished notes (see \cite[Talk Ja-4]{arinkinnotes}). Our construction below is different from his original one.
\end{rem}

\proof Consider the $M(\mCO_x)$-torsor $\Bun_M^{\level_{\infty x}}$ on $\Bun_M$ that classifies $M$-torsors on $X$ equipped with a trivialization over $\Spec(\mCO_x)$. Define
\[
\Gr_{G,x}\twisttimes \Bun_M:= \Gr_{G,x}\mt^{M(\mCO_x)} \Bun_M^{\level_{\infty x}}:
\]
this space classifies triples $(\mCF_G,\mCF_M, \alpha)$ where $\mCF_G$ (resp. $\mCF_M$) is a $G$-torsor (resp. $M$-torsor) on $X$ and $\alpha: \mCF_G|_{X\setminus x} \to G\mt^M \mCF_M|_{X\setminus x}$ is an isomorphism of $G$-torsors on $X\setminus x$. We have a map 
\[\pi: \Gr_{G,x}\twisttimes \Bun_M\to \Bun_G\]
that remembers only the $G$-torsor $\mCF_G$. Let $d:=\mathrm{dim.rel.}(\Bun_P,\Bun_M)$ be the locally constant function on $\Bun_M$ whose values are the dimensions of fibers of $\Bun_P\to \Bun_M$. The functor $\Eis_{P,\ext,x}$ is defined to be the (unique) colimit-preserving functor such that
\begin{equation}\label{eqn-defn-ext-eis}
\Eis_{P,\ext,x}( \mCF\boxt_{\Sph_{M,x}} \mCM ) := \pi_!( \mCF\widetilde\boxtimes \mCM  )[-2d], 
\end{equation}
for any $\mCF\in \SI_{P,x}$ and $\mCM \in \DMod(\Bun_M)$. To justify well-definedness, we need to check:
\begin{itemize}
	\item[(a)] $\pi_!$ is well-defined on $\mCF\widetilde\boxtimes \mCM $;
	\item[(b)] For any $\mCK\in \Sph_{M,x}$, we have canonical equivalence
	\[\Eis_{P,\ext,x}( (\mCF\convolve_{\Sph_{M,x}} \mCK) \boxt_{\Sph_{M,x}} \mCM ) \simeq \Eis_{P,\ext,x}( \mCF  \boxt_{\Sph_{M,x}}(\mCK \convolve_{\Sph_{M,x}}  \mCM ) ).\]
\end{itemize}

We first prove (a). Since $\SI_{P,x}$ is generated under colimits by the image of $\iota_{M,x,!}$, we can assume $\mCF \simeq \iota_{M,x,!}(\mCG)$ for $\mCG\in \Sph_{M,x}$. Consider the following diagram
\[
\xymatrix{
	\Gr_{G,x}\twisttimes \Bun_M  \ar[d]^-\pi
	& \Gr_{P,x}\twisttimes \Bun_M \ar[l]_-{u} \ar[r]^-{v} \ar[d]^-{\rho}
	& \Gr_{M,x}\twisttimes \Bun_M \ar[d]^-{\sigma}
	\\ 
	\Bun_G 
	& \Bun_P \ar[l]_-p \ar[r]^-q
	& \Bun_M.
}
\]
A diagram chase shows
\[  \pi_!( \mCF\widetilde\boxtimes \mCM  ) \simeq  p_!\circ \rho_! \circ v^!(  \mCG\widetilde\boxtimes \mCM).\]
Consider the Cartesian square
\[
	\begin{tikzcd}
    \Bun_P\times_{\Bun_M}( \Gr_{M,x}\twisttimes \Bun_M )
    \arrow[r,"q'"] \arrow[d,"\sigma'"]
    \arrow[dr, phantom, "\lrcorner", very near start]
    & \Gr_{M,x}\twisttimes \Bun_M \arrow[d,"\sigma"]  \\
      \Bun_P\arrow[r,"q"]
    & \Bun_M.
	\end{tikzcd}
	\]
By Lemma \ref{lem-UHC-for-U-modification} below, the map
\[\theta: \Gr_{P,x}\twisttimes \Bun_M  \to \Bun_P\mt_{\Bun_M}( \Gr_{M,x}\twisttimes \Bun_M )  \]
is universally homologically contractible. In particular, $\theta^!$ is fully faithful and hence $\theta_!\circ\theta^! \simeq \Id$. 
It follows that
\[  \pi_!( \mCF\widetilde\boxtimes \mCM  ) \simeq  p_!\circ \rho_! \circ v^!(  \mCG\widetilde\boxtimes \mCM) \simeq p_!\circ (\rho')_! \circ \theta_! \circ\theta^! \circ {q'}^!(  \mCG\widetilde\boxtimes \mCM) \simeq p_!\circ (\rho')_! \circ {q'}^!(  \mCG\widetilde\boxtimes \mCM).\]
Since $\Bun_P\to \Bun_M$ is smooth, the base-change isomorphism yields
\[  \pi_!( \mCF\widetilde\boxtimes \mCM  ) \simeq  p_!\circ (\rho')_! \circ {q'}^!(  \mCG\widetilde\boxtimes \mCM) \simeq  p_!\circ q^!\circ \sigma_!( \mCG\widetilde\boxtimes \mCM  ) .\]
Note that the RHS is well-defined because $\sigma$ is ind-proper and $p_!\circ q^! \simeq p_!\circ q^*[2d] \simeq \Eis_{P,!}[2d]$ is well-defined. This proves (a), as well as the last claim in the proposition (modulo claim (b)). Moreover, note that if $\mCG$ and $\mCM$ are compact, then so is $\pi_!( \mCF\widetilde\boxtimes \mCM  )$.

Now (b) follows from the fact that the diagonal $M(\mCK_x)$-action on $\Gr_{G,x} \times \Bun_M^{\level_{\infty x}}$ stabilizes the fibers of $\pi$.

By the previous discussion, the functor $\Eis_{P,\ext,x}$ preserves compact objects, hence its right adjoint $\CT_{P,\ext,x}$ is continuous. Also, $\Eis_{P,\ext,x}$ has a natural $\Sph_{G,x}$-linear structure because $\pi$ is $G(\mCK_x)$-equivariant. Then Appendix \ref{sect-lax-linear} guarantees that $\CT_{P,\ext,x}$ is strictly $\Sph_{G,x}$-linear.

\qed[Proposition-Construction \ref{propconstr-extended-EisCT}]

The following lemma was used above.

\begin{lem}\label{lem-UHC-for-U-modification} The map
\[\theta: \Gr_{P,x}\twisttimes \Bun_M  \to \Bun_P\mt_{\Bun_M}( \Gr_{M,x}\twisttimes \Bun_M )  \]
induced by the obvious map $\Gr_{P,x} \to \Bun_P\mt_{\Bun_M} \Gr_{M,x}$ is universally homologically contractible.
\end{lem}

\proof The proof is similar to that of Lemma \ref{lem-UHC-unipotent-hecke}.

For any finite type affine test scheme $s: S\to \Bun_P\mt_{\Bun_M}( \Gr_{M,x}\twisttimes \Bun_M )$, we need to show the base-change $\theta': S'\to S$ of $\theta$ along $s$ is homologically contractible. By definition, $s$ corresponds to $(\mCF_P,\mCF_M,\alpha)$ where $\mCF_P$ (resp. $\mCF_M$) is a $P$-torsor (resp. $M$-torsor) on $X\times S$ and 
\[\alpha: M\mt^P \mCF_P|_{(X\setminus x)\times S} \to \mCF_M|_{(X\setminus x)\times S}  \] is an isomorphism of $M$-torsors on $(X\setminus x)\times S$. Hence, for any finite type affine test $S$-scheme $T$, the groupoid (which is actually a set) $\Maps_{S}(T,S')$ classifies isomorphisms 
\[\beta: \mCF_P|_{(X\setminus x) \times T} \to P\mt^M \mCF_M|_{(X\setminus x)\times T} \] of $P$-torsors on $(X\setminus x)\times T$ such that $M\mt^P \beta$ is equal to the composition of $\alpha|_{(X\setminus x)\times T}$ with $\mCF_M|_{(X\setminus x)\times T} \simeq M\mt^P (P\mt^M \mCF_M)|_{(X\setminus x)\times T} $. It follows that $\Maps_{S}(T,S')$ is the set of $P$-equivariant maps $\mCF_P|_{(X\setminus x) \times T} \to P/M$. Since $(X\setminus x)\times T$ is affine, any $P$-torsor on it has an $M$-structure. Hence $\Maps_{S}(T,S')$ is a non-empty torsor for the group $\Maps((X\setminus x)\times T , P/M )$. Since $P/M$ is an affine space and $X\setminus x$ is affine, the functor $T\mapsto \Maps((X\setminus x)\times T , P/M )$ is represented by $S\times \mBA^\infty$, where $\mBA^\infty$ is the infinite dimensional affine space associated to the vector space $\Gamma(X\setminus x, \mCO)^{\oplus \dim (P/M)}$. So we have a free and transitive action of $\mBA^\infty$ on the fibers of $\theta':S'\to S$, which implies $\theta
'$ is homologically contractible.

\qed[Lemma \ref{lem-UHC-for-U-modification}]

We also need the following lemma:

\begin{lem} \label{lem-SIP-global-model}
We have a pair of inverse functors:
\[ \Av_!^{M(\mCO^-)/M}: \SI_P \adjoint \DMod(\Gr_G)^{U(\mCK)M(\mCO^-)}: \Av_*^{M(\mCO)/M}. \]
More precisely, $\Av_*^{M(\mCO)/M}$ is the composition
\[ \DMod(\Gr_G)^{U(\mCK)M(\mCO^-)} \xrightarrow{\oblv} \DMod(\Gr_G)^{U(\mCK)M}\xrightarrow{\Av_*} \DMod(\Gr_G)^{U(\mCK)M(\mCO)} =: \SI_P\]
and $\Av_!^{M(\mCO^-)/M}$ is defined similarly.
\end{lem}

\proof We have adjoint pairs similar to
\[
\iota_{M,!}: \DMod(\Gr_M)^{M(\mCO)} \adjoint \DMod(\Gr_G)^{U(\mCK)M(\mCO)} : \iota_M^!,
\]
but with $M(\mCO)$-invariance replaced by $M(\mCO^-)$-invariance (resp. $M$-invariance). We abuse notation by writing these adjoint pairs as $(\iota_{M,!},\iota_M^!)$. By construction, we have a commutative diagram
\[
\xymatrix{
	\DMod(\Gr_G)^{U(\mCK)M(\mCO^-)} \ar[r]^-\oblv \ar[d]^-{\iota_M^!} &
	\DMod(\Gr_G)^{U(\mCK)M} \ar[r]^-{\Av_*}\ar[d]^-{\iota_M^!}  &
	\DMod(\Gr_G)^{U(\mCK)M(\mCO)}\ar[d]^-{\iota_M^!}  \\
	\DMod(\Gr_M)^{M(\mCO^-)} \ar[r]^-\oblv &
	\DMod(\Gr_M)^{M} \ar[r]^-{\Av_*} &
	\DMod(\Gr_M)^{M(\mCO)}
}
\]
that is left adjointable along vertical directions. Indeed, the left adjointability follows from the fact that $\iota_{M,!}:\DMod(\Gr_M) \to \DMod(\Gr_G)^{U(\mCK)}$ is $M(\mCK)$-linear. Recall the composition of the bottom row is an equivalence: it is just Lafforgue's equivalence for the reductive group $M$. Hence the composition of the top row is also an equivalence by the Barr--Beck--Lurie theorem on monadicity.

\qed[Lemma \ref{lem-SIP-global-model}]

\proof[Proof of Proposition \ref{prop-local-global-CT-P1}.] By Proposition-Construction \ref{propconstr-extended-EisCT}, we only need to construct the following commutative diagram
\[
\xymatrix{
	\Sph_G \ar[d]^-{\simeq}_-{  \Av_!^{G(\mCO^-)/G} } \ar[r]^-{\Av_!^{U(\mCK)/U(\mCO)}} &
	\SI_P \ar[d]^-\simeq_-{\Id \ot \Av_!^{M(\mCO^-)/M}   } \ar[r]^-{\iota_P^!} &
	\Sph_M \ar[d]^-{\simeq}_-{  \Av_!^{M(\mCO^-)/M} } \\
	\DMod(\Bun_G(\mBP^1)) \ar[r]_-{\CT_{P,*,\ext}} & 
	\SI_P\ot_{\Sph_M} \DMod(\Bun_M(\mBP^1)) \ar[r]_-{\iota_M^! \ot\Id}
	&
	\DMod(\Bun_M(\mBP^1)).
}
\]
The right square is obvious. For the left one, recall $\Eis_{P,!,\ext}$ is defined by the formula
\[ \Eis_{P,!,\ext}( \mCF\boxt_{\Sph_{M,x}} \mCM ) := \pi_!( \mCF\widetilde\boxtimes \mCM  ), \]
where $\mCF\in \SI_P$ and $\mCM\in \DMod(\Bun_M(\mBP^1))\simeq \DMod(\Gr_M/M(\mCO^-))$, so that
\[ \mCF\widetilde\boxtimes \mCM \in \DMod( \Gr_G\twisttimes \Bun_M(\mBP^1) ) \simeq \DMod(\Gr_G \mt^{M(\mCO)} M(\mCK) )^{M(\mCO^-)}. \]
Via the above equivalence, $\pi_!$ is the partially defined functor
\[ \DMod(\Gr_G \mt^{M(\mCO)} M(\mCK) )^{M(\mCO^-)} \xrightarrow{\mathrm{push}} \DMod(\Gr_G )^{M(\mCO^-)} \xrightarrow{\Av_!} \DMod(\Gr_G )^{G(\mCO^-)} \simeq \DMod(\Bun_G(\mBP^1)).\]
Since $\mCF$ is $U(\mCK)$-equivariant, the image of $\mCF\widetilde\boxtimes \mCM$ under the first functor in the above composition is $U(\mCK)$-equivariant. In particular, $\Eis_{P,!,\ext}$ factors as
\[\SI_P\ot_{\Sph_M} \DMod(\Bun_M(\mBP^1)) \xrightarrow{c}  \DMod(\Gr_G )^{U(\mCK)M(\mCO^-)} \xrightarrow{\Av_!^{G(\mCO^-)/P(\mCO^-)}} \DMod(\Bun_G(\mBP^1)),\]
where $\Av_!^{G(\mCO^-)/P(\mCO^-)}$ is the partially defined functor
\[  \DMod(\Gr_G )^{U(\mCK)M(\mCO^-)} \xrightarrow{\oblv} \DMod(\Gr_G)^{P(\mCO^-)} \xrightarrow{\Av_!} \DMod(\Gr_G)^{G(\mCO^-)}\simeq \DMod(\Bun_G(\mBP^1)). \]
It is easy to check the composition
\[ \SI_P \simeq \SI_P\ot_{\Sph_M} \DMod(\Bun_M(\mBP^1))\xrightarrow{c}  \DMod(\Gr_G )^{U(\mCK)M(\mCO^-)}  \]
is equivalent to the functor $\Av_!^{M(\mCO^-)/M}$ in Lemma \ref{lem-SIP-global-model}. In particular, $c$ is an equivalence and therefore $\Av_!^{G(\mCO^-)/P(\mCO^-)}$ is well-defined. It follows $\CT_{P,*,\ext}$ is equivalent to
\[ \DMod(\Bun_G(\mBP^1)) \xrightarrow{\Av_*^{U(K)/U(\mCO^-)}} \DMod(\Gr_G )^{U(\mCK)M(\mCO^-)}  \xrightarrow{c^{-1}} \SI_P\ot_{\Sph_M} \DMod(\Bun_M(\mBP^1)).  \]
Hence, it remains to check that the following diagram is naturally commutative:
\[
\xymatrix{
	\Sph_G \ar[d]^-{\simeq}_-{  \Av_!^{G(\mCO^-)/G} } \ar[r]^-{\Av_!^{U(\mCK)/U(\mCO)}} &
	\SI_P  \\
	\DMod(\Bun_G(\mBP^1)) \ar[r]_-{\Av_*^{U(K)/U(\mCO^-)}} & 
	 \DMod(\Gr_G )^{U(\mCK)M(\mCO^-)} \ar[u]^-\simeq_-{\Av_*^{M(\mCO)/M}   } .
}
\]
It is easy to check that the composition
\[
\DMod(\Bun_G(\mBP^1)) \xrightarrow{\Av_*^{U(K)/U(\mCO^-)}}  \DMod(\Gr_G )^{U(\mCK)M(\mCO^-)} \xrightarrow{\Av_*^{M(\mCO)/M} } \SI_P \xrightarrow{\Av_*^{Ph_P/P(\mCO)}} \DMod(\Gr_G)^{Ph_P}
\]
is equivalent to
\[\Av_*^{Ph_P/P}: \DMod(\Bun_G(\mBP^1)) \simeq \DMod(\Gr_G)^{G(\mCO^-)} \xrightarrow{\oblv} \DMod(\Gr_G)^{P} \xrightarrow{\Av_*} \DMod(\Gr_G)^{Ph_P}. \]
Hence, by Lemma \ref{lem-semi-inf-vs-parahoric} and Lemma \ref{lem-!-avg-vs-oblv-via-Raskin}, we only need to show that
\[\Sph_G \xrightarrow{ \Av_!^{G(\mCO^-)/G} } \DMod(\Bun_G(\mBP^1)) \xrightarrow{\Av_*^{Ph_P/P}}  \DMod(\Gr_G)^{Ph_P} \]
is equivalent to the forgetful functor. Since $G(O)/G \simeq Ph_P/P$, the functor $\Av_*^{Ph_P/P}$ factors as
\[\DMod(\Bun_G(\mBP^1)) \simeq \DMod(\Gr_G)^{G(\mCO^-)} \xrightarrow{\Av_*^{G(\mCO)/G}} \DMod(\Gr_G)^{G(\mCO)} \xrightarrow{\oblv} \DMod(\Gr_G)^{Ph_P}.\]
Then we are done because $\Av_!^{G(\mCO^-)/G}$ and $\Av_*^{G(\mCO)/G}$ are inverse to each other.

\qed[Proposition \ref{prop-local-global-CT-P1}]

\subsection{Hyperbolic localizations and temperedness: II}
\label{ssec-HL-temper-2}

In this subsection, we prove Proposition \ref{prop-HL*-temper}, which says that the functor $\HL^*_P: \Sph_G\to \Sph_M$ detects anti-temperedness.

Consider the Borel subgroup $B\cap M$ of $M$ and the corresponding functor $\HL^*_{B_M}:\Sph_M \to \Sph_T$. We claim $\HL^*_B \simeq \HL^*_{B_M}\circ \HL^*_{P}$. Indeed, this follows from Remark \ref{rem-hyperbolic-localization}, the base-change isomorphism, and the identification $\Gr_{B_M} \mt_{\Gr_M} \Gr_P \simeq \Gr_B$. 

Hence, to prove the proposition, we can assume $P=B$ is the Borel subgroup. Recall that any object in $\Sph_T$ is $T$-tempered. Hence we only need to show that an object $\mCF\in \Sph_G$ is $G$-anti-tempered iff $\HL^*_B(\mCF)\simeq 0$. Equivalently, by Lemma \ref{lem-atemp-t-structure}, we only need to show that an object $\mCF\in \Sph_G$ has vanishing cohomologies (with respect to the perverse t-structure) iff $\HL^*_B(\mCF)\simeq 0$. 

Consider the locally constant function $d$ on $\Gr_T$ whose value on the connected component $\Gr_{T,\lambda}$ is $\langle 2\check \rho,\lambda\rangle$. It is well-known (see e.g. \cite{mirkovic2007geometric}) that for $\mCM\in \Sph_G^{\heartsuit}$, the object $\HL^*_B(\mCM)[d]$ is contained in $\Sph_T^{\heartsuit}$, and that the resulting functor $\Sph_G^{\heartsuit} \to \Sph_T^{\heartsuit}$ can be identified with the forgetful functor $\Rep(\check G)\to \Rep(\check T)$.

We claim this formally implies that the functor $\HL^*_B(-)[d]: \Sph_G\to \Sph_T$ is right t-exact. To prove the claim, recall that any compact object in $\Sph_G$ has bounded cohomologies and that the perverse t-structure on $\Sph_G$ is compatible with filtered colimits. Now, any object $\mCF\in \Sph_G^{\le 0}$ can be written as a filtered colimit $\mCF \simeq \colim_i F_i$ with $F_i$ compact. Then $\mCF \simeq \tau^{\le 0}\mCF \simeq \colim_i \tau^{\le 0} F_i$. Since $\tau^{\le 0} F_i$ has bounded cohomologies, it is a finite colimit of $H^{-j}(\tau^{\le 0} F_i)[j]$ for $j\ge 0$. Hence the functor $\HL^*_B(-)[d]$ sends each $\tau^{\le 0} F_i$ to a finite colimit of objects in $\Sph_T^{\le 0}$. It follows that $\HL^*_B(\mCF)[d] \in \Sph_T^{\le 0}$ as desired. This proves the claim.

Combining the claim with the fact that the functor $\Sph_G^{\heartsuit} \to \Sph_T^{\heartsuit}$ is conservative, we see that $\mCF$ has vanishing cohomologies iff $\HL^*_B(\mCF)$ does. But the latter condition is equivalent to $\HL^*_B(\mCF)\simeq 0$ because $\Sph_T^{\le -\infty}$ is the zero category.

\qed[Proposition \ref{prop-HL*-temper}]

\subsection{Proof of Theorem \ref{thm-iota-temperedness}: a local-to-global argument}
\label{ssec-SI-to-IGP-temper}

In this subsection, we deduce Theorem \ref{thm-iota!*-temperedness} from Theorem \ref{thm-iota!*-temperedness-local}. The key ingredient is the following construction.

\begin{propconstr} \label{propconstr-extended-EisCT-gen} There are canonical $\Sph_{G,x}$-linear adjoint functors
\[ \Eis_{P,\ext,x}^{P\hgen}: \SI_{P,x}\ot_{\Sph_{M,x}} \DMod(\Bun_M) \adjoint \mCI(G,P): \CT_{P,\ext,x}^{P\hgen}
\]
such that the composition
\[\DMod(\Bun_M)\simeq \Sph_{M,x}\ot_{\Sph_{M,x}} \DMod(\Bun_M) \adjoint  \SI_{P,x}\ot_{\Sph_{M,x}} \DMod(\Bun_M) \adjoint \mCI(G,P)
\]
is equivalent to $(\iota_{M,!},\iota_M^!)$. 
\end{propconstr}

\proof The proof is similar to that of Proposition-Construction \ref{propconstr-extended-EisCT}, where we defined $\Sph_{G,x}$-linear adjoint functors 
\[ \Eis_{P,\ext,x}: \SI_{P,x}\ot_{\Sph_{M,x}} \DMod(\Bun_M) \adjoint \DMod(\Bun_G): \CT_{P,\ext,x}
\]
using $!$-pushforward along the map $\pi: \Gr_{G,x}\twisttimes \Bun_M\to \Bun_G$. By definition, the map $\pi$ factors as
\[ \Gr_{G,x}\twisttimes \Bun_M\to \Bun_G^{P\hgen} \to \Bun_G. \]
Hence $\Eis_{P,\ext,x}$ factors through as
\[ \SI_{P,x}\ot_{\Sph_{M,x}} \DMod(\Bun_M)\to \DMod(\Bun_G^{P\hgen}) \xrightarrow{!\hpush}\DMod(\Bun_G). \]
Consider the composition
\[ \DMod(\Bun_M)\to \SI_{P,x}\ot_{\Sph_{M,x}} \DMod(\Bun_M)\to \DMod(\Bun_G^{P\hgen}). \]
As in the proof of Proposition-Construction \ref{propconstr-extended-EisCT}, one can show that this functor is equivalent to 
\[ \DMod(\Bun_M) \xrightarrow{*\hpull}\DMod(\Bun_P) \xrightarrow{!\hpush}\DMod(\Bun_G^{P\hgen}) .\]
Hence, by Lemma \ref{lem-generator-IGP}, it factors through
\[ \DMod(\Bun_M) \xrightarrow{\iota_{M,!}}\mCI(G,P) . \]
We thus have shown that $\Eis_{P,\ext,x}$ factors through a functor
\[ \SI_{P,x}\ot_{\Sph_{M,x}} \DMod(\Bun_M)\to \mCI(G,P), \]
which is the desired functor $\Eis_{P,\ext,x}^{P\hgen}$. We let $\CT_{P,\ext,x}^{P\hgen}$ be its continuous right adjoint.

\qed[Proposition-Construction \ref{propconstr-extended-EisCT-gen}]

Using Proposition-Construction \ref{propconstr-extended-EisCT-gen}, Theorem \ref{thm-iota!*-temperedness}(1) and (2) obviously follow from Theorem \ref{thm-iota!*-temperedness-local}(1) and (2). 

To proceed, we need one more geometric input. We claim we have the following diagram which is Cartesian at least at the reduced level:
\[
  \begin{tikzcd}
   \Gr_{P,x}\widetilde\times \Bun_M
    \arrow[r] \arrow[d]
    \arrow[dr, phantom, "\lrcorner", very near start]
    &\Gr_{G,x}\widetilde\times \Bun_M \arrow[d]  \\
      \Bun_P \arrow[r]
    & \Bun_G^{P\hgen}.
  \end{tikzcd}
\]
Namely, unwinding the definitions, the desired fiber product classifies triples consisting of:
\begin{itemize}
	\item[(i)] an $M$-torsor $\mCF_M$ on $X$;
	\item[(ii)] a $P$-torsor $\mCF_P$ on $X$;
	\item[(iii)] an isomorphism $\alpha: G\mt^P\mCF_P|_{X\setminus x}\to G\mt^M \mCF_M|_{X\setminus x}$ between $G$-torsors on $X\setminus x$
\end{itemize}
with the following property: there exists an open subset $U\subset X$ such that $\alpha$ induces an isomorphism $P\mt^P\mCF_P|_U \simeq P\mt^M \mCF_M|_U$ of $P$-torsors on $U$. Since $G/P$ is separated, this last equivalence can be uniquely extended over $X\setminus x$. Meanwhile, such extended isomorphism also determines $\alpha$. Hence the above data are exactly classified by $\Gr_{P,x}\widetilde\times \Bun_M$. This proves the claim.

Using the base-change isomorphism for the above Cartesian square, up to a cohomological shift (inherited by that in (\ref{eqn-defn-ext-eis})), the functor
\begin{equation} \label{eqn-thm-iota!*-temperedness-1}
\SI_{P,x}\ot_{\Sph_{M,x}} \DMod(\Bun_M)  \xrightarrow{ \Eis_{P,\ext,x}^{P\hgen}} \mCI(G,P) \xrightarrow{\iota_{M}^*} \DMod(\Bun_M) 
\end{equation}
is equivalent to
\begin{equation} \label{eqn-thm-iota!*-temperedness-2}
\SI_{P,x}\ot_{\Sph_{M,x}} \DMod(\Bun_M)  \xrightarrow{ \iota_{M}^*\ot \Id} \Sph_{M,x}\ot_{\Sph_{M,x}} \DMod(\Bun_M) \simeq \DMod(\Bun_M). 
\end{equation}
We also have a similar right adjoint version of this equivalence.

Now one can again deduce Theorem \ref{thm-iota!*-temperedness}(3)(4) from Theorem \ref{thm-iota!*-temperedness-local}(3)(4). We explain this for (3) and leave (4) to the readers. Note that the $\Sph_{G,x}$-linear functor 
\[ \CT_{P,\ext,x}^{P\hgen}: \mCI(G,P)\to  \SI_{P,x}\ot_{\Sph_{M,x}} \DMod(\Bun_M) \]
is conservative because it is a factor of the conservative functor $\iota_M^!$. Hence an object in $\mCI(G,P)$ is $G$-tempered (resp. $G$-anti-tempered) iff it is sent to a $G$-tempered object (resp. $G$-anti-tempered) by the functor $\CT_{P,\ext,x}^{P\hgen}$. Passing to the left adjoints, we see that $\mCI(G,P)^\temp$ (resp. $\mCI(G,P)^\atemp$) is generated under colimits by images of $G$-tempered objects (resp. $G$-anti-tempered objects) under $\Eis_{P,\ext,x}^{P\hgen}$. Hence, to prove Theorem \ref{thm-iota!*-temperedness}(4), we only need to show that the composition (\ref{eqn-thm-iota!*-temperedness-1}) sends $G$-anti-tempered objects to $M$-anti-tempered objects. By the last paragraph, we only need to prove this for (\ref{eqn-thm-iota!*-temperedness-2}). But this follows immediately from Theorem \ref{thm-iota!*-temperedness-local}(4).

\qed[Theorem \ref{thm-iota!*-temperedness}]


\subsection{A (conjectural) spectral explanation of the local statement}
\label{ssec-CTtemmperlocalspectal}

Assuming certain conjectures (which are actually undocumented known results), Theorem \ref{thm-iota!*-temperedness} can also be proved by translating to the spectral side. For completeness, we provide the argument in this remark but assume the readers are familiar with some works on the subject, and in particular with the theory of singular support as in \cite{arinkin2015singular}.

The main theorem of \cite{arkhipov2004quantum} provides an equivalence 
\[\DMod(\Gr_G)^I_{\mathrm{loc.c}} \simeq \Coh((\pt\mt_{\check \mfg} \check \mfn)/{\check B}),\] 
where $I:=G(\mCO)\mt_G B$ is the Iwahori subgroup and $\DMod(\Gr_G)^I_{\mathrm{loc.c}}\subset \DMod(\Gr_G)^I$ is the full subcategory of \emph{locally compact objects}, i.e., those objects whose image in $\DMod(\Gr_G)$ are compact. Therefore 
\[\Ind(\DMod(\Gr_G)^I_{\mathrm{loc.c}} ) \simeq \IndCoh(  (\pt\mt_{\check \mfg} \check \mfn)/{\check B}).\]
The category $\DMod(\Gr_G)^I$ can be obtained from the LHS by a \emph{renormalization} construction, which corresponds to a modification of \emph{singular supports} on the RHS (see \cite[\S 12]{arinkin2015singular} for a similar construction for $\Sph_G$). More precisely, we have
\[ \DMod(\Gr_G)^I \simeq \IndCoh_{\Nilp(\check \mfb)/\check B}(  (\pt\mt_{\check \mfg} \check \mfn)/{\check B}),\]
where\footnote{In this formula, we fix a $\check G$-equivariant self-duality $\check \mfg\simeq \check \mfg^*$.} $\Nilp(\check \mfb)/\check B \subset \check \mfb/\check B  \simeq \Sing( (\pt\mt_{\check \mfg} \check \mfn)/{\check B})$. By Raskin's equivalence $\SI_B \simeq \DMod(\Gr_G)^I$ (see {\cite[Corollary 6.2.3]{raskin2016chiral}}), we obtain a spectral description for $\SI_B$. This description can be generalized to any parabolic subgroups\footnote{During the preparation of this paper, J. Campbell and S. Raskin announced a new proof of the derived Satake equivalence using chiral algebras and their method can also prove this conjecture.} $P$:

\begin{conj} \label{conj-parablic-abg}
For any standard parabolic subgroup $P$ of $G$, let $\check P$ be the corresponding parabolic subgroup of $\check G$. Then there is a canonical equivalence
\[ \SI_P \simeq \IndCoh_{\Nilp(\check \mfp)/\check P}(  (\pt\mt_{\check \mfg} \check \mfu)/{\check P}) \simeq \IndCoh_{\Nilp(\check \mfp)/\check P}(  ( \pt\mt_{\check \mfg} \check \mfp \mt_{\check \mfm} \pt)/{\check P}),\]
where $\check \mfu$ is the Lie algebra of the unipotent radical of $\check P$ and $ \Nilp(\check \mfp)/\check P \subset \check \mfp/\check P \simeq \Sing( (\pt\mt_{\check \mfg} \check \mfu)/{\check P})$. This equivalence is required to satisfy the following conditions:
\begin{itemize}
	\item It is compatible with the actions of the geometric Satake \emph{monoidal} equivalences
\[ \Sph_G \simeq \IndCoh_{\Nilp(\check \mfg)/\check G}(  (\pt\mt_{\check \mfg} \pt)/{\check G}),\; \Sph_M \simeq \IndCoh_{\Nilp(\check \mfm)/\check M}(  (\pt\mt_{\check \mfm} \pt)/{\check M}); \]

	\item The standard object $\Delta_0\in \SI_P$ corresponds to the image of the trivial representation under the functor
\[ \Rep(\check P) \simeq \IndCoh( \pt/\check P ) \xrightarrow{\mathrm{push}} \IndCoh(  (\pt\mt_{\check \mfg} \check \mfu)/{\check P}) \xrightarrow{\Psi} \IndCoh_{\Nilp(\check \mfp)/\check P}(  (\pt\mt_{\check \mfg} \check \mfu)/{\check P}).\]
\end{itemize}
\end{conj}

Using the convolution actions on $\Delta_0$, it is easy to see

\begin{cor} \label{cor-spectral-iota} Assuming Conjecture \ref{conj-parablic-abg}, we have
\begin{itemize}
	\item[(1)] The functor $\iota_{M,!}: \Sph_M \to \SI_P $ corresponds to $!$-pull-$*$-push along the following diagram\footnote{We always mean $!$-pullback along the right arm and then $*$-pushforward along the left arm. Part of the claim says such $!$-pull-$*$-push functors send objects with given singular supports to objects with desired singular supports.
	}
	\[  ( \pt\mt_{\check \mfg} \check \mfp \mt_{\check \mfm} \pt)/{\check P} \gets   (\pt\mt_{\check \mfm} \pt)/{\check P} \to  (\pt\mt_{\check \mfm} \pt)/\check M \]

	\item[(2)] Up to tensoring with a line bundle on $(\pt\mt_{\check \mfm} \pt)/\check M$, the functor $\iota_{M}^!:  \SI_P \to \Sph_M$ corresponds to $!$-pull-$*$-push along the following diagram
	\[  (\pt\mt_{\check \mfm} \pt)/\check M \gets   (\pt\mt_{\check \mfm} \pt)/{\check P} \to ( \pt\mt_{\check \mfg} \check \mfp \mt_{\check \mfm} \pt)/{\check P} ; \]

	\item[(3)] The functor $\iota_{M}^*: \SI_P\to\Sph_M $ corresponds to $!$-pull-$*$-push along the following diagram
	\[   (\pt\mt_{\check \mfm} \pt)/\check M \gets   (\pt\mt_{\check \mfp^-} \pt)/{\check M} \to  (\pt\mt_{\check \mfg} \check \mfu)/{\check P} \]
	
	\item[(4)] Up to tensoring with a line bundle on $(\pt\mt_{\check \mfm} \pt)/\check M$, the functor $\iota_{M,*}: \Sph_M\to\SI_P $ corresponds to $!$-pull-$*$-push along the following diagram
	\[  (\pt\mt_{\check \mfg} \check \mfu)/{\check P}  \gets   (\pt\mt_{\check \mfp^-} \pt)/{\check M} \to  (\pt\mt_{\check \mfm} \pt)/\check M   \]

	\item[(5)] The functor $\Av_{!}^{U(\mCK)/U(\mCO)}: \Sph_G \to \SI_P $ corresponds to $!$-pull-$*$-push along the following diagram
	\[  ( \pt\mt_{\check \mfg} \check \mfp \mt_{\check \mfm} \pt)/{\check P} \gets   (\pt\mt_{\check \mfg} \pt)/{\check P} \to  (\pt\mt_{\check \mfg} \pt)/\check G; \]

	\item[(6)] The functor $\Av_{*}^{G(\mCO)/P(\mCO)}:  \SI_P \to \Sph_G $ corresponds to $!$-pull-$*$-push along the following diagram
	\[ (\pt\mt_{\check \mfg} \pt)/\check G \gets   (\pt\mt_{\check \mfg} \pt)/{\check P} \to  ( \pt\mt_{\check \mfg} \check \mfp \mt_{\check \mfm} \pt)/{\check P} .\]
\end{itemize}

\end{cor}

Then one can deduce Theorem \ref{thm-iota!*-temperedness-local} from the above result by calculating the singular supports using singular differentials.

\section{The building blocks \texorpdfstring{$\mCW(G,P\supset Q)$}{ }} \label{sect-WGPQ}

In this section, we define and study the building blocks $\mCW(G,P\supset Q)$ of the automorphic gluing category. We fix a closed point $x\in X$ and an identification $\mCO_x \simeq \mCO$. Definitions in this section are supposed to be independent of such choices\footnote{We will prove this in a future paper.}.

In \S \ref{ssec-IGP-functorial}, we extend the functorialities of the parabolic categories $\mCI(G,P)$ in such a way that the functors $\iota_{M,!}:\DMod(\Bun_M)\to \mCI(G,P)$ and $\Eis^\enh_{Q\to P}: \mCI(G,Q)\to \mCI(G,P)$ are instances of the same construction. These new functorialities will be used in the definition of $\mCW(G,P\supset Q)$.

In \S \ref{ssec-WGPQ}, for an inclusion of parabolic subgroups $Q\subset P$ of $G$, we define the notion of $P$-tempered objects in the parabolic category $\mCI(G,Q)$, and define $\mCW(G,P\supset Q)$ as the full subcategory of these $P$-tempered objects. We study their functorial properties: these are the main ingredient for the construction of the automorphic gluing category. There are two kinds of functorialities and they are adjoint to each other: for the left adjoint one, $\mCW(G,P\supset Q)$ is covariant in $Q$ and $G$, but contravariant in $P$. 

In \S \ref{ssec-IGP-proof-2}-\ref{ssec-IGP-proof-5}, we provide proofs for results in \S \ref{ssec-WGPQ}. We suggest that the readers skip them first and proceed to \S \ref{sect-aut-gluing-thm}, where we use the building blocks $\mCW(G,P\supset Q)$ to state our main theorem, the \emph{automorphic gluing theorem}.

\begin{rem} The letter $\mCW$ in $\mCW(G,P\supset Q)$ is due to its close relation with the extended Whittaker categories defined in \cite{gaitsgory2015outline}. This relation has already been indicated in \cite[Theorem 1.4.8]{beraldo2021geometric} and in the introduction. We do not further explore this topic in this paper.

\end{rem}

\subsection{More functorialities of \texorpdfstring{$\mCI(G,P)$}{I(G,P)}}
\label{ssec-IGP-functorial}
In this subsection, we extend the functorialities of the parabolic categories $\mCI(G,P)$.

\begin{notn} For any parabolic subgroups $Q\subset R$ of $G$, write
\[ \mCI(R,Q):= \mCI(R/U_R, Q/U_R), \]
where $R/U_R$ is the Levi quotient group of $R$ and $Q/U_R$ is the image of $Q$ under the map $R\to R/U_R$.
\end{notn}

For any parabolic subgroups $Q\subset R$ of $G$, consider the diagram
\[ 
 \Bun_G^{Q\hgen} \xleftarrow{p^{Q\hgen}} \Bun_R^{Q\hgen}\xrightarrow{q^{Q\hgen}}  \Bun_{R/U_R}^{Q/U_R\hgen} .
\]
We have:
\begin{lem} \label{lem-Eis!-with-Q-gen}
The map $q^{Q\hgen}$ is relatively represented by quasi-compact smooth algebraic stacks. Also, the functor $(p^{Q\hgen})_!\circ (q^{Q\hgen})^*$ is well-defined on $\mCI(R, Q)$ and sends it into $\mCI(G,Q)$.
\end{lem}

\proof It is easy to see the canonical map $\mBB Q\to \mBB(Q/U_R) \mt_{\mBB (R/U_R)} \mBB R$ is an isomorphism. Hence it induces an isomorphism
\[ \Bun_R^{Q\hgen} \to \Bun_{R/U_R}^{Q/U_R\hgen}\mt_{\Bun_{R/U_R}} \Bun_R  .\]
Then the first claim follows from the fact that $q:\Bun_R\to \Bun_{R/U_R}$ is relatively represented by quasi-compact smooth algebraic stacks.

For the second claim, by Lemma \ref{lem-generator-IGP} (applied to $\mCI(R/U_R,Q/U_R)$), we can pre-compose with the functor $\iota_{M_Q,!}: \DMod(\Bun_{M_Q}) \to \mCI(R,Q)$. Using the base-change isomorphism, this composition is given by $*$-pull-$!$-push along the diagram
\[\Bun_G^{Q\hgen} \gets \Bun_{Q/U_R}\mt_{ \Bun_{R/U_R}^{Q/U_R\hgen} } \Bun_R^{Q\hgen}  \to  \Bun_{M_Q} .\]
Again, the isomorphism $\mBB Q\to \mBB(Q/U_R)\mt_{\mBB (R/U_R)} \mBB R$ induces an isomorphism
\[ \Bun_Q \to \Bun_{Q/U_R}\mt_{ \Bun_{R/U_R}^{Q/U_R\hgen} } \Bun_R^{Q\hgen} . \]
Hence, it remains to show that the $*$-pull-$!$-push along the diagram
\[ \Bun_G^{Q\hgen} \gets \Bun_Q  \to  \Bun_{M_Q}\]
is well-defined and has image in $\mCI(G,Q)$. But this is just Lemma \ref{lem-generator-IGP}.

\qed[Lemma \ref{lem-Eis!-with-Q-gen}]

\begin{defn} For any inclusion of parabolic subgroups $Q\subset R$ of $G$, Lemma \ref{lem-Eis!-with-Q-gen} provides adjoint functors
\[  \Eis_{R,!}^{Q\hgen}: \mCI(R, Q) \adjoint \mCI(G,Q): \CT_{R,*}^{Q\hgen}.\]
Note that the right adjoint $\CT_{R,*}^{Q\hgen}$ is continuous because $\Eis_{R,!}^{Q\hgen}$ preserves compact objects.
\end{defn}

We have the following analogue of Theorem \ref{thm-2nd-adjointness-DG}:

\begin{prop}[{2nd adjointness for $\Bun_G^{Q\hgen}$}] \label{prop-2nd-adjointness-IGQ}
For any \emph{standard} parabolic subgroups $Q\subset P$ of $G$, let $P^-$ be the parabolic subgroup of $G$ opposite to $P$ and $M_P:=P\cap P^-$ be the Levi \emph{subgroup}. Then the functor $\CT_{P,*}^{Q\hgen}$ is canonically equivalent to $*$-pull-$!$-push along the diagram
\[
\Bun_{M_P}^{Q\cap M_P\hgen} \gets \Bun_{P^-}^{Q\cap M_P\hgen} \to \Bun_{G}^{Q\hgen}  .
\]
\end{prop}

\proof We only need to prove that the functor in the statement is right adjoint to $\Eis_{P,!}^{Q\hgen}$. When $Q=P$, this is just \cite[Theorem 1.3.1]{chen2020deligne}. The proof there also works in the general case. Namely, let $\mBG_m\to T$ be a group homomorphism such that the adjoint $\mBG_m$-action on $G$ has attractor, repeller and fixed loci given by $P, P^-$ and $M_P$. The restricted action on $Q$ has attractor, repeller and fixed loci given by $Q,  Q\cap P^-$ and $Q\cap M_P$. Note that we have  $Q\cap P^-=Q\cap M_P$ because $Q\subset P$. Hence, as in \emph{loc.cit.}, we can construct a Drinfeld input (see \cite[\S 3.2-3.3]{chen2020deligne} and \cite[Appendix C]{drinfeld2013algebraic}) with the two correspondences are given by
\[\begin{aligned}
\Bun_{M_P}^{Q\cap M_P\hgen} &\gets& \Bun_{P^-}^{Q\cap M_P\hgen} & \to&  \Bun_{G}^{Q\hgen}  ;\\
 \Bun_{G}^{Q\hgen} &\gets& \Bun_{P}^{Q\hgen} & \to& \Bun_{M_P}^{Q\cap M_P\hgen} .  
\end{aligned}
\]{}
Then we are done by applying Drinfeld's theorem on adjunctions (see \cite[Theorem 3.2.8]{chen2020deligne} and \cite[Appendix C]{drinfeld2013algebraic}).

\qed[Proposition \ref{prop-2nd-adjointness-IGQ}]

\begin{propconstr}\label{propconstr-IGP-functorial}
 Let $\Arr(\Par_G)$ be the category of arrows of parabolic subgroups of $G$. There is a canonical functor
\[
\mCI(-,-)_{\Eis}:\Arr(\Par_G) \to \DGCat,\; [ R\supset Q] \mapsto \mCI(R,Q)
\]
such that
\begin{itemize}
	\item For fixed $Q$, a morphism $ [ R_1\supset Q] \to  [ R_2\supset Q] $ is sent to the functor
	\[ \Eis_{R_1/U_{R_2},!}^{Q/U_{R_2}\hgen}: \mCI(R_1,Q) \to  \mCI(R_2,Q);\]
	\item For fixed $R$, a morphism $ [ R\supset Q_1] \to  [ R\supset Q_2] $ is sent to the functor
	\[ \Eis_{ Q_1/U_R\to Q_2/U_R }^{\enh}: \mCI(R,Q_1) \to  \mCI(R,Q_2) .\]
\end{itemize}
In particular, we obtain a functor
\[
\mCI(-,-)_{\CT}:\Arr(\Par_G)^\op \to \DGCat,\;  [ R\supset Q] \mapsto \mCI(R,Q)
\]
by passing to right adjoints.

\end{propconstr}

\proof By definition, we only need to construct a functor 
\[\Arr(\Par_G) \to \Corr(\PreStk) ,\;  [ R\supset Q]\mapsto \Bun_{R/U_R}^{Q/U_R\hgen} \] 
such that a morphism $ [ R_1\supset Q_1]\to  [ R_2\supset Q_2]$ is sent to the correspondence
\[ \Bun_{R_1/U_{R_1}}^{Q_1/U_{R_1}\hgen}  \gets \Bun_{R_1/U_{R_2}}^{Q_1/U_{R_2}\hgen}   \to \Bun_{R_2/U_{R_2}}^{Q_2/U_{R_2}\hgen} . \]
To do this, it suffices to check that composition of morphisms is well-defined, i.e., for a chain $[ R_1\supset Q_1]\to  [ R_2\supset Q_2] \to  [ R_3\supset Q_3]$, we need the canonical map
\[
\Bun_{R_1/U_{R_3}}^{Q_1/U_{R_3}\hgen} \to \Bun_{R_1/U_{R_2}}^{Q_1/U_{R_2}\hgen}\mt_{  \Bun_{R_2/U_{R_2}}^{Q_2/U_{R_2}\hgen} }  \Bun_{R_2/U_{R_3}}^{Q_2/U_{R_3}\hgen}
\]
to be an isomorphism. But this follows from the canonical isomorphisms
\[ \mBB( R_1/U_{R_3} ) \to \mBB(R_1/U_{R_2})\mt_{  \mBB(R_2/U_{R_2}) } \mBB(R_2/U_{R_3}),\; \mBB( Q_1/U_{R_3} ) \to \mBB(Q_1/U_{R_2})\mt_{  \mBB(Q_2/U_{R_2}) } \mBB(Q_2/U_{R_3}). \]

\qed[Proposition-Construction \ref{propconstr-IGP-functorial}]

\begin{notn} 
For a morphism $[ R_1\supset Q_1]\to  [ R_2\supset Q_2]$ in $\Arr(\Par_G)$, we write 
\[ \Eis_{[ R_1\supset Q_1]\to  [ R_2\supset Q_2]}: \mCI(R_1,Q_1) \to \mCI(R_2,Q_2):\CT_{[ R_2\supset Q_2]\gets  [ R_1\supset Q_1]} \]
 for the adjoint functors provided by Proposition-Construction \ref{propconstr-IGP-functorial}. We omit the subscripts if there is no ambiguity.
\end{notn}

\subsection{Tempered objects in \texorpdfstring{$\mCI(G,Q)$}{I(G,Q)}}
\label{ssec-WGPQ}

In this subsection, we study various temperedness conditions in $\mCI(G,Q)$ and use them to define the building blocks $\mCW(G,P\supset Q)$.

\begin{defn}
For any parabolic subgroup $P$ of $G$, consider the $\Sph_G$-action on $\mCI(G,P)$ given by Proposition \ref{prop-sph-acts-on-IGP} and the identification $\mCO_x\simeq \mCO$. Let $\mCW(G,G\supset P)\subset \mCI(G,P)$ be the full subcategory of $G$-tempered objects, i.e., 
\[ \mCW(G,G\supset P):= \mCI(G,P)^{G\mathrm{-}\temp}.\]
\end{defn}

\begin{defn} For parabolic subgroups $P\subset R$ of $G$, we write
\[ \mCW(R,R\supset P):= \mCW(R/U_R, R/U_R\supset P/U_R).\]
This is a full subcategory of $\mCI(R,P)=\mCI(R/U_R, P/U_R)$. Objects in $\mCI(R,P)$ are said to be \emph{$R$-tempered} if they are contained in $\mCW(R,R\supset P)$.
\end{defn}

\begin{defn} For parabolic subgroups $Q\subset P\subset R$ of $G$, an object $\mCF\in  \mCI(R,Q)$ is said to be \emph{$P$-tempered} if its image under the functor $\CT:  \mCI(R,Q)\to  \mCI(P,Q)$ is $P$-tempered. Let $\mCW(R,P\supset Q)\subset \mCI(R,Q)$ be the full subcategory of $P$-tempered objects, i.e., it is defined by the following Cartesian diagram
\begin{equation} \label{eqn-def-P-temp-IGQ}
  \begin{tikzcd}
   \mCW(R,P\supset Q)
    \arrow[r,"\CT"] \arrow[d,"\subset"]
    \arrow[dr, phantom, "\lrcorner", very near start]
    &\mCW(P, P\supset Q) \arrow[d,"\subset"]  \\
      \mCI(R,Q) \arrow[r,"\CT"]
    & \mCI(P, Q).
  \end{tikzcd}
\end{equation}
We also write $\mCW(R,P):= \mCW(R,P\supset P)$.
\end{defn}

\begin{exam} \label{exam-always-B-temper}
We have $\mCW(R,B\supset B)=\mCI(R,B)$ because of Example \ref{exam-T-temper}.
\end{exam}

\begin{warn} For $Q\subset P\subset R$, there is \emph{no} $\Sph_{M_P}$-action on the category $\mCI(R,Q)$ unless $P=R$, hence we can \emph{not} define the notion of $M_P$-tempered objects in this category using the theory in \S \ref{ssec-temper}. This is the reason we use ``$P$-tempered'' rather than ``$M_P$-tempered'' in the above definition.
\end{warn}

To describe the functorialities of $\mCW(-,-\supset -)$, we need the following category:

\begin{defn} Let $\Tw^{[2]}(\Par_G)$ be the category of twisted chains $[R\supset P\supset Q]$ of parabolic subgroups of $G$, i.e., the objects of $\Tw^{[2]}(\Par_G)$ are standard parabolic subgroups $[R\supset P\supset Q]$, and there is a unique morphism from $[R_1\supset P_1\supset Q_1]$ to $[R_2\supset P_2\supset Q_2]$ iff we have following commutative diagram
\[
\xymatrix{
  R_1  \ar[d]^-\subset & P_1 \ar[l]^-\supset & Q_1 \ar[l]^-\supset \ar[d]^-\subset\\
  R_2 & P_2  \ar[u]^-\subset \ar[l]^-\supset & Q_2 \ar[l]^-\supset.
}
\]
\end{defn}

We prove the following proposition in \S \ref{ssec-IGP-proof-4}.

\begin{thm} \label{thm-functorial-W-Eis}
For any morphism from $[R_1\supset P_1\supset Q_1]$ to $[R_2\supset P_2\supset Q_2]$ in $\Tw^{[2]}(\Par_G)$, the functor
\[\Eis:\mCI(R_1,Q_1)\to \mCI(R_2,Q_2)\]
restricts to a functor
\[ \Eis: \mCW(R_1, P_1\supset Q_1) \to  \mCW(R_2, P_2\supset Q_2). \]
In particular, the functor (see Proposition-Construction \ref{propconstr-IGP-functorial})
\[ \mCI(-,-)_{\Eis}:\Arr(\Par_G) \to \DGCat,\; [ R\supset Q] \mapsto \mCI(R,Q)\]
induces a functor
\[ \mCW(-,-\supset-)_{\Eis}: \Tw^{[2]}(\Par_G)\to \DGCat,\;[ R\supset P\supset Q] \mapsto \mCW(R,P\supset Q). \]

\end{thm}

The following two results logically follow from Theorem \ref{thm-functorial-W-Eis}. But in fact the proof of the theorem uses them. We prove them in \S \ref{ssec-IGP-proof-2} and \S \ref{ssec-IGP-proof-3}.

\begin{prop}[Case {$R_1=P_1=P_2=P, Q_1=Q_2=Q, R_2=R$}] \label{prop-left-adjointable-of-def-P-temp-IGQ}
For parabolic subgroups $Q\subset P\subset R$ of $G$, the functor $\Eis: \mCI(P, Q)\to  \mCI(R,Q) $ preserves $P$-tempered objects. In other words, it induces a functor
\[ \Eis:\mCW(P, P\supset Q)\to \mCW(R,P\supset Q).  \]
\end{prop}

\begin{rem}
Because of the lack of $\Sph_{M_P}$-action on $\mCI(R,Q)$, the above proposition is \emph{non-obvious}.
\end{rem}

\begin{prop}[Case {$R_1=R_2=R, Q_1=Q_2=Q$}] \label{prop-G-temp-imply-M-temp-global}
For parabolic subgroups $Q\subset P_1\subset P_2\subset R$ of $G$, any $P_2$-tempered object in $\mCI(R,Q)$ is $P_1$-tempered, i.e., $\mCW(R, P_2\supset Q) \subset \mCW(R,P_1\supset Q)$.

\end{prop}

\begin{rem} 
When $Q=P_1, P_2=R=G$, the above proposition can be viewed as a global analogue of Corollary \ref{cor-G-temp-imply-M-temp-local}.
\end{rem}

The proof of the following corollary is independent of Theorem \ref{thm-functorial-W-Eis}. In fact, the proof of the theorem uses it.

\begin{cor} \label{cor-cpt-gen-W}
For parabolic subgroups $Q\subset P\subset R$ of $G$, the category $\mCW(R,P\supset Q)$ is compactly generated.
\end{cor}

\proof The functor $\Eis:\mCW(P, P\supset Q)\to \mCW(R,P\supset Q)$ in Proposition \ref{prop-left-adjointable-of-def-P-temp-IGQ} preserves compact objects and its image generates the target under colimits because its right adjoint is continuous and conservative. Hence, we just need to show $\mCW(P, P\supset Q)$ is compactly generated. Replacing $G$ by $P/U_P$, we can assume $P=G$. In other words, it remains to show that $\mCI(G,Q)^{G\mathrm{-temp}}$ is compactly generated: this follows from Lemma \ref{lem-generator-IGP} and Lemma \ref{lem-cpt-gen-temp-general}.

\qed[Corollary \ref{cor-cpt-gen-W}]

\begin{cor} \label{cor-P-temperization}
 For parabolic subgroups $Q\subset P\subset R$ of $G$, the inclusion functor $ \mCW(R,P\supset Q)\to \mCI(R,Q)$ has a continuous right adjoint
\[ \temp_P: \mCI(R,Q)\to \mCW(R,P\supset Q). \]
Furthermore, the commutative diagram (\ref{eqn-def-P-temp-IGQ}) induces a commutative diagram
\[
\begin{tikzcd}
   \mCW(R,P\supset Q)
    \arrow[r,"\CT"] 
    &\mCW(P, P\supset Q)   \\
      \mCI(R,Q) \arrow[r,"\CT"] \arrow[u,"\temp_P"]
    & \mCI(P, Q). \arrow[u,"\temp_P"]
  \end{tikzcd}
  \]
\end{cor}

\proof To prove the first claim, we show that $ \mCW(R,P\supset Q)\to \mCI(R,Q)$ preserves compact objects. Using Proposition \ref{prop-left-adjointable-of-def-P-temp-IGQ}, we can reduce to the case $R=P$. Replacing $G$ by $P/U_P$, we can assume $G=R=P$. Then Lemma \ref{lem-cpt-gen-temp-general} implies $ \mCW(G,G\supset Q)\to \mCI(G,Q)$ preserves compact objects.

The second claim follows from Proposition \ref{prop-left-adjointable-of-def-P-temp-IGQ} by passing to right adjoints.

\qed[Corollary \ref{cor-P-temperization}]

We will also need the right adjoint version of Theorem \ref{thm-functorial-W-Eis}. To state it, we will rely on the following lemma.

\begin{lem} \label{lem-CT-preserve-P-temper}
For parabolic subgroups $Q_1\subset Q_2\subset P\subset R_1\subset R_2$, the functor $\CT:\mCI(R_2,Q_2)\to \mCI(R_1,Q_1)$ preserves $P$-tempered objects. In other words, it restricts to a functor
\[  \CT:\mCW(R_2,P\supset Q_2)\to \mCW(R_1,P\supset Q_1).\]
\end{lem}

\proof By definition, it suffices to verify that the functor $\CT:\mCI(R_2,Q_2)\to \mCI(P,Q_1)$ preserves $P$-tempered objects. This functor factors as $\mCI(R_2,Q_2)\to \mCI(P,Q_2) \to \mCI(P,Q_1)$. The first functor preserves $P$-tempered objects by definition. The second functor preserves $P$-tempered objects because $P$-temperedness for the source and target is defined as $M_P$-temperedness, and this functor is $\Sph_{M_P,x}$-linear.

\qed[Lemma \ref{lem-CT-preserve-P-temper}]

\begin{cor} \label{cor-functorial-W-CT}
For any morphism from $[R_1\supset P_1\supset Q_1]$ to $[R_2\supset P_2\supset Q_2]$ in $\Tw^{[2]}(\Par_G)$, the functor
\[ \Eis: \mCW(R_1, P_1\supset Q_1) \to  \mCW(R_2, P_2\supset Q_2) \]
in Theorem \ref{thm-functorial-W-Eis} has a continuous right adjoint
\[ \CT_{\temp}:  \mCW(R_2, P_2\supset Q_2)\to  \mCW(R_1, P_1\supset Q_1)  \]
 equivalent to both
\[  \mCW(R_2, P_2\supset Q_2) \xrightarrow{\CT} \mCW(R_1, P_2\supset Q_1) \xrightarrow{\temp_{P_1}} \mCW(R_1, P_1\supset Q_1)  \]
and 
\[  \mCW(R_2, P_2\supset Q_2) \xrightarrow{\temp_{P_1}} \mCW(R_2, P_1\supset Q_2) \xrightarrow{\CT} \mCW(R_1, P_1\supset Q_1) . \]
In particular, the functor
\[ \mCW(-,-\supset-)_{\Eis}: \Tw^{[2]}(\Par_G)\to \DGCat,\;[ R\supset P\supset Q] \mapsto \mCW(R,P\supset Q)\]
induces a functor
\[ \mCW(-,-\supset-)_{\CT_{\temp}}: \Tw^{[2]}(\Par_G)^\op\to \DGCat,\;[ R\supset P\supset Q] \mapsto \mCW(R,P\supset Q) \]
by passing to right adjoints.
\end{cor}

\proof By Theorem \ref{thm-functorial-W-Eis}, the functor $ \Eis: \mCW(R_1, P_1\supset Q_1) \to  \mCW(R_2, P_2\supset Q_2) $ is equivalent to both
\[ \mCW(R_1, P_1\supset Q_1) \subset  \mCW(R_1, P_2\supset Q_1) \xrightarrow{\Eis}  \mCW(R_2, P_2\supset Q_2) \]
and
\[\mCW(R_1, P_1\supset Q_1)  \xrightarrow{\Eis} \mCW(R_2, P_1\supset Q_2) \subset  \mCW(R_2, P_2\supset Q_2) . \]
Then we are done by Corollary \ref{cor-P-temperization} and Lemma \ref{lem-CT-preserve-P-temper}.

\qed[Corollary \ref{cor-functorial-W-CT}]

The following result is proved in \S \ref{ssec-IGP-proof-5}:

\begin{prop} \label{prop-W-preserve-by-Sph}
For parabolic subgroups $Q\subset P\subset R$ of $G$, let $M_R:=R/U_R$ be the Levi quotient group of $R$. Then the $\Sph_{M_R,x}$-action on $\DMod(\Bun_{M_R}^{Q/U_R\hgen})$ preserves the full subcategories
\[ \mCW(R,P\supset Q) \subset \mCI(R,Q) \subset  \DMod(\Bun_{M_R}^{Q/U_R\hgen}).\]
\end{prop}

\subsection{Proof of Proposition \ref{prop-left-adjointable-of-def-P-temp-IGQ}}
\label{ssec-IGP-proof-2}

We can assume $R=G$ and $Q\subset P$ are standard parabolics. Define $P^-$ and $M:=P\cap P^-$ as usual. 

By definition, we just need to check that the endo-functor $ \mCI(P, Q)\xrightarrow{\Eis}  \mCI(G,Q) \xrightarrow{\CT}\mCI(P,Q)$ preserves $P$-tempered objects. Using 2nd adjointness on $\Bun_G^{Q\hgen}$ (Proposition \ref{prop-2nd-adjointness-IGQ}), the functor $\CT$ is $*$-pull-$!$-push along the diagram
\[ \Bun_{M}^{Q\cap M\hgen}\gets  \Bun_{P^-}^{Q\cap M\hgen} \to \Bun_G^{Q\hgen} .\]
On the other hand, the functor $\Eis$ is $*$-pull-$!$-push along the diagramm
\[\Bun_G^{Q\hgen} \gets  \Bun_{P}^{Q\hgen} \to \Bun_{M}^{Q\cap M\hgen} \]
Hence, we have reduced the proposition to the following result:

\begin{lem} \label{lem-Zastav-preserve-temper}
Using the above notations, the $*$-pull-$!$-push functor along the diagram
\[ \Bun_{M}^{Q\cap M\hgen}\gets  \Bun_{P^-}^{Q\cap M\hgen} \mt_{\Bun_G^{Q\hgen} } \Bun_P^{Q\hgen} \to \Bun_{M}^{Q\cap M\hgen}\]
preserves $M$-tempered objects.

\end{lem}

To prove the lemma, we need to study the fiber product
\[\Bun_{P^-}^{Q\cap M\hgen} \mt_{\Bun_G^{Q\hgen} } \Bun_P^{Q\hgen}\simeq \mathbf{Maps}_{\mathrm{gen}}(X,P^-\backslash G/P \gets Q\cap M\backslash Q/Q )\]
and review some geometric objects related to it.

\begin{defn} Recall the \emph{relative open Zastava space} defined in\footnote{Using the notation of \cite[\S 3]{braverman2002intersection}, it is the disjoint union of $_0\!Z_{\Bun_M}^\theta$ for $\theta\in \Lambda_{G,P}^{\pos}$.} \cite{braverman2002intersection}, \cite{feigin1999semiinfinite} (also see \cite[\S 6.3]{schieder2016geometric}:
\[  _0\!\Zas_{P,\rel}:= \mathbf{Maps}_{\mathrm{gen}}(X,P^-\backslash G/P \gets \mBB M ). \]
In more familiar words, it classifies $(\mCP_G, \mCP_M^l,\mCP_M^r,\alpha_V,\beta_V)$ where:
\begin{itemize}
  \item $\mCP_G$ is a $G$-torsor on $X$;
  \item $\mCP_M^r$ and $\mCP_M^l$ are $M$-torsors on $X$;
  \item For any finite dimensional $G$-representation $V$ and the corresponding $M$-respresentation $W:= V^{U_P}$, $\alpha_V: W_{\mCP_{M^r}} \to V_{\mCP_G} $ is an injection of vector bundles on $X$, and $\beta_V: V_{\mCP_G}\to W_{\mCP_{M^r}}$ is an surjection of vector bundles;
  \item The collections $\{\alpha_V\}$ and $\{\beta_V\}$ satisfy the Drinfeld-Plucker relations (which we do not spell out);
  \item The composition $\beta_V\circ \alpha_V:  W_{\mCP_{M^r}}\to W_{\mCP_{M^l}} $ is required to be an \emph{injection between coherent sheaves}.
\end{itemize}
Its connected components are all smooth. Also recall the absolute open Zastava space $_0\!\Zas_{P}$: it is defined similarly by replacing $\mCP_M^r$ with the trivial $M$-torsor.  We have
\[   _0\!\Zas_{P,\rel}\simeq \, _0\!\Zas_{P} \widetilde{\mt}\Bun_M. \]

We define
\[ _0\!\Zas_{P,\rel}^{Q\cap M\hgen}:= \mathbf{Maps}_{\mathrm{gen}}(X,P^-\backslash G/P \gets Q\cap M\backslash Q/Q ).
\]
In more familiar words, it classifies $(\mCP_G, \mCP_M^1,\mCP_M^2,\alpha_V,\beta_V,\mCP_{Q\cap M}^{\mathrm{gen}})$, where:
\begin{itemize}
  \item $(\mCP_G, \mCP_M^l,\mCP_M^r,\alpha_V,\beta_V)$ is as in the definition of $_0\!\Zas_{P,\rel}$;
  \item $\mCP_{Q\cap M}^{\mathrm{gen}}$ is a generic $(Q\cap M)$-reduction of both $\mCP_{M^r}$ and $\mCP_{M^l}$, compatible with the maps $\beta_V\circ \alpha_V$.
\end{itemize}
More precisely, the maps $\beta_V\circ \alpha_V$ are isomorphisms at the generic point of $X$, hence they provide an identification of the generic $M$-torsors of $\mCP_{M^r}$ and $\mCP_{M^l}$; then we require $\mCP_{Q\cap M}^{\mathrm{gen}}$ to be compatible with this identification.

We also define the absolute version $_0\!\Zas_{P}$ by replacing $\mCP_M^r$ in the above definition by the trivial $M$-torsor. We have
\[   _0\!\Zas_{P,\rel}^{Q\cap M\hgen}\simeq \, _0\!\Zas_{P} \widetilde{\mt}\Bun_M^{Q\cap M\hgen}. \]
\end{defn}

\begin{defn} Recall the \emph{positive part of the Hecke ind-stack for $M$-torsors}, which we denote\footnote{Using the notation of \cite[\S 3]{braverman2002intersection}, it is the disjoint union of $\mathrm{Mod}_{M}^{+,\theta}$ for $\theta\in \Lambda_{G,P}^{\pos}$.} by $\Hecke_{M,G\hyphen\pos}$. It classifies $(\mCP_M^l,\mCP_M^r,\gamma_V)$ where
\begin{itemize}
  \item $\mCP_M^r$ and $\mCP_M^l$ are $M$-torsors on $X$;
  \item For any finite dimensional $G$-representation $V$ and the corresponding $M$-representation $W:= V^{U_P}$, the map $\gamma_V: W_{\mCP_{M^r}}\to W_{\mCP_{M^l}}  $ is an \emph{injection between coherent sheaves};
  \item The collection of $\{\gamma_V\}$ should satisfy certain relations which we do not spell out.
\end{itemize}
Also recall $\Gr_{M,G\hyphen\pos}$ is defined similarly by replacing $\mCP_M^r$ by the trivial $M$-torsor. We have
\[ \Hecke_{M,G\hyphen\pos}\simeq \Gr_{M,G\hyphen\pos}\widetilde{\mt}\Bun_M. \]

There is a map
\[  _0\!\Zas_{P,\rel}\to \Hecke_{M,G\hyphen\pos} \]
sending $(\mCP_G, \mCP_M^l,\mCP_M^r,\alpha_V,\beta_V)$ to $( \mCP_M^l,\mCP_M^r,\beta_V\circ \alpha_V)$.

We can similarly define $\Hecke_{M,G\hyphen\pos}^{Q\cap M\hgen}$ and the map
\[  _0\!\Zas_{P,\rel}^{Q\cap M\hgen}\to \Hecke_{M,G\hyphen\pos}^{Q\cap M\hgen}. \]
We also have
\[ \Hecke_{M,G\hyphen\pos}^{Q\cap M\hgen} \simeq \Gr_{M,G\hyphen\pos} \widetilde{\mt}\Bun_M^{Q\cap M\hgen}.\]
\end{defn}

\begin{constr} Recall there is a map
\begin{equation} \label{eqn-map-pos-Hecke-to-divisor}  \Hecke_{M,G\hyphen\pos} \to  \Div_{\Lambda_{G,P}^{\pos}}:= \bigsqcup_{\theta\in \Lambda_{G,P}^{\pos}} X^\theta, 
\end{equation}
where
\begin{itemize}
  \item $\Lambda_G^{\pos}\subset \Lambda_G$ is the submonoid spanned by all positive simple coroots in the coweight lattice $\Lambda_G$;
  \item $\Lambda_{G,P}$ is the quotient of $\Lambda_G$ by the $\mBZ$-span of simple coroots contained in $M$;
  \item $\Lambda_{G,P}^{\pos}$ is the image of $\Lambda_G^{\pos}$ inside $\Lambda_{G,P}$;
  \item Each $\theta\in \Lambda_{G,P}^{\pos}$ can be uniquely written as the image of $\sum n_i\alpha_i$ for $n_i\in \mBN^{\ge 0}$ and simple coroots $\alpha_i$ \emph{not} contained in $M$; we write $X^\theta:= \prod X^{(n_i)}$, where $X^{(n_i)}$ is the $n_i$-th symmetric product of $X$;
  \item Any point of $ \Div_{\Lambda_{G,P}^{\pos}}$ can be uniquely written as $\sum \theta_k x_k$ with $\theta_k\in \Lambda_{G,P}^{\pos}$ and $x_k$ being distinct points on $X$, i.e., is a $\Lambda_{G,P}^{\pos}$-colored divisor on $X$.
  \item For $(\mCP_M^l,\mCP_M^r,\gamma_V)$ as in the definition of $\Hecke_{M,G\hyphen\pos}$, the map (\ref{eqn-map-pos-Hecke-to-divisor}) sends it to the unique point $\sum \theta_k x_k$ such that when $W=V^{U_P}$ is the 1-dimensional representation given by character $\check \lambda$, the map $\gamma_V: W_{\mCP_{M^r}}\to W_{\mCP_{M^l}}$ induces an isomorphism $W_{\mCP_{M^r}}( \sum \langle \theta_k, \check \lambda \rangle\cdot x_k )\simeq W_{\mCP_{M^l}}$.
\end{itemize}
The similarly defined map 
\[  \Gr_{M,G\hyphen\pos} \to \Div_{\Lambda_{G,P}^{\pos}},
\]
satisfies the following \emph{factorization property}. Let
\[ ( \Div_{\Lambda_{G,P}^{\pos}} \mt \Div_{\Lambda_{G,P}^{\pos}})_\disj \subset \Div_{\Lambda_{G,P}^{\pos}}\mt \Div_{\Lambda_{G,P}^{\pos}}\]
be the open subscheme containing points $(\sum \theta_k x_k,\sum \theta'_l x'_l)$ such that $x_k,x'_l$ are distinct. Let
\[ (\Gr_{M,G\hyphen\pos} \mt \Gr_{M,G\hyphen\pos} )_\disj \subset \Gr_{M,G\hyphen\pos} \mt \Gr_{M,G\hyphen\pos}  \]
be the pre-image of this open subscheme. Consider the map
\[\mathrm{sum}: ( \Div_{\Lambda_{G,P}^{\pos}} \mt \Div_{\Lambda_{G,P}^{\pos}})_\disj \to  \Div_{\Lambda_{G,P}^{\pos}}  \]
sending $(\sum \theta_k x_k,\sum \theta'_l x'_l)$ to $\sum \theta_k x_k + \sum \theta'_l x'_l$. Then there is a canonical Cartesian square
\[
\begin{tikzcd}
     (\Gr_{M,G\hyphen\pos} \mt \Gr_{M,G\hyphen\pos} )_\disj
    \arrow[r] \arrow[d]
    \arrow[dr, phantom, "\lrcorner", very near start]
    &\Gr_{M,G\hyphen\pos} \arrow[d]  \\
      ( \Div_{\Lambda_{G,P}^{\pos}} \mt \Div_{\Lambda_{G,P}^{\pos}})_\disj \arrow[r,"\mathrm{sum}"]
    &  \Div_{\Lambda_{G,P}^{\pos}},
  \end{tikzcd}
\]
where the top arrow is defined using a standard re-gluing construction.

The composition 
\[ _0\Zas_{P}\to  \Gr_{M,G\hyphen\pos} \to \Div_{\Lambda_{G,P}^{\pos}} \]
also satisfies factorization property.
\end{constr}

With these geometric constructions in our toolbox, we are ready to prove Lemma \ref{lem-Zastav-preserve-temper} and therefore finish the proof of Proposition \ref{prop-left-adjointable-of-def-P-temp-IGQ}.

\medskip

\noindent\emph{Proof of Lemma \ref{lem-Zastav-preserve-temper}.} Consider the commutative diagram
\[
\xymatrix{
  & _0\!\Zas_{P,\rel}^{Q\cap M\hgen} \ar[d]^-{\pi} \ar[ld]^-{p_l} \ar[rd]^-{p_r} \\
  \Bun_{M}^{Q\cap M\hgen} & \Hecke_{M,G\hyphen\pos}^{Q\cap M\hgen} \ar[l]^-{h_l} \ar[r]_-{h_r} &  \Bun_{M}^{Q\cap M\hgen} .
}
\]
We just need to check that $p_{l,!}\circ p_r^*$ preserves $M$-tempered objects. Using the projection formula, this functor is equivalent to $h_{l,!}( \mCK \ot^* h_r^*(-) )$, where $\mCK$ is the $!$-pushforward of the constant sheaf along $\pi$. In other words, it is the functor given by kernel $ \mCK$.

Recall we have a map
\[  \Hecke_{M,G\hyphen\pos}^{Q\cap M\hgen}\to  \bigsqcup_{\theta\in \Lambda_{G,P}^{\pos}} X^\theta.\]
Let $\mCK|_{X^\theta}$ be the restriction of $\mCK$ on the pre-image of $X^\theta$. Then it remains to prove that the functor given by kernel $\mCK|_{X^\theta}$ preserves $M$-tempered objects.

Note that $X^\theta$ has a stratification given by 
\[\bigcup_{\mu\in 0\cup \Lambda_{G,P}^\pos, \mu\le \theta} (X\setminus x)^{\theta-\mu}\mt x^{\mu}.\] 
This induces a stratification on $\Hecke_{M,G\hyphen\pos}^{Q\cap M\hgen}|_{X^\theta}$. Let $\mCK|_{(X\setminus x)^{\theta-\mu}\mt x^{\mu}}$ be the $*$-restriction of $\mCK|_{X^\theta}$ to this stratum. Then we just need to prove that the functor given by kernel $\mCK|_{(X\setminus x)^{\theta-\mu}\mt x^{\mu}}$ preserves $M$-tempered objects.

By the factorization properties, we have
\[ \Hecke_{M,G\hyphen\pos}^{Q\cap M\hgen,\theta,\mu\cdot x } \simeq 
\Hecke_{M,G\hyphen\pos}^{Q\cap M\hgen }|_{(X\setminus x)^{\theta-\mu}} \mt_{\Bun_M^{Q\cap M\hgen}} \Hecke_{M,G\hyphen\pos}^{Q\cap M\hgen }|_{\mu \cdot x},
   \]
and $\mCK|_{(X\setminus x)^{\theta-\mu}\mt x^{\mu}}$ is the $*$-tensor product of $*$-pullbacks of $\mCK|_{(X\setminus x)^{\theta-\mu}}$ and $\mCK|_{\mu\cdot x}$. We will prove that the functors given by kernels $\mCK|_{(X\setminus x)^{\theta-\mu}}$ and $\mCK|_{\mu\cdot x}$ both preserve $M$-tempered objects.

The claim for the first functor is obvious: Hecke modifications \emph{away} from $x$ commute with Hecke modifications at $x$, hence the functor given by kernel $\mCK|_{(X\setminus x)^{\theta-\mu}}$ commutes with the $\Sph_{M,x}$-action on $\DMod(\Bun_M^{Q\cap M\hgen})$.

For the claim for the second functor, note that we have a map
\[ \Hecke_{M,G\hyphen\pos}^{Q\cap M\hgen }|_{\mu \cdot x} \to \Hecke_{M,x}^{Q\cap M\hgen }.\]
Hence the second functor is equivalent to $\mCL\convolve_{\Sph_{M,x}}-$, where $\mCL$ is the $!$-pushforward of $\mCK|_{\mu\cdot x}$ along this map. Then we are done because $\Sph_{M,x}^{\temp}$ is a two-sided ideal of $\Sph_{M,x}$.

\qed[Lemma \ref{lem-Zastav-preserve-temper}]

\qed[Proposition \ref{prop-left-adjointable-of-def-P-temp-IGQ}]
\subsection{Proof of Proposition \ref{prop-G-temp-imply-M-temp-global}}
\label{ssec-IGP-proof-3}

We will deduce the proposition from the following result:

\begin{prop} \label{prop-CT-preserve-temperedness-general} For parabolic subgroups $Q\subset P\subset R$ of $G$, the functor $\CT: \mCI(R,Q) \to \mCI(P, Q)$
sends $R$-tempered objects to $P$-tempered objects. In other words, it restricts to a functor
\[\CT: \mCW(R,R\supset Q) \to \mCW(P, P\supset Q).\]
\end{prop}

\proof We can assume $R=G$. The proof is similar to that of Theorem \ref{thm-iota!*-temperedness} in \S \ref{ssec-SI-to-IGP-temper} and we only sketch the main ideas. We just need to construct $\Sph_{G,x}$-linear adjoint functors
\[ \SI_{P,x}\ot_{\Sph_{M_P,x}} \mCI(P, Q) \adjoint \mCI(G,Q) \]
such that the composition
\[ \mCI(P, Q)\adjoint \SI_{P,x}\ot_{\Sph_{M_P,x}} \mCI(P, Q) \adjoint \mCI(G,Q) \]
is equivalent to the adjoint pair $(\Eis,\CT)$. As before, we have a canonical map
\[ \pi: \Gr_{G,x} \twisttimes \Bun_{M_P}^{Q_P\hgen} \to \Bun_G^{Q\hgen} \]
and the desired left adjoint is defined by sending $\mCF\boxt_{\Sph_{M_P,x}} \mCM$ to $\pi_!(\mCF \widetilde\boxtimes \mCM)$. 

\qed[Proposition \ref{prop-CT-preserve-temperedness-general}]

Let us finish the proof of Proposition \ref{prop-G-temp-imply-M-temp-global}. Unraveling the definitions, we need to show that if an object in $\mCI(R,Q)$ is sent to a $P_2$-tempered object by the functor $\CT:\mCI(R,Q) \to \mCI(P_2, Q)$, then it is sent to a $P_1$-tempered object by the functor $\CT:\mCI(R,Q) \to \mCI(P_1, Q)$. Note that the constant term functor $\mCI(R,Q) \to \mCI(P_1, Q)$ factors as the following composition of two constant term functors: $\mCI(R,Q) \to \mCI(P_2, Q) \to \mCI(P_1,Q)$. Hence the above claim follows from Proposition \ref{prop-CT-preserve-temperedness-general} (applied to the second functor in the above composition).

\qed[Proposition \ref{prop-G-temp-imply-M-temp-global}]
\subsection{Proof of Theorem \ref{thm-functorial-W-Eis}}
\label{ssec-IGP-proof-4}

Any morphism $[R_1\supset P_1\supset Q_1]\to [R_2\supset P_2\supset Q_2]$ in $\Tw^{[2]}(\Par_G)$ factors as
\[ [R_1\supset P_1\supset Q_1] \to [R_1\supset P_1\supset Q_2] \to [R_1\supset P_2\supset Q_2] \to [R_2\supset P_2\supset Q_2].\]
Hence we just need to prove Theorem \ref{thm-functorial-W-Eis} in the following three special cases: (1) when $R_1=R_2=R, P_1=P_2=P$; (2) when $R_1=R_2=R, Q_1=Q_2=Q$; (3) when $P_1=P_2=P, Q_1=Q_2=Q$.

For case (1), by definition we need to show that if an object $\mCF$ in $\mCI(R,Q_1)$ is sent to a $P$-tempered object in $\mCI(P,Q_1)$ by the functor $\CT:\mCI(R,Q_1)\to \mCI(P,Q_1)$, then its image $\Eis(\mCF)$ in $\mCI(R,Q_2)$ is sent to a $P$-tempered object in $\mCI(P,Q_2)$ by the functor $\CT:\mCI(R,Q_2)\to \mCI(P,Q_2)$. By Lemma \ref{lem-Eis-CT-adjointable} below, it remains to check that the functor $\Eis:  \mCI(P,Q_1)\to \mCI(P,Q_2)$ preserves $P$-tempered objects. But this is obvious because $P$-temperedness for these categories is defined as $M_P$-temperedness, and the functor in question is $\Sph_{M_P,x}$-linear.

Case (2) is just Proposition \ref{prop-G-temp-imply-M-temp-global}.

For case (3), we use the fact that $\mCW(R_1,P\supset Q)$ is generated under colimits by the functor $\Eis: \mCW(P,P\supset Q)\to\mCW(R_1,P\supset Q)  $, see the proof of Corollary \ref{cor-cpt-gen-W}. Hence it suffices to prove that functor $\Eis:\mCI(P,Q)\to \mCI(R_2,Q)$ preserves $P$-tempered objects. But this is just Proposition \ref{prop-left-adjointable-of-def-P-temp-IGQ}.

\qed[Theorem \ref{thm-functorial-W-Eis}]

The following ``base-change'' lemma was used in the above proof and will be used again later:

\begin{lem} \label{lem-Eis-CT-adjointable}
For any parabolic subgroups $Q_1\subset Q_2\subset P\subset R$ of $G$, the commutative square
\[
\xymatrix{
	 \mCI(P, Q_1)  \ar[d]^-{\Eis}
	\ar[r]^-\Eis &
	 \mCI(R, Q_1) \ar[d]^-{\Eis}
	\\
	 \mCI(P, Q_2)	 \ar[r]^-\Eis&
	\mCI(R,Q_2)
}
\]
is right adjointable along the horizontal direction, i.e., we have a canonical commutative square
\[
\xymatrix{
	 \mCI(P, Q_1)  \ar[d]^-{\Eis}
	 &
	 \mCI(R, Q_1) \ar[d]^-{\Eis} \ar[l]^-\CT
	\\
	 \mCI(P, Q_2)	&
	\mCI(R,Q_2) \ar[l]^-\CT.
}
\]
\end{lem}

\proof We can assume $R=G$ and the parabolic subgroups are standard. Define $P^-$ and $M_P$ as usual. By 2nd adjointness for $\Bun_G^{Q_i\hgen}$ (Proposition \ref{prop-2nd-adjointness-IGQ}), the assertion boils down to proving that the composition of the correspondences
\[ \Bun_{M_P}^{Q_1\cap M_P\hgen} \gets \Bun_{P^-}^{Q_1\cap M_P\hgen} \to \Bun_G^{Q_1\hgen}  \]
and
\[ \Bun_{M_P}^{Q_2\cap M_P\hgen}\gets \Bun_{M_P}^{Q_1\cap M_P\hgen} \to \Bun_{M_P}^{Q_1\cap M_P\hgen}  \]
is isomorphic to the composition of the correspondences
\[ \Bun_{G}^{Q_2\hgen}\gets \Bun_{G}^{Q_1\hgen} \to \Bun_{G}^{Q_1\hgen}  \]
and
\[ \Bun_{M_P}^{Q_2\cap M_P\hgen} \gets \Bun_{P^-}^{Q_2\cap M_P\hgen} \to \Bun_G^{Q_2\hgen} . \]
The first composition is
\[ \Bun_{M_P}^{Q_2\cap M_P\hgen} \gets \Bun_{P^-}^{Q_1\cap M_P\hgen} \to \Bun_G^{Q_1\hgen} . \]
Hence, it remains to prove that the map
\[ \Bun_{P^-}^{Q_1\cap M_P\hgen}\to \Bun_{G}^{Q_1\hgen}\mt_{\Bun_{G}^{Q_2\hgen}} \Bun_{P^-}^{Q_2\cap M_P\hgen}\]
is an isomorphism. This follows from the isomorphism $\mBB(Q_1\cap M_P)\to \mBB Q_1\mt_{\mBB Q_2} \mBB(Q_2\cap M_P) $.

\qed[Lemma \ref{lem-Eis-CT-adjointable}]
\subsection{Proof of Proposition \ref{prop-W-preserve-by-Sph}}
\label{ssec-IGP-proof-5}

We can assume $R=G$. In view of Proposition \ref{prop-sph-acts-on-IGP}, the full subcategory $\mCI(G,Q) \subset  \DMod(\Bun_{G}^{Q\hgen})$ is preserved by the $\Sph_{G,x}$-action. To prove the claim about $\mCW(G,P\supset Q)$, we show that, for any $\mCM\in \mCI(G, Q)$ and $\mCF\in \Sph_{G,x}$, if the image of $\mCM$ under the functor $\CT: \mCI(G,Q)\to \mCI(P,Q)$ is $M_P$-tempered, then so is the image of $\mCF\convolve_{\Sph_{G,x}}\mCM$. Recall in the proof of Proposition \ref{prop-CT-preserve-temperedness-general}, we have factorized this functor $\CT$ as
\[ \mCI(G,Q) \xrightarrow{u} \SI_{P,x}\ot_{\Sph_{M_P,x}} \mCI(P, Q)  \xrightarrow{v} \mCI(P, Q). \]
By construction, the left adjoint of the first functor $u$ is $\Sph_{G,x}$-linear, hence so is $u$ (see Appendix \ref{sect-lax-linear}). Thus, it remains to show that, for any $\mCN\in \SI_{P,x}\ot_{\Sph_{M_P,x}} \mCI(P, Q) $ and $\mCF\in \Sph_{G,x}$, if $v(\mCN)$ is $M_P$-tempered, then so is $v( \mCF\convolve_{\Sph_{G,x}}\mCN )$.

The proof below works for any $\Sph_{M_P,x}$-module category $\mCC$ and in particular for $\mCC:=\mCI(P,Q)$. Consider the obvious commutative diagram
\[
\xymatrix{
	\SI_{P,x}\ot_{\Sph_{M_P,x}} \mCC \ar[r]^-{\iota_{P,x}^!\ot\Id} &
	\mCC \\
	\SI_{P,x}\ot_{\Sph_{M_P,x}} \mCC^{M_P\hyphen\temp} \ar[r]^-{\iota_{P,x}^!\ot \Id} \ar[u]^-{\Id\ot \oblv_{M_P}} &
	\mCC^{M_P\hyphen\temp}. \ar[u]^-{ \oblv_{M_P}}
}
\]
This square is left adjointable along the horizontal direction because $\iota_{P,x}^!$ has a $\Sph_{M_P,x}$-linear left adjoint $(\iota_{P,x})_!$. It follows that this square is also right adjointble along the vertical direction. In other words, we have a commutative square
\[
\xymatrix{
	\SI_{P,x}\ot_{\Sph_{M_P,x}} \mCC \ar[r]^-{\iota_{P,x}^!\ot\Id} \ar[d]^-{\Id\ot \temp_{M_P}}
	&
	\mCC \ar[d]^-{\temp_{M_P}} \\
	\SI_{P,x}\ot_{\Sph_{M_P,x}} \mCC^{M_P\hyphen\temp} \ar[r]^-{\iota_{P,x}^!\ot \Id}  &
	\mCC^{M_P\hyphen\temp}. 
}
\]
Now, for any object $\mCN\in \SI_{P,x}\ot_{\Sph_{M_P,x}} \mCC $, the object $(\iota_{P,x}^!\ot\Id)(\mCN)$ is $M_P$-tempered iff the monad $\oblv_{M_P}\circ \temp_{M_P}$ acts on it trivially, iff the morphism
\[ (\Id\ot \oblv_{M_P}\circ \temp_{M_P})(\mCN)\to \mCN \]
is sent to an isomorphism by $\iota_{P,x}^!\ot\Id$. Since the last functor is conservative (because its left adjoint generate the target under colimits), we see $(\iota_{P,x}^!\ot\Id)(\mCN)$ is $M_P$-tempered iff $(\Id\ot \oblv_{M_P}\circ \temp_{M_P})(\mCN)\to \mCN $ is an isomorphism. The latter condition is clearly preserved by the $\Sph_{G,x}$-action on $\mCN$, hence we are done.

\qed[Proposition \ref{prop-W-preserve-by-Sph}]
\section{The automorphic gluing theorem} \label{sect-aut-gluing-thm}

In this section, we state our main theorem and perform several reduction steps.

In \S \ref{ssec-mainthm}, we give the precise formulation of the main theorem (Theorem \ref{thm-main-thm}), which states that $\DMod(\Bun_G)$ can be glued from $\mCW(G,P\supset Q)$, where $[P\supset Q]$ ranges over \emph{twisted arrows} in the category of standard parabolic subgroups of $G$.

We also give an anti-tempered version of the main theorem (Theorem \ref{thm-main-thm-atemp}), which says that the category $\DMod(\Bun_G)^{G\hyphen\atemp}$ can be glued from $\mCI(G,R)^{G\hyphen\atemp}$, where $R$ ranges over \emph{proper} standard parabolic subgroups of $G$.

We will reduce the main theorem to its anti-tempered version, and prove the latter in \S \ref{sect-gluing-atemp-proof}.

In \S \ref{ssec-gluing-IGR}, to achieve this reduction, we state the automorphic gluing for the parabolic category, which says that $\mCI(G,R)$ can be glued from $\mCW(G,P \supset Q)$, where $[P \supset Q]$ ranges over twisted arrows in the category of standard parabolic subgroups of $G$ contained in $R$. We prove (see Proposition \ref{prop-gluing-IGR}) that this result follows from our main theorem applied to the Levi quotient group of $R$, which is a reductive group whose rank is smaller than $G$.

In \S \ref{ssec-proof-mainthm}, using induction and the above result, we deduce our main theorem from its anti-tempered version.

\subsection{Statement of the main theorem}
\label{ssec-mainthm}

\begin{constr} \label{constr-gluing-functor}
Let $\Tw(\Par_G^\st)$ be the category of twisted arrows between \emph{standard} parabolic subgroups of $G$, i.e., the objects of $\Tw(\Par_G^\st)$ are standard parabolic subgroups $[P\supset Q]$, and there is a unique morphism from $[P_1\supset Q_1]$ to $[P_2\supset Q_2]$ iff $Q_1\subset Q_2\subset P_2\subset P_1$.

By Theorem \ref{thm-functorial-W-Eis} and Proposition \ref{prop-W-preserve-by-Sph}, we have a functor
\[ \mCW(G,-\supset -)_{\Eis}: \Tw(\Par_G^\st) \to \Sph_{G,x}\hmod,\; [P\supset Q]\mapsto \mCW(G,P\supset Q),
\]
where a morphism from $[P_1\supset Q_1]$ to $[P_2\supset Q_2]$ in $\Tw(\Par_G)$ is sent to the $\Sph_{G,x}$-linear functor 
\[ \Eis: \mCW(G,P_1\supset Q_1) \to \mCW(G,P_2\supset Q_2). \]
View $ \mCW(G,-\supset -)_{\Eis}$ as a diagram in $\Sph_{G,x}\hmod$ and consider its colimit. We have a functor
\[ \gamma_G^L: \colim \mCW(G,-\supset -)_{\Eis} \to \mCI(G,G)= \DMod(\Bun_G) \]
whose restriction to each $ \mCW(G,P\supset Q)$ is the $\Sph_{G,x}$-linear functor
\[\mCW(G,P\supset Q) \subset \mCI(G,Q) \xrightarrow{\Eis}  \mCI(G,G).\]

By Corollary \ref{cor-functorial-W-CT} and Appendix \ref{sect-lax-linear}, we can pass to right adjoints and obtain a functor
\[ \mCW(G,-\supset -)_{\CT_{\temp}}: \Tw(\Par_G^\st)^\op \to\Sph_{G,x}\hmod,\; [P\supset Q]\mapsto \mCW(G,P\supset Q),
\]
where a morphism from $[P_1\supset Q_1]$ to $[P_2\supset Q_2]$ in $\Tw(\Par_G)$ is sent to the $\Sph_{G,x}$-linear functor 
\[ \CT_{\temp}:  \mCW(G, P_2\supset Q_2)\to  \mCW(G, P_1\supset Q_1)  .\]
We have a functor
\[ \gamma_G: \DMod(\Bun_G)= \mCI(G,G)\to\lim \mCW(G,-\supset -)_{\CT_{\temp}}\]
whose value at each $ \mCW(G,P\supset Q)$ is the $\Sph_{G,x}$-linear functor
\[ \mCI(G,G) \xrightarrow{\CT}  \mCI(G,Q) \xrightarrow{\temp_{P}} \mCW(G,P\supset Q).\]

Define
\[ \Glue_G:= \colim \mCW(G,-\supset -)_{\Eis} \simeq \lim \mCW(G,-\supset -)_{\CT_{\temp}}. \]
Then we obtain $\Sph_{G,x}$-linear adjoint functors
\[ \gamma_G^L :\Glue_G \adjoint \DMod(\Bun_G):\gamma_G \]
The functor $\gamma_G$ (resp. $\gamma_G^L$) is called the \emph{automorphic gluing functor} (resp. \emph{automorphic de-gluing functor}).
\end{constr}

The following result is the main theorem of this paper:

\begin{thm}[Automorphic gluing] \label{thm-main-thm}
The adjoint functors 
\[ \gamma_G^L :\Glue_G \adjoint \DMod(\Bun_G):\gamma_G \]
are inverse to each other.
\end{thm}

In \S \ref{ssec-proof-mainthm}, we will reduce the theorem to its $G$-anti-tempered version, which is stated below.

\begin{constr} \label{constr-gluing-atemp}
Define
\[
\Glue_{G\hyphen\atemp}:= \colim_{ R\in (\Par_G^\st), R\neq G } \mCI(G,R)^{G\hyphen\atemp}\simeq \lim_{ R\in (\Par_G^\st)^\op, R\neq G } \mCI(G,R)^{G\hyphen\atemp} ,
\]
where the connecting functors in the colimit (resp. limit) are the functors $\Eis^{G\hyphen\atemp}$ (resp. $\CT^{G\hyphen\atemp}$).

We have a functor
\[\beta_G^L:  \colim_{ R\in (\Par_G^\st), R\neq G } \mCI(G,R)^{G\hyphen\atemp} \to \mCI(G,G)^{G\hyphen\atemp} \]
whose restriction to each $  \mCI(G,R)^{G\hyphen\atemp}$ is the functor
\[\Eis^{G\hyphen\atemp}: \mCI(G,R)^{G\hyphen\atemp}\to   \mCI(G,G)^{G\hyphen\atemp}.\]
We similarly define the functor 
\[ \beta_G: \mCI(G,G)^{G\hyphen\atemp} \to\lim_{ R\in (\Par_G^\st)^\op, R\neq G } \mCI(G,R)^{G\hyphen\atemp} \]
and obtain adjoint functors
\[ \beta_R^L :\Glue_{G\hyphen\atemp}  \adjoint\mCI(G,G)^{G\hyphen\atemp}:\beta_R. \]
\end{constr}

The following theorem is proved in $\S$ \ref{sect-gluing-atemp-proof}.

\begin{thm}[Automorphic gluing for anti-tempered objects] \label{thm-main-thm-atemp}
The adjoint functors
\[ \beta_R^L :\Glue_{G\hyphen\atemp}  \adjoint\mCI(G,G)^{G\hyphen\atemp}:\beta_R \]
are inverse to each other.
\end{thm}

\begin{warn} In Construction \ref{constr-gluing-atemp}, we require $R\neq G$. Otherwise the theorem is trivial.
\end{warn}

\begin{warn} The functor \[\DMod(\Bun_G)\to \lim_{ R\in (\Par_G^\st)^\op, R\neq G } \mCI(G,R) \]
is \emph{not} an equivalence. For example, when $G=T$ is a torus, the RHS is the zero category. 
\end{warn}

\begin{warn} Recall $\Glue_G$ is also equipped with a $\Sph_{G,x}$-action, hence we can define its full subcategory $\Glue_G^{G\hyphen\atemp}$ of $G$-anti-tempered objects, which \emph{a priori} is different from $\Glue_{G\hyphen\atemp}$. However, once Theorem \ref{thm-main-thm} and Theorem \ref{thm-main-thm-atemp} are proven, it will be clear that these two categories are equivalent.

\end{warn}

\subsection{Automorphic gluing for the parabolic categories \texorpdfstring{$\mCI(G,R)$}{I(G,R)}}
\label{ssec-gluing-IGR}

In this subsection, we fix $R\in \Par_G^\st$ with $R\neq G$ and let $M_R:=R/U_R$ be its Levi quotient group.

To deduce Theorem \ref{thm-main-thm} from Theorem \ref{thm-main-thm-atemp}, we need to de-glue $\mCI(G,R)^{G\hyphen\atemp}$ into pieces. In fact, such de-gluing exists even for $\mCI(G,R)$. In this subsection, we describe and prove such de-gluing under the assumption that Theorem \ref{thm-main-thm} is correct for all smaller rank reductive groups.

\begin{constr} Similarly to Construction \ref{constr-gluing-functor}, let $\Tw(\Par_R^\st)\subset \Tw(\Par_G^\st)$ be the full category of twisted arrows between standard parabolic subgroups of $G$ contained in $R$. Define
\[ \Glue_R:= \colim \mCW(G,-\supset -)_{\Eis}|_{\Tw(\Par_R^\st) } \simeq \lim \mCW(G,-\supset -)_{\CT_{\temp}}|_{\Tw(\Par_R^\st)}. \]
In other words,
\[ \Glue_R\simeq \colim_{ [P\supset Q]\in \Tw(\Par_R^\st)  } \mCW(G,P\supset Q) \simeq  \lim_{ [P\supset Q]\in \Tw(\Par_R^\st)^\op  } \mCW(G,P\supset Q), \]
where the connecting functors in the colimit (resp. limit) are the functors $\Eis$ (resp. $\CT_\temp$).

We have a functor
\[ \gamma_R^L:\colim_{ [P\supset Q]\in \Tw(\Par_R^\st)  } \mCW(G,P\supset Q) \to \mCI(G,R) \]
whose restriction to each $ \mCW(G,P\supset Q)$ is the $\Sph_{G,x}$-linear functor
\[\mCW(G,P\supset Q) \subset \mCI(G,Q) \xrightarrow{\Eis}  \mCI(G,R).\]
We similarly define the functor 
\[ \gamma_R: \mCI(G,R)\to \lim_{ [P\supset Q]\in \Tw(\Par_R^\st)^\op  } \mCW(G,P\supset Q) \]
and obtain $\Sph_{G,x}$-linear adjoint functors
\[ \gamma_R^L :\Glue_R \adjoint \mCI(G,R):\gamma_R. \]
\end{constr}

\begin{prop}[Automorphic Gluing for {$\mCI(G,R)$}]\label{prop-gluing-IGR} Suppose the adjoint functors 
\[ \gamma_{M_R}^L :\Glue_{M_R} \adjoint \DMod(\Bun_{M_R}):\gamma_{M_R} \]
are inverse to each other (i.e. Theorem \ref{thm-main-thm} is known for $M_R$), then the adjoint functors
\[ \gamma_R^L :\Glue_R \adjoint \mCI(G,R):\gamma_R \]
are inverse to each other.
\end{prop}

\proof According to our notation, $\mCI(R,R)=\DMod(\Bun_{M_R})$ and we have adjoint functors
\[ \Eis: \mCI(R,R)\adjoint \mCI(G,R):\CT. \]
By Lemma \ref{lem-generator-IGP}, this adjoint pair is monadic.

On the other hand, by Theorem \ref{thm-functorial-W-Eis}, we have a functor
\[\colim_{ [P\supset Q]\in \Tw(\Par_R^\st)  } \mCW(R,P\supset Q)\to   \colim_{ [P\supset Q]\in \Tw(\Par_R^\st)  } \mCW(G,P\supset Q)\]
induced by $\Eis: \mCW(R,P\supset Q)\to  \mCW(G,P\supset Q)$. By definition, this gives a functor $\Glue_{M_R}\to \Glue_{R}$, which we still denote by
\[ \Eis: \Glue_{M_R}\to \Glue_{R}. \]
Passing to right adjoints, by Corollary \ref{cor-functorial-W-CT}, we obtain a functor
\[\lim_{ [P\supset Q]\in \Tw(\Par_R^\st)^\op  } \mCW(G,P\supset Q)\to   \lim_{ [P\supset Q]\in \Tw(\Par_R^\st)^\op  } \mCW(R,P\supset Q)\]
induced by $\CT:  \mCW(G,P\supset Q)\to  \mCW(R,P\supset Q)$. By definition, this gives a functor $\Glue_{R}\to \Glue_{M_R}$, which we still denote by
\[ \CT: \Glue_{R}\to \Glue_{M_R}. \]
By construction, these two functors are adjoint to each other:
\[ \Eis:  \Glue_{M_R}\adjoint \Glue_{R}:\CT.\]
We claim the right adjoint is conservative. To prove the claim, we only need to show each $\CT:  \mCW(G,P\supset Q)\to  \mCW(R,P\supset Q)$ is conservative. By definition, we only need to show $\CT:\mCI(G,Q)\to \mCI(R,Q)$ is conservative. But this follows from the fact that both $\CT:\mCI(G,Q)\to\mCI(Q,Q)$ and $\CT:\mCI(R,Q)\to\mCI(Q,Q)$ are conservative (Lemma \ref{lem-generator-IGP}). This proves the claim that $\CT:  \Glue_{M_R}\to \Glue_{R}$ is conservative. In other words, we have shown that 
\[ \CT: \Glue_{R}\to \Glue_{M_R} \]
is also monadic.

Combining the above two paragraphs, we only need to show the monads 
\[ \mCI(R,R)\to \mCI(G,R)\to \mCI(R,R),\; \Glue_{M_R}\to \Glue_R\to \Glue_{M_R} \]
are compatible with the equivalence $\gamma_{M_R}:\mCI(R,R)\simeq \Glue_{M_R}$. More precisely, by Corollary \ref{cor-functorial-W-CT}, we have a commuative diagram
\[
\xymatrix{
\colim_{ [P\supset Q]\in \Tw(\Par_R^\st)  } \mCW(R,P\supset Q)  \ar[r]^-{\gamma_{M_R}^L} \ar[d]^-\Eis
& \mCI(R,R) \ar[d]^-\Eis\\
 \colim_{ [P\supset Q]\in \Tw(\Par_R^\st) } \mCW(G,P\supset Q)  \ar[r]^-{\gamma_{R}^L}
& \mCI(G,R) ,
}
\]
and we only need to show it is right adjointable along the vertical direction. By definition, we only need to show that, for any $[P_1\supset Q_1] \to [P_2\supset Q_2]  $ in $\Tw(\Par_R^\st) $, the commutative squares
\[
\xymatrix{
\mCW(R,P_1\supset Q_1)  \ar[r]^-{\Eis} \ar[d]^-\Eis
& \mCW(R,P_2\supset Q_2)  \ar[r]^-{\Eis} \ar[d]^-\Eis
& \mCI(R,R) \ar[d]^-\Eis\\
\mCW(G,P_1\supset Q_1)  \ar[r]^-{\Eis}
&\mCW(G,P_2\supset Q_2)  \ar[r]^-{\Eis}
&  \mCI(G,R) 
}
\]
are right adjointable along the vertical direction. By Lemma \ref{lem-CT-preserve-P-temper}, we only need to prove similar thing for
\[
\xymatrix{
\mCI(R, Q_1)  \ar[r]^-{\Eis} \ar[d]^-\Eis
& \mCI(R, Q_2)  \ar[r]^-{\Eis} \ar[d]^-\Eis
& \mCI(R,R) \ar[d]^-\Eis\\
\mCI(G, Q_1)  \ar[r]^-{\Eis}
&\mCI(G, Q_2)  \ar[r]^-{\Eis}
&  \mCI(G,R) .
}
\]
Then we are done by Lemma \ref{lem-Eis-CT-adjointable}.

\qed[Proposition \ref{prop-gluing-IGR}]
\subsection{Deducing the main theorem from its anti-tempered version}
\label{ssec-proof-mainthm}

In this subsection, we deduce Theorem \ref{thm-main-thm} from Theorem \ref{thm-main-thm-atemp}.

To prove Theorem \ref{thm-main-thm}, we proceed by induction on the (semi-simple) rank of the reductive group $G$. When $\Rank(G)=0$, i.e. $G=T$ is a torus, we only need to show $\DMod(\Bun_T) \simeq \mCW(T,T\supset T)$. But this is obvious because any object is $T$-tempered. 

Now we assume $\Rank(G)>0$ and that Theorem \ref{thm-main-thm} is correct for any reductive group with smaller rank. In particular, the assumption of Proposition \ref{prop-gluing-IGR} is satisfied.

To continue, we need the following useful lemma:

\begin{lem} \label{lem-equi-of-adjoint-can-be-check-on-temp&atemp}
Let $\mCC$ and $\mCD$ be $\Sph_G$-module categories and \[L:\mCC \adjoint \mCD:R\] 
be $\Sph_G$-linear adjoint functors. Then $L$ and $R$ are equivalences iff the induced functors 
\[L^{G\hyphen\temp}:\mCC^{G\hyphen\temp} \adjoint \mCD^{G\hyphen\temp}:R^{G\hyphen\temp},\;L^{G\hyphen\atemp}:\mCC^{G\hyphen\atemp} \adjoint \mCD^{G\hyphen\atemp}:R^{G\hyphen\atemp} \]
are equivalences.
\end{lem}

\proof The ``only if'' part is obvious. For the ``if'' part, we just need to prove that the adjunction natural transformations $L\circ R\to \Id_\mCC$ and $\Id_\mCD\to R\circ L$ are invertible. Since any object $x\in \mCC$ (resp. $y\in \mCD$) is an extension of its anti-temperization by its temperization, we only need to show $L\circ R(x)\to x$ (resp. $y\to R\circ L(y)$) is invertible when $x$ (resp. $y$) is either tempered or anti-tempered. If $x$ is tempered, we have $L\circ R(x)\simeq L^{G\hyphen\temp}\circ R^{G\hyphen\temp}(x)$, and the RHS is equivalent to $x$ by assumption. The other cases are similar.

\qed[Lemma \ref{lem-equi-of-adjoint-can-be-check-on-temp&atemp}]

By the lemma, we just need to show that both
\[ (\gamma_G^L)^{G\hyphen\temp} :\Glue_G^{G\hyphen\temp} \adjoint \DMod(\Bun_G)^{G\hyphen\temp}:(\gamma_G)^{G\hyphen\temp}  \]
and
\[ (\gamma_G^L)^{G\hyphen\atemp} :\Glue_G^{G\hyphen\atemp} \adjoint \DMod(\Bun_G)^{G\hyphen\atemp}:(\gamma_G)^{G\hyphen\atemp} \]
are mutually inverse functors. The $G$-tempered case is almost obvious:

\begin{lem} \label{lem-main-thm-temper}
The adjoint functors 
\[ (\gamma_G^L)^{G\hyphen\temp} :\Glue_G^{G\hyphen\temp} \adjoint \DMod(\Bun_G)^{G\hyphen\temp}:(\gamma_G)^{G\hyphen\temp} \]
are inverse to each other.

\end{lem}

\proof By Proposition \ref{prop-G-temp-imply-M-temp-global}, for any $[P\supset Q]\in \Tw(\Par_G^\st)$, we have $\mCW(G,G\supset Q)\subset \mCW(G,P\supset Q)\subset \mCI(G,Q)$. By definition $\mCW(G,G\supset Q)=\mCI(G,Q)^{G\hyphen\temp}$. It follows that $\mCW(G,G\supset Q)\simeq \mCW(G,P\supset Q)^{G\hyphen\temp} $ as full subcategories of $\mCW(G,P\supset Q)$.

Hence we have
\[ \Glue_G^{G\hyphen\temp}\simeq \colim_{[P\supset Q]\in \Tw(\Par_G^\st) } \mCW(G,P\supset Q)^{G\hyphen\temp} \simeq  \colim_{[P\supset Q]\in \Tw(\Par_G^\st) } \mCW(G,G\supset Q).\]
Note that the terms in the RHS do not depend on $P$. In other words, we have a colimit diagram that is constant along the fibers of the forgetful functor $\Tw(\Par_G^\st)\to \Par_G^\st,\; [P\supset Q]\mapsto Q$. Since these fibers are weakly contractible (because they have final objects), we have 
\[ \colim_{[P\supset Q]\in \Tw(\Par_G^\st) } \mCW(G,G\supset Q) \simeq \colim_{Q\in \Par_G^\st} \mCW(G,G\supset Q) \simeq \mCW(G,G\supset G)=\DMod(\Bun_G)\]
as desired.

\qed[Lemma \ref{lem-main-thm-temper}]

For the $G$-anti-tempered case, by Theorem \ref{thm-main-thm-atemp} and Proposition \ref{prop-gluing-IGR}, we have
\begin{eqnarray*}
\DMod(\Bun_G)^{G\hyphen\atemp} &\simeq& \lim_{ R\in (\Par_G^\st)^\op, R\neq G } \mCI(G,R)^{G\hyphen\atemp} \\
&\simeq& \lim_{ R\in (\Par_G^\st)^\op, R\neq G } \Glue_R^{G\hyphen\atemp} \\
&\simeq& \lim_{ R\in (\Par_G^\st)^\op, R\neq G } \big( \lim_{ [P\supset Q]\in \Tw(\Par_R^\st)^\op  } \mCW(G,P\supset Q)^{G\hyphen\atemp} \big)\\
&\simeq&  \lim_{ [P\supset Q]\in \Tw(\Par_G^\st)^\op} \big( \lim_{ R\in (\Par_G^\st)^\op, R\supset P, R\neq G } \mCW(G,P\supset Q)^{G\hyphen\atemp}  \big).
\end{eqnarray*}
Note that the inner limit diagram in the RHS is a constant one (i.e. does not depend on $R$), and the index category $\big\{R\in  (\Par_G^\st)^\op\, \big|\, R\supset P, R\neq G \big \}\subset  (\Par_G^\st)^\op$ is weakly contractible because $P$ is a final object. Hence we have
\begin{eqnarray*}
\DMod(\Bun_G)^{G\hyphen\atemp} &\simeq&  \lim_{ [P\supset Q]\in \Tw(\Par_G^\st)^\op} \big( \lim_{ R\in (\Par_G^\st)^\op, R\supset P, R\neq G } \mCW(G,P\supset Q)^{G\hyphen\atemp}  \big) \\
& \simeq & \lim_{ [P\supset Q]\in \Tw(\Par_G^\st)^\op}\mCW(G,P\supset Q)^{G\hyphen\atemp}  \\
&\simeq& \Glue_G^{G\hyphen\atemp}.
\end{eqnarray*}
It follows from the construction that the underlying functor of this equivalence is $\gamma_G^{G\hyphen\atemp}$. Hence
\[ (\gamma_G^L)^{G\hyphen\atemp} :\Glue_G^{G\hyphen\atemp} \adjoint \DMod(\Bun_G)^{G\hyphen\atemp}:(\gamma_G)^{G\hyphen\atemp} \]
are mutually inverse equivalences as desired.

\qed[Theorem \ref{thm-main-thm}]

\section{Proof of the anti-tempered automorphic gluing theorem}
\label{sect-gluing-atemp-proof}

In this section, we prove Theorem \ref{thm-main-thm-atemp}. We need to show that $\beta^L_G\circ \beta_G \to \Id$ and $\Id\to \beta_G\circ \beta_G^L$ are both equivalences.

In \S \ref{ssec-ff-anti-temper}, we prove the equivalence $\beta^L_G\circ \beta_G \simeq \Id$ using the \emph{Deligne--Lusztig duality on $\Bun_G$}, which is reviewed in \S \ref{ssec-DL}.

The rest of this section (from \S \ref{ssec-surj-anti-temper-prep} to \S \ref{ssec-surj-anti-temper-combo}) is devoted to proving the equivalence $\Id\to \beta_G\circ \beta_G^L$.

In \S \ref{ssec-surj-anti-temper-prep}, we reduce it to the calculation of a certain colimit. More precisely, for any family $(\mCM_P \in \mCI(G,P)^{G\hyphen\atemp})$ compatible with the $\CT$ functors, we need to show that
\[\iota_L^!(\mCM_R)\simeq \colim_{P\in (\Par_G^\st), P\neq G} \CT_{R,*}\circ \Eis_{P\to G}^\enh (\mCM_P) \in \DMod(\Bun_L).\]
Here $R$ is a maximal proper parabolic subgroup of $G$, with Levi quotient group $L$, and $\iota_L^!:\mCI(G,R)\to \DMod(\Bun_M)$ is the functor defined in \S \ref{ssec-IGP}.

Let us first describe our strategy for this calculation:

\begin{strategy} \label{strategy:filtration}
The strategy consists of three major steps.

For each proper parabolic $P \in \Par_G^\st$, write $\alpha_P:=\CT_{R,*}\circ \Eis_{P\to G}^\enh$. In the first step, we define a compatible\footnote{that is, functorial in $P$} exhaustive (increasing) filtration $\alpha_P^{\le \overline{w}}$ on each $\alpha_P$. The filtration is labelled by elements $\overline{w}$ in the poset $W_G/W_R$, where $W_G$ (resp. $W_R$) is the Weyl group of $G$ (resp. $L$) and the partial order on $W_G/W_R$ is the standard one induced from the Bruhat order.

In the second step, we prove the \emph{cancellation lemma}: for any triple $(Q\subset P,\overline{w})$ satisfying certain purely combinatorial conditions, the morphism $\alpha_Q^{= \overline{w}}(\mCM_Q)\to \alpha_P^{= \overline{w}}(\mCM_P)$ is invertible. Here, $ \alpha_P^{= \overline{w}}$ denotes the $\overline{w}$-graded piece of the filtration on $\alpha_P$.

Finally, in the third step, we use the cancellation lemma and some more Weyl combinatorics to show that 
\[\colim_{P\in (\Par_G^\st), P\neq G} \alpha_P^{=\overline{w}}(\mCM_P)
\simeq
\begin{cases}
0 & \mbox{for $\overline{w} \neq \overline{1}$}; \\
\iota_L^!(\mCM_R) & \mbox{for $\overline{w}=\overline{1}$}.
\end{cases}
\]
\end{strategy}

\begin{rem}
Thus, the starting point of the above strategy is the filtration on $\alpha_P:=\CT_{R,*}\circ \Eis_{P\to G}^\enh$. However, its construction is subtle: we will actually construct two filtrations, the \emph{Weyl filtration} and the \emph{dual Weyl filtration}. The former is more natural, but for technical reasons (see below) we are not able to prove the relative cancellation lemma. On the other hand, the dual Weyl filtration is more complicated to define (it requires miraculous parabolic duality), but it allows us to successfully apply the above strategy. It is quite possible that these two filtrations coincide.
\end{rem}

Let us continue with the summary of the contents of the present section. In \S \ref{ssec-surj-anti-temper-weyl}, we define the \emph{Weyl filtration} on each $\alpha_P$. This filtration is induced by a filtration on the prestack 
\[\Bun_G^{P\hgen}\mt_{\Bun_G} \Bun_R.\]
We will then provide a \emph{non-rigorous proof} of the corresponding cancellation lemma. It is non-rigorous because we assume that certain averaging functors (which are right adjoints of the standard functors, i.e., the forgetful functors) commute with \emph{$!$-pushforward functors} (which are left adjoints of the standard functors, i.e., the $!$-pullback functors). In general, it is hard to prove commutativity between left adjoints and right adjoints. For this reason, we \emph{cannot} make this proof rigorous, but we decided to include it anyway for pedagogical reasons and because we believe it can made it work once more technology is established.

To overcome the above issues, we use categorical duality. Namely, in \S \ref{ssec-surj-anti-temper-dual-weyl}, we consider the \emph{dual} Weyl filtration $\alpha_P'^{\le \overline{w}}$ on $\alpha_P:=\CT_{R,*}\circ \Eis_{P\to G}^\enh$ induced by the Weyl filtration on its \emph{dual functor} $\alpha_P^\vee \simeq   \CT_{G\gets P^-}^\enh\circ  \Eis_{R^-,!}$. Here this identification is due to the \emph{parabolic miraculous duality} $\mCI(G,P)^\vee \simeq \mCI(G,P^-)$ of \cite{chen2021thesis}, which we review in \S \ref{ssec-IGP-miraculous}; the Weyl filtration on $ \CT_{G\gets P^-}^\enh\circ  \Eis_{R^-,!}$ is again induced by a filtration of geometric objects.

In \S \ref{ssec-surj-anti-temper-cancel}, we state and prove the cancellation lemma for the dual Weyl filtrations. This can be rigorously proved because we only need to use $!$-pullback and $*$-pushforward functors.

In \S \ref{ssec-surj-anti-temper-combo}, we calculate each graded piece:
\[
\colim_{P\in (\Par_G^\st), P\neq G} \alpha_P'^{=\overline{w}}(\mCM_P)
\simeq
\begin{cases}
\iota_L^!(\mCM_R) & \mbox{for $\overline{w}=\overline{1}$}; \\
0 & \mbox{for $\overline{w} = \overline{w}_\circ$}; \\
0 & \mbox{otherwise}. \\
\end{cases}
\]
As mentioned in Strategy \ref{strategy:filtration}, this colimit will be calculated (in the first and third case) using the cancellation lemma and Weyl combinatorics. The case when $\overline{w}=\overline{w}_\circ$ is the longest element deserves special attention and quite different methods: we will deduce its vanishing from Deligne--Lusztig duality and Theorem \ref{thm-iota!*-temperedness}.

\subsection{Recollection: the Deligne--Lusztig duality}
\label{ssec-DL}

In this subsection, we recall the Deligne--Lusztig duality on $\Bun_G$ proved in \cite{chen2020deligne}. This is one of the main ingredients in the proof of the equivalence $\beta^L_G\circ \beta_G \simeq \Id$.

\begin{constr} Consider the objects
\[ \Delta_!(k), \Delta_*(\omega) \in \DMod(\Bun_G\mt \Bun_G), \]
where $\Delta: \Bun_G\to \Bun_G\mt \Bun_G$ is the diagonal map and $k$ (resp. $\omega$) is the constant (resp. dualizing) D-module. Via the equivalence
\[\DMod(\Bun_G\mt \Bun_G) \simeq \DMod(\Bun_G)\ot_k \DMod(\Bun_G) \simeq \LFun_k(\DMod(\Bun_G)^\vee,\DMod(\Bun_G)) ,\]
we obtain functors between DG categories:
\[ \PsId_{G,!} ,\PsId_{G,*}: \DMod(\Bun_G)^\vee\to \DMod(\Bun_G). \]
\end{constr}

The following result is the main theorem of \cite{gaitsgory2017strange}:

\begin{thm}[Miraculous duality on {$\Bun_G$}] \label{thm-miraculous-duality}
The functor $\PsId_{G,!}: \DMod(\Bun_G)^\vee\to \DMod(\Bun_G)$ is an equivalence.
\end{thm}

As a consequence, we obtain an endo-functor
\[  \PsId_{G,*}\circ \PsId_{G,!}^{-1} \]
on the category $\DMod(\Bun_G)=\mCI(G,G)$. 

On the other hand, we have the following construction:

\begin{constr}
The functor (see Proposition-Construction \ref{propconstr-IGP-functorial})
\[\mCI(G,-)_{\Eis}: \Par_G^\st \to \DGCat,\; R \mapsto \mCI(G,R)\]
induces an endo-functor
\[ \colim_{ R \in \Par_G^\st, R\neq G }  \Eis_{R\to G }^\enh \circ   \CT_{G\gets R}^\enh \]
on $\mCI(G,G)$ equipped with a natural transformation to the identity functor $\Id$.

\end{constr}

The main theorem of \cite{chen2020deligne} states:

\begin{thm}[Deligne--Lusztig duality on {$\Bun_G$}] \label{thm-DL} Up to a cohomological shift, the functor $ \PsId_{G,*}\circ \PsId_{G,!}^{-1} $ is equivalent to
\[\mathrm{DL}_G:=  \coFib(\colim_{ R \in \Par_G^\st, R\neq G }  \Eis_{R\to G }^\enh \circ   \CT_{G\gets R}^\enh \to \Id). \]

\end{thm}

\subsection{\texorpdfstring{The equivalence $\beta^L_G\circ \beta_G \simeq \Id$}{Fully faithful of the G-anti-tempered gluing functor}}
\label{ssec-ff-anti-temper}

The goal of this section is to prove $\beta^L_G\circ \beta_G \simeq \Id$. 

We need to show that, for $\mCF\in \mCI(G,G)^{G\hyphen\atemp}$, the cofiber of $\beta^L_G\circ \beta_G(\mCF)\to \mCF$ vanishes. We prove this by showing that such cofiber is both $G$-tempered and $G$-anti-tempered.

By definition, we have
\[\beta_G^L\circ \beta_G(\mCF) \simeq \colim_{ R \in \Par_G^\st, R\neq G }  \Eis_{R\to G }^\enh \circ   \CT_{G\gets R}^\enh (\mCF).\]
Hence by the Deligne--Lusztig duality (Theorem \ref{thm-DL}), the desired cofiber is isomorphic to $\PsId_{G,*}\circ \PsId_{G,!}^{-1}(\mCF)$ up to a shift.

By \cite[\S 4.4.3]{drinfeld2015compact}, the image of the functor $\PsId_{G,*}$ is contained in $\DMod(\Bun_G)^{{*\hgen}}$. By \cite[Theorem B]{beraldo2021geometric}, we have $\DMod(\Bun_G)^{{*\hgen}}\subset \DMod(\Bun_G)^{G\htemp}$. Hence the object $\PsId_{G,*}\circ \PsId_{G,!}^{-1}(\mCF)$ is $G$-tempered as desired.

On the other hand, since $ \Eis_{R\to G }^\enh$  and $ \CT_{G\gets R}^\enh (\mCF)$ are $\Sph_{G,x}$-linear, the objects $\mCF$ and $ \beta^L_G\circ \beta_G(\mCF)$ are both $G$-anti-tempered. Hence their cofiber is also $G$-anti-tempered as desired.

Combining the above two paragraphs, we obtain that $\beta^L_G\circ \beta_G(\mCF)\to \mCF$ is an isomorphism for any $\mCF\in \mCI(G,G)^{G\hyphen\atemp}$.

\subsection{\texorpdfstring{The equivalence $\Id\to \beta_G\circ \beta_G^L$}{Fully faithful of the G-anti-tempered de-gluing functor}: reduction to a colimit calculation}
\label{ssec-surj-anti-temper-prep}

The goal of the rest of this section is to prove $\Id\to \beta_G\circ \beta_G^L$. In this subsection, we reduce it to a colimit calculation (see Goal \ref{goal-surj-anti-temper-colimit} below).

By definition, an object of
\[\Glue_{G\hyphen\atemp}\simeq \lim_{ P\in (\Par_G^\st)^\op, P\neq G } \mCI(G,P)^{G\hyphen\atemp}\]
is a family of objects $\mCM_P\in \mCI(G,P)^{G\hyphen\atemp}$ for $P\in (\Par_G^\st)^\op, P\neq G,$ compatible with the CT functors. From now on, we fix such a family and denote the corresponding object by $\mCM_\bullet\in \Glue_{G\hyphen\atemp}$. Our goal is to show that the arrow $\mCM_\bullet\to \beta_G\circ \beta_G^L(\mCM_\bullet)$ is invertible. For this, we need to prove that $\mCM_P \to \ev_P\circ \beta_G\circ \beta_G^L(\mCM_\bullet)$ is invertible for any $P\in (\Par_G^\st)^\op, P\neq G$, where
\[\ev_P: \Glue_{G\hyphen\atemp}\simeq \lim_{ P\in (\Par_G^\st)^\op, P\neq G } \mCI(G,P)^{G\hyphen\atemp}\to \mCI(G,P)^{G\hyphen\atemp} \]
is the evaluation functor. Clearly, it suffices to check that
\[\mCM_R \to \ev_R\circ \beta_G\circ \beta_G^L(\mCM_\bullet)\]
is invertible for any \emph{maximal (proper)} standard parabolic subgroup $R$ of $G$. From now on, we fix such an $R$ and denote its Levi quotient group by $L$.

Note that $\ev_P$ has a left adjoint (the inserting functor)
\[ \ins_P: \mCI(G,R)^{G\hyphen\atemp}\to \colim_{ P\in (\Par_G^\st), P\neq G } \mCI(G,P)^{G\hyphen\atemp} \simeq \Glue_{G\hyphen\atemp}.\]
It follows from the definitions that
\[\mCM_\bullet \simeq \colim_{ P\in (\Par_G^\st), P\neq G } \ins_P\circ \ev_P( \mCM_\bullet) \simeq \colim_{ P\in (\Par_G^\st), P\neq G } \ins_P(\mCM_P).\]
Thus, we obtain
\[ \ev_R\circ \beta_G\circ \beta_G^L(\mCM_\bullet) \simeq \colim_{ P\in (\Par_G^\st), P\neq G } \ev_R\circ \beta_G\circ \beta_G^L\circ \ins_P(\mCM_P) \simeq \colim_{P\in (\Par_G^\st), P\neq G} \CT_{G\gets R}^\enh\circ \Eis_{P\to G}^\enh (\mCM_P),
 \]
where for $P_1\subset P_2$, the connecting morphism is given by
\[
\begin{aligned}
\CT_{G\gets R}^\enh\circ \Eis_{P_1\to G}^\enh (\mCM_{P_1}) \simeq \CT_{G\gets R}^\enh\circ \Eis_{P_1\to G}^\enh \circ \CT_{P_2\gets P_1}^\enh (\mCM_{P_2}) \simeq\\
\simeq  \CT_{G\gets R}^\enh\circ \Eis_{P_2\to G}^\enh  \circ \Eis_{P_1\to P_2}^\enh \circ \CT_{P_2\gets P_1}^\enh (\mCM_{P_2}) \to  \CT_{G\gets R}^\enh\circ \Eis_{P_2\to G}^\enh  (\mCM_{P_2}) .
\end{aligned}
\]
Hence, we only need to show that the composition 
\[ \mCM_R\to  \CT_{G\gets R}^\enh\circ \Eis_{R\to G}^\enh (\mCM_R)\to  \colim_{P\in (\Par_G^\st), P\neq G} \CT_{G\gets R}^\enh\circ \Eis_{P\to G}^\enh (\mCM_P)\]
is invertible. Since the functor $\iota_L^!:\mCI(G,R)\to \DMod(\Bun_{L})$ is conservative, it suffices to check that the arrow
\[ \iota_L^! (\mCM_R)\to \colim_{P\in (\Par_G^\st), P\neq G} \iota_L^!\circ \CT_{G\gets R}^\enh\circ \Eis_{P\to G}^\enh (\mCM_P) \] 
is invertible. By Proposition-Construction \ref{propconstr-IGP-functorial}, $ \iota_L^!\circ \CT_{G\gets R}^\enh\simeq \CT_{R,*}$. Hence, we have reduced our proof to the following goal:

\begin{goal}\label{goal-surj-anti-temper-colimit}
 We need to show the morphism 
 \[\iota_L^! (\mCM_R)\to \colim_{P\in (\Par_G^\st), P\neq G}  \CT_{R,*}\circ \Eis_{P\to G}^\enh (\mCM_P)\]
 is invertible.
\end{goal}

\subsection{\texorpdfstring{The equivalence $\Id\to \beta_G\circ \beta_G^L$}{Fully faithful of the G-anti-tempered de-gluing functor}: the Weyl filtration}
\label{ssec-surj-anti-temper-weyl}

In this subsection, we introduce the Weyl filtration on the functor 
\[\alpha_P:=\CT_{R,*}\circ \Eis_{P\to G}^\enh\]
and provide a non-rigorous proof for the cancellation lemma. 

By definition, $\alpha_P$ is the restriction of
\begin{equation*}
\DMod(\Bun_G^{P\hgen}) \xrightarrow{!\mathrm{-push}} \DMod(\Bun_G) \xrightarrow{!\mathrm{-pull}} \DMod(\Bun_R) \xrightarrow{*\mathrm{-push}} \DMod(\Bun_L)
\end{equation*}
to the full subcategory $\mCI(G,P) \subseteq \DMod(\Bun_G^{P\hgen})$. By Remark \ref{rem-BunG-Pgen-to-BunG-pseudo-proper} and the pseudo-proper base-change theorem (see \cite[Corollary 1.5.4]{gaitsgory2015atiyah}), the above composition is equivalent to
\[
\begin{aligned}
 \delta_P: \DMod(\Bun_G^{P\hgen}) \xrightarrow{!\mathrm{-pull}} \DMod(\Bun_G^{P\hgen}\mt_{\Bun_G} \Bun_R) \to\\
 \xrightarrow{!\mathrm{-push}} \DMod(\Bun_R) \xrightarrow{*\mathrm{-push}} \DMod(\Bun_L). 
\end{aligned}
\]
Clearly, we have:
\[\Bun_G^{P\hgen}\mt_{\Bun_G} \Bun_R \simeq \mathbf{Maps}_{\mathrm{gen}}(X,\mBB R\gets P\backslash G/R ) . \]
We are going to define a filtration on this prestack. In plain words, this prestack classifies $G$-torsors on $X$ equipped with an $R$-structure and a generic $P$-structure; then, roughly speaking, the filtration constructed below measures the relative position of the generic $R$-torsor and $P$-torsor inside the generic $G$-torsor.

\begin{constr} \label{constr-filtration-geo}
The double quotient $P\backslash G/R$ has a stratification labelled by $W_P\backslash W_G/W_R$, where $W_P$ (resp. $W_R$, $W_G$) is the Weyl group for the its Levi subgroup. Namely, for any $w\in W_G$, we have a stratum\footnote{To simplify the notations, we fix a section $W_G=N_G(T)/T \to N_G(T)$.} $P\backslash PwR/R$ which only depends on the image of $w$ in $W_P\backslash W_G/W_R$. 

Consider the standard partial order on $W_G/W_R$ induced by the Bruhat order of $W_G$. In other words, for $\overline{w}_1,\overline{w}_2\in W_G/W_R$, we have the following equivalent conditions:
\begin{itemize}
	\item $\overline{w}_1\le \overline{w}_2$;
	\item For any representatives $w_1,w_2\in W_G$ of $\overline{w}_1,\overline{w}_2$, the orbit $Bw_1R$ is contained in the closure of $Bw_2R$;
	\item For representatives $w_1,w_2\in W_G$ of $\overline{w}_1,\overline{w}_2$ that are maximal for the Bruhat order, $w_1\le w_2$.
\end{itemize}
For any $\overline{w}\in W_G/W_R$ and any representative $w\in W_G$ of it, consider the closure $\overline{PwR}$ of $PwR$ in $G$. Note that it only depends on $\overline{w}$. Define
\[ (\Bun_G^{P\hgen}\mt_{\Bun_G} \Bun_R)^{\le \overline{w}} := \mathbf{Maps}_{\mathrm{gen}}(X,\mBB R\gets P\backslash \overline{PwR}/R )
\]
Then for $\overline{w}_1\le \overline{w}_2$, we have a map
\[ (\Bun_G^{P\hgen}\mt_{\Bun_G} \Bun_R)^{\le \overline{w}_1} \to (\Bun_G^{P\hgen}\mt_{\Bun_G} \Bun_R)^{\le \overline{w}_2} . \]
This is a closed embedding (up to non-reduced structures) because $\overline{Pw_1R}\to \overline{Pw_2R}$ is a closed embedding. Note that, for the longest element $w_\circ$ in $W_G$, we have
\[ (\Bun_G^{P\hgen}\mt_{\Bun_G} \Bun_R)^{\le \overline{w}_\circ} =  (\Bun_G^{P\hgen}\mt_{\Bun_G} \Bun_R).  \]
Hence, we obtain a filtration of $\Bun_G^{P\hgen}\mt_{\Bun_G} \Bun_R $ indexed by the poset $ W_G/W_R$.

\end{constr}

\begin{constr} \label{constr-weyl-filtration-alpha}
The above filtration on $\Bun_G^{P\hgen}\mt_{\Bun_G} \Bun_R $ induces a filtration on the functor $\delta_P$ such that $\delta_P^{\le \overline{w}}$ is the composition
\[
\begin{aligned}
 \DMod(\Bun_G^{P\hgen}) \xrightarrow{!\mathrm{-pull}} \DMod((\Bun_G^{P\hgen}\mt_{\Bun_G} \Bun_R)^{\le \overline{w}}) \to\\
 \xrightarrow{!\mathrm{-push}} \DMod(\Bun_R) \xrightarrow{*\mathrm{-push}} \DMod(\Bun_L).
\end{aligned}
\]
 Note that we have $\delta_P^{\le \overline{w}_\circ}= \delta_P$. We obtain a filtration on $\alpha_P$ by restricting to $\mCI(G,P)\subset \DMod(\Bun_G^{P\hgen})  $.
\end{constr}

\begin{constr} \label{constr-weyl-filtration-functorial}
The above filtrations are compatible with the connecting morphisms in $\colim_P  \alpha_P (\mCM_P)$, i.e., compatible with the natural transfomations
\begin{equation} \label{eqn-nt-alpha}
\begin{aligned}
 \alpha_Q \circ \CT_{P\gets Q}^\enh = \CT_{R,*}\circ \Eis_{Q\to G}^\enh \circ \CT_{P\gets Q}^\enh \simeq \CT_{R,*}\circ \Eis_{P\to G}^\enh \circ \Eis_{Q\to P}^\enh \circ \CT_{P\gets Q}^\enh \to\\ 
\to \CT_{R,*}\circ \Eis_{P\to G}^\enh= \alpha_P.
\end{aligned}
\end{equation}
 In other words, we have canonical natural transformations
\[\alpha_Q^{\le \overline{w}} \circ \CT_{P\gets Q}^\enh \to \alpha_P^{\le \overline{w}}. \]

In more detail, for $Q\subset P$, consider the map $p_{Q\to P}^\enh:\Bun_G^{Q\hgen}\to \Bun_G^{P\hgen}$. Recall (see Definition \ref{defn-Eis-CT-enh}) that
\[
\CT_{P\gets Q}^\enh\simeq \Av_*^{U_Q(\mBA)}\circ (p_{Q\to P}^\enh)^!|_{\mCI(G,P)}.
\]
Hence we have a natural transformation
\begin{equation} \label{eqn-proof-surj-anti-temper-3}
 \oblv^{U_Q(\mBA)} \circ \CT_{P\gets Q}^\enh \to (p_{Q\to P}^\enh)^!|_{\mCI(G,P)}
\end{equation}
between functors $\mCI(G,P)\to \DMod( \Bun_G^{Q\hgen})$.

On the other hand, the obvious commutative sqaure
\[
\xymatrix{
	(\Bun_G^{Q\hgen}\mt_{\Bun_G} \Bun_R)^{\le \overline{w}} \ar[r] \ar[d] &
	\Bun_G^{Q\hgen} \ar[d]\\
	(\Bun_G^{P\hgen}\mt_{\Bun_G} \Bun_R)^{\le \overline{w}} \ar[r] &
	\Bun_G^{P\hgen}
}
\]
induces a natural transformation 
\begin{equation} \label{eqn-proof-surj-anti-temper-2}
\delta_Q^{\le \overline{w}}\circ  (p_{Q\to P}^\enh)^!\to \delta_P^{\le \overline{w}}.
\end{equation}
Combining with the previous paragraph, we obtain the desired natural transformation
\begin{equation} \label{eqn-proof-surj-anti-temper-4}
 \alpha_Q^{\le \overline{w}} \circ \CT_{P\gets Q}^\enh \simeq \delta_Q^{\le \overline{w}}\circ \oblv^{U_Q(\mBA)} \circ \CT_{P\gets Q}^\enh \to \delta_Q^{\le \overline{w}}\circ (p_{Q\to P}^\enh)^!|_{\mCI(G,P)} \to \delta_P^{\le \overline{w}}|_{\mCI(G,P)} \simeq \alpha_P^{\le \overline{w}} . 
 \end{equation}
\end{constr}

\begin{constr} \label{constr-gr-weyl}
For any $\overline{w}\in W_G/W_R$, we define
\[ \delta_P^{< \overline{w}}:= \colim_{ \overline{w}'<\overline{w} }  \delta_P^{\le \overline{w}'},\; \delta_P^{= \overline{w}}:= \mathrm{coFib}(\delta_P^{< \overline{w}}\to \delta_P^{\le \overline{w}}  ).\]
It follows from the definitions that:
\begin{itemize}
	\item If $\overline{w}$ is minimal in the fiber of $W_G/W_R\to W_P\backslash W_G/W_R $ containing it, then 
	\[
	\begin{aligned}
	\delta_P^{= \overline{w}}: \DMod(\Bun_G^{P\hgen}) \xrightarrow{!\mathrm{-pull}} \DMod((\Bun_G^{P\hgen}\mt_{\Bun_G} \Bun_R)^{= \overline{w}}) \to\\
	\xrightarrow{*\mathrm{-push}} \DMod(\Bun_G^{P\hgen}\mt_{\Bun_G} \Bun_R)  \xrightarrow{!\mathrm{-push}} \DMod(\Bun_R) \xrightarrow{*\mathrm{-push}} \DMod(\Bun_L);
	\end{aligned}
	\]
	moreover
	\[ (\Bun_G^{P\hgen}\mt_{\Bun_G} \Bun_R)^{= \overline{w}}\simeq  \mathbf{Maps}_{\mathrm{gen}}(X,\mBB R\gets P\backslash PwR /R ).\]

	\item If $\overline{w}$ is not minimal in the fiber of $W_G/W_R\to W_P\backslash W_G/W_R $ containing it, then $\delta_P^{= \overline{w}}\simeq 0$.
\end{itemize}
Restricting to $\mCI(G,P)$, we obtain the functors $\alpha_P^{= \overline{w}}$.
\end{constr}

\begin{rem} \label{rem-geometric-stratum-as-Bungen}
Write $R^w:=wRw^{-1}$. Note that $ P\backslash PwR /R\simeq \mBB (P\cap R^w)$ and $\mBB R\simeq \mBB R^w$. Hence we have
\begin{equation} \label{eqn-geometric-stratum}
(\Bun_G^{P\hgen}\mt_{\Bun_G} \Bun_R)^{= \overline{w}} \simeq \Bun_{R^w}^{(P\cap R^w)\hgen}.
\end{equation}
Since any inner automorphism of $G$ induces the identity map on $\mBB G$, the isomorphism $\mBB R\simeq \mBB R^w$ is defined over $\mBB G$. It follows that via (\ref{eqn-geometric-stratum}), the map $(\Bun_G^{P\hgen}\mt_{\Bun_G} \Bun_R)^{= \overline{w}}\to \Bun_G^{P\hgen}$ is given by the natural map $\Bun_{R^w}^{(P\cap R^w)\hgen}\to \Bun_G^{P\hgen}$ induced by $R^w\to G$.

\end{rem}

\begin{defn} For fixed $\overline{w}\in W_G/W_R$, let $\Par_G^{\overline{w}}\subset \Par_G^\st\setminus\big\{ G\big\}$ be the sub-poset of parabolics $P$ such that $\overline{w}$ is minimal in the fiber of $W_G/W_R\to W_P\backslash W_G/W_R $ containing it.
\end{defn}

Now we state the cancellation lemma (or in fact conjecture) for the Weyl filtration on the functor $\alpha_P:=\CT_{R,*}\circ \Eis_{P\to G}^\enh$. We write it as a conjecture and provide a non-rigorous proof as a remark.

\begin{conj}  [Cancellation lemma for the Weyl filtration] \label{conj-cancel-weyl} 
For fixed $\overline{w}\in W_G/W_R$ and any representative $w$ of it, let $Q,P\in \Par_G^{\overline{w}}$ with $Q\subset P$ satisfy the following conditions:
\begin{itemize}
	\item[(i)] $Q\cap R^w=P\cap R^w$;
	\item[(ii)] $U_Q\cap M_P\subset U_R^w$, where $U_Q$ (resp. $U_R$) is the unipotent radical of $Q$ (resp. $R$), $M_P$ is the Levi subgroup of $P$ and $U_R^w:=wU_Rw^{-1}$.
\end{itemize}
Then (\ref{eqn-proof-surj-anti-temper-4}) induces an equivalence
\[\alpha_Q^{= \overline{w}}\circ \CT_{P\gets Q}^\enh  \simeq \alpha_P^{= \overline{w}}\]
between functors $\mCI(G,P)\to \DMod(\Bun_L)$

\end{conj}

\begin{rem}[A non-rigorous proof] \label{rem-non-rigorous-proof}
By Condition (i) and Remark \ref{rem-geometric-stratum-as-Bungen}
\[(\Bun_G^{P\hgen}\mt_{\Bun_G} \Bun_R)^{= \overline{w}}\simeq (\Bun_G^{Q\hgen}\mt_{\Bun_G} \Bun_R)^{= \overline{w}}.\]
We claim this implies the natural transformation 
\[\delta_Q^{= \overline{w}}\circ  (p_{Q\to P}^\enh)^!\to \delta_P^{= \overline{w}}\]
induced by (\ref{eqn-proof-surj-anti-temper-2}) is invertible. Indeed, unwinding the definitions, we only need to show that the following commutative square (given by the usual base-change isomorphism)
\[
\xymatrix{
	\DMod( (\Bun_G^{Q\hgen}\mt_{\Bun_G} \Bun_R)^{= \overline{w}} ) \ar[d]^-{*\mathrm{-push}}  &
	\DMod( (\Bun_G^{P\hgen}\mt_{\Bun_G} \Bun_R)^{= \overline{w}} )\ar[l]^-{!\mathrm{-pull}}  \ar[d]^-{*\mathrm{-push}} \\
	\DMod( \Bun_G^{Q\hgen}\mt_{\Bun_G} \Bun_R) &
	\DMod( \Bun_G^{P\hgen}\mt_{\Bun_G} \Bun_R) \ar[l]^-{!\mathrm{-pull}} 
}
\]
is left adjointable along horizontal direction. But this follows from the construction of pseudo-proper $!$-pushforward (see \cite[Corollary 1.5.4]{gaitsgory2015atiyah}).

Hence it remains to show that the natural transformation 
\[\delta_Q^{= \overline{w}}\circ \oblv^{U_Q(\mBA)} \circ \CT_{P\gets Q}^\enh \simeq \delta_Q^{= \overline{w}}\circ \oblv^{U_Q(\mBA)} \circ \Av_*^{U_Q(\mBA)} \circ  (p_{Q\to P}^\enh)^!|_{\mCI(G,P)} \to \delta_Q^{= \overline{w}}\circ  (p_{Q\to P}^\enh)^!|_{\mCI(G,P)} \]
induced by (\ref{eqn-proof-surj-anti-temper-3}) is invertible. Now \emph{imagine} the adelic group $U_Q(\mBA)$ exists in algebraic geometry and the analogy in Remark \ref{rem-adelic-picture} can be made precise. For any object $\mCM_P\in \mCI(G,P)$, since it is $U_P(\mBA)$-equivariant, so is $(p_{Q\to P}^\enh)^!(\mCM_P)$. Note that $U_Q/U_P\simeq U_Q\cap M_P$. Hence 
\[ \oblv^{U_Q(\mBA)}\circ  \Av_*^{U_Q(\mBA)}\circ (p_{Q\to P}^\enh)^!(\mCM_P) \simeq \oblv^{U_Q\cap M_P(\mBA)}\circ \Av_*^{U_Q\cap M_P(\mBA)}\circ (p_{Q\to P}^\enh)^!(\mCM_P).\]
Hence we only need to show 
\[ \delta_Q^{= \overline{w}}\to  \delta_Q^{= \overline{w}} \circ \oblv_*^{U_Q\cap M_P(\mBA)} \circ \Av_*^{U_Q\cap M_P(\mBA)}\] is invertible. Note that the map 
\[\Bun_{R^w}^{(Q\cap R^w)\hgen} \simeq  (\Bun_G^{Q\hgen}\mt_{\Bun_G} \Bun_R)^{= \overline{w}}\to \Bun_L\] factors through $\Bun_{L^w}^{(Q\cap L^w)\hgen}$. We can replace $\delta_Q^{= \overline{w}}$ by the composition
\[\delta': \DMod(\Bun_G^{P\hgen}) \xrightarrow{!\mathrm{-pull}} \DMod( \Bun_{R^w}^{(Q\cap R^w)\hgen} )   \xrightarrow{\mathrm{push}} \DMod( \Bun_{L^w}^{(Q\cap L^w)\hgen}) .\]
By Condition (ii), $U_Q\cap M_P\subset \ker(Q\cap R^w\to Q\cap L^w)$, hence $U_Q\cap M_P(\mBA)$ also acts on $ \Bun_{R^w}^{(Q\cap R^w)\hgen}$ and stabilizes the fibers of $ \Bun_{R^w}^{(Q\cap R^w)\hgen}\to \Bun_{L^w}^{(Q\cap L^w)\hgen}$. This indicates
\[\delta'\simeq \delta'\circ \oblv_*^{U_Q\cap M_P(\mBA)}\circ \Av_*^{U_Q\cap M_P(\mBA)}\].

The above argument is non-rigorous at least for the following two reasons:
\begin{itemize}
	\item[(1)] We did not define the pushforward functor along $\Bun_{R^w}^{(Q\cap R^w)\hgen} \to \Bun_{L^w}^{(Q\cap L^w)\hgen}$, which is neither pseudo-proper nor represented by algebraic stacks;

	\item[(2)] We did not define the adelic-actions and verify various properties of the averaging functors.
\end{itemize}
Problem (1) is not essential. But Problem (2) seems out of reach with our current techniques.

\end{rem}

\subsection{Recollection: the parabolic miraculous duality \texorpdfstring{$\mCI(G,P)^\vee \simeq \mCI(G,P^-)$}{ }}
\label{ssec-IGP-miraculous}

As explained in the beginning of this section, we need to use the duality between $\mCI(G,P)$ and $\mCI(G,P^-)$. The following result is a global analogue of Theorem \ref{thm-inv-inv-duality} and a generalization of the miraculous duality on $\Bun_G$ (Theorem \ref{thm-miraculous-duality}):

\begin{thm}[\!\!{\cite[Theorem E, Theorem 5.3.5]{chen2021thesis}}] \label{thm-inv-inv-duality-global}
Let $P$ and $P^-$ be opposite parabolic subgroups of $G$ and $M:=P\cap P^-$ be their Levi subgroup. Then the categories $\mCI(G,P)$ and $\mCI(G,P^-)$ are canonically dual to each other. Moreover,
\begin{itemize}
	\item[(1)] In the case $P=P^-=G$, we recover the miraculous duality on $\Bun_G$;
	\item[(2)] Via this duality and the miraculous duality on $\Bun_G$, the functors
	\[ \Eis_{P\to G}^\enh:\mCI(G,P)\to \DMod(\Bun_G),\;\; \Eis_{P^-\to G}^\enh:\mCI(G,P^-)\to \DMod(\Bun_G) \]
	are conjugate to each other;
	\item[(3)] Via this duality and the miraculous duality on $\Bun_M$, the functors
	\[ \iota_{M,!}:\DMod(\Bun_M)\to \mCI(G,P),\;\; \iota_{M,!}^-:\DMod(\Bun_M)\to \mCI(G,P^-) \]
	are conjugate to each other.
\end{itemize}
\end{thm}

As in the local case (see \S \ref{ssec-SI}), the above theorem formally implies:

\begin{cor} \label{cor-inv-inv-duality-functors-global}
Via the above dualities, we have:
\begin{itemize}
	\item[(1)] The functors 
	\begin{eqnarray*}
		\Eis_{P\to G}^\enh:\mCI(G,P)\to \DMod(\Bun_G),\;\; \CT_{G\gets P^-}^\enh:\DMod(\Bun_G)\to \mCI(G,P^-)
	\end{eqnarray*}
	are dual to each other;
	\item[(2)] The functors 
	\begin{eqnarray*}
		\iota_{M}^!: \mCI(G,P) \to \DMod(\Bun_M) ,\;\; \iota_{M,!}^-: \DMod(\Bun_M) \to \mCI(G,P^-)
	\end{eqnarray*}
	are dual to each other;
	\item[(3)] The functors 
	\begin{eqnarray*}
		\iota_{M}^*: \mCI(G,P) \to \DMod(\Bun_M) ,\;\; \iota_{M,*}^-: \DMod(\Bun_M) \to \mCI(G,P^-)
	\end{eqnarray*}
	are dual to each other.
\end{itemize}
\end{cor}

\begin{rem} \label{rem-inv-inv-duality-functors-global}
In fact, for any\footnote{Except in this remark, we always use $Q$ to denote a parabolic subgroup \emph{smaller} than $P$. We give up this convention here to be compatible with \cite{chen2021thesis}.} $ P\subset Q$ and $ P^-\subset Q^-$, the functors
	\begin{eqnarray*}
		\Eis_{P\to Q}^\enh:\mCI(G,P)\to \mCI(G,Q),\;\; \CT_{Q^-\gets P^-}^\enh:\mCI(G,Q^-)\to \mCI(G,P^-)
	\end{eqnarray*}
are also dual to each other via the above dualities. This is implicit in \cite{chen2021thesis}. Let us explain the main ideas but treat the constructions in \emph{loc.cit.} as a blackbox.

We can assume $P$ and $Q$ are standard. In \cite[\S 5.3]{chen2021thesis} (see also \cite{chen2020deligne}), the author considered a certain algebraic stack $\overline{\Bun}_G$ equipped with a closed filtration $(\overline{\Bun}_G)_{\le P}$ labelled by elements $P\in \Par_G^{\st}$. There is a certain functor 
\[\mbK: \Par_G^\st\to \DMod( \overline{\Bun}_G)\]
such that, for each $P\in \Par_G^{\st}$, the object $\mbK(P)$ is supported on the closed substack $(\overline{\Bun}_G)_{\le P}$ of $\overline{\Bun}_G$. There are maps
\[ \overline{\Delta}^\enh_{\le P}: (\overline{\Bun}_G)_{\le P}\to \Bun_G^{P\hgen}\mt \Bun_G^{P^-\hgen} \]
functorial in $P$. The $!$-pushforward of $\mbK(P)$ along this map, which we denote by $\mCK_P$, is contained in the full subcategory
\[ \mCI(G,P)\ot \mCI(G,P^-) \simeq \mCI(G\mt G,P\mt P^-) \subset \DMod( \Bun_G^{P\hgen}\mt \Bun_G^{P^-\hgen}  ) .\]
Then \cite[Theorem 5.3.5]{chen2021thesis} says that $\mCK_P$ is the unit object for a duality between $\mCI(G,P)$ and $\mCI(G,P^-)$. 

Hence, to show that $\Eis_{P\to Q}^\enh$ and $\CT_{Q^-\gets P^-}^\enh$ are mutually dual, we only need to provide a canonical equivalence
\[ \Eis_{P\to Q}^\enh\ot \Id( \mCK_P ) \simeq \Id\ot \CT_{Q^-\gets P^-}^\enh(\mCK_Q) \]
in $\mCI(G,Q)\ot \mCI(G,P^-)$. This can be constructed as follows. Recall that the maps $\overline{\Delta}^\enh_{\le P}$ and the objects $\mbK(P)$ are functorial in $P$. Hence we obtain a morphism
\[ \Eis_{P\to Q}^\enh\ot \Eis_{P^-\to Q^-}^\enh (\mCK_P)\to \mCK_Q,\]
which induces a morphism
\[ \Eis_{P\to Q}^\enh\ot \Id( \mCK_P ) \to \Id\ot \CT_{Q^-\gets P^-}^\enh(\mCK_Q) .\]
We only need to show that this morphism is invertible.

When $Q=G$, this is exactly the content of the proof of \cite[Theorem 5.3.5]{chen2021thesis} in \cite[\S 6.13.2]{chen2021thesis}. The proof there reduces the claim to \cite[Goal 6.4.3]{chen2021thesis}, and then to \cite[Goal 6.4.17]{chen2021thesis}. In the general case, we can proceed as in \emph{loc.cit.} and reduce it to \cite[Goal 6.4.18]{chen2021thesis}.

\end{rem}

\begin{rem} It is very possible that the dualities in Theorem \ref{thm-inv-inv-duality-global} are compatible with all the functorialities of $\mCI(-,-)$ in Proposition-Construction \ref{propconstr-IGP-functorial}. Up to homotopy coherence, this follows from Corollary \ref{cor-inv-inv-duality-functors-global} and Remark \ref{rem-inv-inv-duality-functors-global}. However, one needs more geometric inputs to give a homotopy coherent proof, e.g., one needs to relate $\overline{\Bun}_G$ to $\overline{\Bun}_M$.
\end{rem}

\begin{rem} It is very possible that the dualities in Theorem \ref{thm-inv-inv-duality-global} are compatible with the $\Sph_{G,x}$-actions on $\mCI(G,P)$ and $\mCI(G,P^-)$. This is \emph{not} obvious because there is no $\Sph_{G,x}$-action on $\overline{\Bun}_G$.
\end{rem}
\subsection{\texorpdfstring{The equivalence $\Id\to \beta_G\circ \beta_G^L$}{Fully faithful of the G-anti-tempered de-gluing functor}: the dual Weyl filtration}
\label{ssec-surj-anti-temper-dual-weyl}

In this subsection, we introduce the dual Weyl filtration on the functor 
\[\alpha_P: \mCI(G,P) \xrightarrow{\Eis_{P\to G}^\enh} \DMod(\Bun_G) \xrightarrow{\CT_{G\gets R}^\enh} \mCI(G,R) \xrightarrow{\iota_L^!} \DMod(\Bun_L).\]

Via the dualities in Corollary \ref{cor-inv-inv-duality-functors-global}, the above functor is dual to
\[  \DMod(\Bun_L) \xrightarrow{\iota_{L,!}^-} \mCI(G,R^-) \xrightarrow{\Eis_{R^-\to G}^\enh} \DMod(\Bun_G)\xrightarrow{\CT_{G\gets P^-}^\enh}  \mCI(G,P^-). \]
In other words, we have
\[\alpha_P^\vee \simeq \CT_{G\gets P^-}^\enh\circ \Eis_{R^-,!}.\]
Consider the commutative diagram
\[
\xymatrix{
\Bun_G^{P^-\hgen}\mt_{\Bun_G} \widetilde{\Bun}_{R^-} \ar[d] \ar[r] &
\Bun_G^{P^-\hgen} 
 \ar[d] \\
\widetilde{\Bun}_{R^-}\ar[r] & \Bun_G.
}
\]
By the base-change isomorphism, the functor $\CT_{G\gets P^-}^\enh\circ \Eis_{R^-,!}$ is equivalent to the composition of\footnote{Here we do not need to distinguish the $!$-pushforward and the $*$-pushforward along $\Bun_G^{P^-\hgen}\mt_{\Bun_G} \widetilde{\Bun}_{R^-}\to \Bun_G^{P^-\hgen}$ because it is a disjoint union of proper maps.}
\[
\begin{aligned}
 \epsilon_{P^-}: \DMod(\Bun_L) \xrightarrow{*\hpull} \DMod(\Bun_{R^-}) \xrightarrow{!\hpush} \DMod( \widetilde{\Bun}_{R^-} ) \to\\ \xrightarrow{!\hpull} \DMod(\Bun_G^{P^-\hgen}\mt_{\Bun_G} \widetilde{\Bun}_{R^-}) \xrightarrow{\mathrm{push}} \DMod(\Bun_G^{P^-\hgen})
\end{aligned}
 \]
with
\[ \Av_*^{U_P^-(\mBA)}: \DMod(\Bun_G^{P^-\hgen})\to \mCI(G,P^-).  \]

By definition, we have
\[\Bun_G^{P^-\hgen}\mt_{\Bun_G} \widetilde{\Bun}_{R^-} \simeq \mathbf{Maps}_{\mathrm{gen}}(X, G\backslash \overline{G/U_R^-}/L \gets P^-\backslash G/R^- ) . \]
Similarly to Construction \ref{constr-filtration-geo}, we can define a filtration on this fiber product. This filtration is still labelled by elements $\overline{w}\in W_G/W_R$ because the definition of the poset $W_G/W_R$ is invariant under the Cartan involution.

\begin{constr} For any $\overline{w}\in W_G/W_R$ and any representative $w\in W_G$ of it, consider the closure $\overline{P^-wR^-}$ of $P^-wR^-$ in $G$. Note that it only depends on $\overline{w}$. Define
\[ (\Bun_G^{P^-\hgen}\mt_{\Bun_G} \widetilde{\Bun}_{R^-})^{\le \overline{w}} := \mathbf{Maps}_{\mathrm{gen}}(X,G\backslash \overline{G/U_R^-}/L \gets P^-\backslash \overline{P^-wR^-}/R^- )
\]
Then we obtain a filtration of $\Bun_G^{P^-\hgen}\mt_{\Bun_G} \widetilde{\Bun}_{R^-} $ by the poset $ W_G/W_R$. We also consider
\[ (\Bun_G^{P^-\hgen}\mt_{\Bun_G} \widetilde{\Bun}_{R^-})^{= \overline{w}} := \mathbf{Maps}_{\mathrm{gen}}(X,G\backslash \overline{G/U_R^-}/L \gets P^-\backslash P^-wR^-/R^- )
\]
As in Remark \ref{rem-geometric-stratum-as-Bungen}, we have
\begin{equation} \label{eqn-geometric-stratum-as-Bungen-dual}
  (\Bun_G^{P^-\hgen}\mt_{\Bun_G} \widetilde{\Bun}_{R^-})^{= \overline{w}}\simeq \widetilde{\Bun}_{(R^-)^w}^{P^-\cap (R^-)^w\hgen}, 
 \end{equation}
where $(R^-)^w:= wR^- w^{-1}$, which is also the image of $R^w$ under the Cartan involution.
\end{constr}

\begin{constr} \label{constr-description-gr-dual-weyl}
As in Construction \ref{constr-weyl-filtration-alpha} and Construction \ref{constr-gr-weyl}, we obtain a filtration on $\epsilon_{P^-}$ such that 
\[
\begin{aligned}
 \epsilon_{P^-}^{\le \overline{w}}: \DMod(\Bun_L) \xrightarrow{*\hpull} \DMod(\Bun_{R^-}) \xrightarrow{!\hpush} \DMod( \widetilde{\Bun}_{R^-} ) \to\\ \xrightarrow{!\hpull} \DMod((\Bun_G^{P^-\hgen}\mt_{\Bun_G} \widetilde{\Bun}_{R^-})^{\le \overline{w}}) \xrightarrow{\mathrm{push}} \DMod(\Bun_G^{P^-\hgen}).
\end{aligned}
 \]
 We also have a similar description for $\epsilon_{P^-}^{= \overline{w}}$. Namely, 
 \begin{itemize}
	\item If $P\in \Par_G^{\overline{w}}$, then 
	\[
	\begin{aligned}
	\epsilon_{P^-}^{= \overline{w}}: \DMod(\Bun_L) \xrightarrow{*\hpull} \DMod(\Bun_{R^-}) \xrightarrow{!\hpush} \DMod( \widetilde{\Bun}_{R^-} ) \to\\ \xrightarrow{!\hpull} \DMod((\Bun_G^{P^-\hgen}\mt_{\Bun_G} \widetilde{\Bun}_{R^-})^{= \overline{w}}) \xrightarrow{*\hpush} \DMod(\Bun_G^{P^-\hgen}),
	\end{aligned}
	\]
  which by (\ref{eqn-geometric-stratum-as-Bungen-dual}) is equivalent to
  \[ 
  \begin{aligned}
  \DMod(\Bun_L) \simeq \DMod(\Bun_{L^w}) \xrightarrow{*\hpull} \DMod(\Bun_{(R^-)^w}) \xrightarrow{!\hpush} \DMod( \widetilde{\Bun}_{(R^-)^w} ) \to\\ \xrightarrow{!\hpull} \DMod( \widetilde{\Bun}_{(R^-)^w}^{P^-\cap (R^-)^w\hgen}) \xrightarrow{*\hpush} \DMod(\Bun_G^{P^-\hgen}). 
  \end{aligned}
  \]
	\item If $P\notin \Par_G^{\overline{w}}$, then $\epsilon_{P^-}^{= \overline{w}}\simeq 0$.
\end{itemize} 

 Composing with $\Av_*^{U_P^-(\mBA)}$, we obtain a filtration on $\alpha_P^\vee$.
\end{constr}

In view of Remark \ref{rem-inv-inv-duality-functors-global}, the natural transformation (\ref{eqn-nt-alpha}) induces a natural transformation
\[ \Eis_{Q^-\to P^-}^\enh\circ \alpha_Q^\vee \simeq (\alpha_Q \circ \CT_{P\gets Q}^\enh)^\vee\to (\alpha_P)^\vee  \]
which, by construction, is given by 
\begin{equation} \label{eqn-nt-alpha-dual}
\begin{aligned}
\Eis_{Q^-\to P^-}^\enh\circ \alpha_Q^\vee \simeq \Eis_{Q^-\to P^-}^\enh\circ \CT_{G\gets Q^-}^\enh\circ \Eis_{R^-,!} \simeq \Eis_{Q^-\to P^-}^\enh\circ \CT_{P^-\gets Q^-}^\enh \circ \CT_{G\gets P^-}^\enh\circ \Eis_{R^-,!} \to \\
\to \CT_{G\gets P^-}^\enh\circ \Eis_{R^-,!} \simeq (\alpha_P)^\vee.
\end{aligned}
\end{equation}

\begin{constr} Similarly to Construction \ref{constr-weyl-filtration-functorial}, the natural transformation (\ref{eqn-nt-alpha-dual}) is compatible with the filtrations.

In more detail, consider the map $p_{Q^-\to P^-}^\enh:\Bun_G^{Q^-\hgen}\to \Bun_G^{P^-\hgen}$. By definition
\[
 \oblv^{U_{P}^-(\mBA)} \circ \Eis_{Q^-\to P^-}^\enh \simeq (p_{Q^-\to P^-}^\enh)_!\circ \oblv^{U_{Q}^-(\mBA)}.
\]
Hence we have a natural transformation
\begin{equation} \label{eqn-proof-surj-anti-temper-3-dual} 
\Eis_{Q^-\to P^-}^\enh\circ \Av_*^{U_{Q}^-(\mBA)} \to \Av_*^{U_{P}^-(\mBA)} \circ (p_{Q^-\to P^-}^\enh)_!
\end{equation}
between functors $\DMod(\Bun_G^{Q^-\hgen})\to \mCI(G,P^-)$.

On the other hand, the obvious commutative sqaure
\[
\xymatrix{
	(\Bun_G^{Q^-\hgen}\mt_{\Bun_G} \widetilde{\Bun}_{R^-})^{\le \overline{w}} \ar[r] \ar[d] &
	\Bun_G^{Q^-\hgen} \ar[d]\\
	(\Bun_G^{P^-\hgen}\mt_{\Bun_G} \widetilde{\Bun}_{R^-})^{\le \overline{w}}\ar[r] &
	\Bun_G^{P^-\hgen}
}
\]
induces a natural transformation 
\begin{equation} \label{eqn-proof-surj-anti-temper-2-dual}
  (p_{Q^-\to P^-}^\enh)_!\circ \epsilon_{Q^-}^{\le \overline{w}}\to \epsilon_{P^-}^{\le \overline{w}}.
\end{equation}
Combining with the previous paragraph, we obtain the desired natural transformation
\begin{equation} \label{eqn-proof-surj-anti-temper-4-dual}
\begin{aligned}
 \Eis_{Q^-\to P^-}^\enh\circ (\alpha_Q^\vee)^{\le \overline{w}} \simeq  \Eis_{Q^-\to P^-}^\enh\circ   \Av_*^{U_Q^-(\mBA)} \circ \epsilon_{Q^-}^{\le \overline{w}} \to \Av_*^{U_{P}^-(\mBA)} \circ (p_{Q^-\to P^-}^\enh)_! \circ \epsilon_{Q^-}^{\le \overline{w}} \to \\
 \to \Av_*^{U_{P}^-(\mBA)} \circ  \epsilon_{P^-}^{\le \overline{w}} \simeq (\alpha_P^\vee)^{\le \overline{w}}
 \end{aligned}
 \end{equation}
\end{constr}

\begin{constr} \label{constr-dual-weyl-filtration}
Passing back to the dual side, we obtain the \emph{dual Weyl filtration} $ \alpha_P'^{\le \overline{w}}$ on $\alpha_P$, where
\[  \alpha_P'^{\le \overline{w}}:= ((\alpha_P^\vee)^{\le \overline{w}})^\vee.\]
The natural transformation (\ref{eqn-proof-surj-anti-temper-4-dual}) induces 
\begin{equation} \label{eqn-proof-surj-anti-temper-5}
 \alpha_P'^{\le \overline{w}}\circ \CT_{P\gets Q}^\enh \to \alpha_Q'^{\le \overline{w}}.
 \end{equation}
\end{constr}

\begin{rem} We have not checked whether the dual Weyl filtration coincides with the Weyl filtration, i.e., $\alpha_P'^{\le \overline{w}}\simeq \alpha_P^{\le \overline{w}}$. Probably this can be proven by using the description of the unit object of the duality between $\mCI(G,P)$ and $\mCI(G,P^-)$ in \cite{chen2021thesis}.
\end{rem}

\subsection{\texorpdfstring{The equivalence $\Id\to \beta_G\circ \beta_G^L$}{Fully faithful of the G-anti-tempered de-gluing functor}: the cancellation lemma}
\label{ssec-surj-anti-temper-cancel}

In this subsection, we prove the cancellation lemma. The proof is very technical and we believe future techniques can simplify it.

\begin{prop}  [Cancellation lemma for the dual Weyl filtration] \label{prop-cancel-dual-weyl} 

For fixed $\overline{w}\in W_G/W_R$ and any representative $w$ of it, let $Q,P\in \Par_G^{\overline{w}}$ with $Q\subset P$ satisfy the following conditions:

\begin{itemize}
	\item[(i)] $Q\cap R^w=P\cap R^w$;
	\item[(ii)] $U_Q\cap M_P\subset U_R^w$, where $U_Q$ (resp. $U_R$) is the unipotent radical of $Q$ (resp. $R$), $M_P$ is the Levi subgroup of $P$ and $U_R^w:=wU_Rw^{-1}$.
\end{itemize}
Then (\ref{eqn-proof-surj-anti-temper-4}) induces an equivalence
\[\alpha_Q'^{= \overline{w}}\circ \CT_{P\gets Q}^\enh  \simeq \alpha_P'^{= \overline{w}}\]
between functors $\mCI(G,P)\to \DMod(\Bun_L)$.
\end{prop}

\begin{rem} In Corollary \ref{cor-condition-cancellation-lemma}, we will show Condition (ii) actually follows from Condition (i) and the assumption $Q\subset P\in \Par_G^{\overline{w}}$.
\end{rem}

\proof 

By duality, we need to show that the natural transformation
\[
\Eis_{Q^-\to P^-}^\enh\circ (\alpha_Q^\vee)^{= \overline{w}} \to (\alpha_P^\vee)^{= \overline{w}}
\]
is an equivalence.

Using the Cartan involution on $G$, Condition (i) and (ii) are equivalent to 
\begin{itemize}
	\item[(i-)] $Q^-\cap (R^-)^w=P^-\cap (R^-)^w$;
	\item[(ii-)] $U_Q^-\cap M_P\subset (U_R^-)^w$.
\end{itemize}
Hence, by (\ref{eqn-geometric-stratum-as-Bungen-dual}), we have
\[   (\Bun_G^{P^-\hgen}\mt_{\Bun_G} \widetilde{\Bun}_{R^-})^{= \overline{w}}\simeq  (\Bun_G^{Q^-\hgen}\mt_{\Bun_G} \widetilde{\Bun}_{R^-})^{= \overline{w}}.\]
It follows that (\ref{eqn-proof-surj-anti-temper-2-dual}) induces an equivalence
\[
  (p_{Q^-\to P^-}^\enh)_!\circ \epsilon_{Q^-}^{= \overline{w}}\simeq \epsilon_{P^-}^{= \overline{w}}.
\]
It remains to show the natural transformation
\[  \Eis_{Q^-\to P^-}^\enh\circ   \Av_*^{U_Q^-(\mBA)} \circ \epsilon_{Q^-}^{=\overline{w}} \to \Av_*^{U_{P}^-(\mBA)} \circ (p_{Q^-\to P^-}^\enh)_! \circ \epsilon_{Q^-}^{= \overline{w}}   \]
induced by (\ref{eqn-proof-surj-anti-temper-3-dual}) is an equivalence. To simplify the notations, we apply the Cartan involution and show that
\[  \Eis_{Q\to P}^\enh\circ   \Av_*^{U_Q(\mBA)} \circ \epsilon_{Q}^{=\overline{w}} \to \Av_*^{U_{P}(\mBA)} \circ (p_{Q\to P}^\enh)_! \circ \epsilon_{Q}^{= \overline{w}}   \]
is an equivalence.

Similarly to Remark \ref{rem-non-rigorous-proof}, our strategy is to show that the images of the functor $\epsilon_{Q}^{= \overline{w}}$ are ``$U_Q\cap M_P(\mBA)$-equivariant'' objects, and then show the natural transformation (\ref{eqn-proof-surj-anti-temper-3-dual}) sends any such object to an invertible morphism.

For this, we first need to define such ``$U_Q\cap M_P(\mBA)$-equivariant objects''.  Consider the functor
\begin{equation} \label{eqn-proof-cancellation-1}
   \DMod(\Bun_G^{Q\hgen}) \xrightarrow{!\hpull} \DMod( \Bun_Q ) \xrightarrow{*\hpush} \DMod( \Bun_{Q/U_P} ).
\end{equation}

\begin{defn}
Define $\DMod(\Bun_G^{Q\hgen})^{U_Q\cap M_P(\mBA)}\subset \DMod(\Bun_G^{Q\hgen}) $ to be the full subcategory sitting in the following Cartesian diagram
\[
\begin{tikzcd}
  \DMod(\Bun_G^{Q\hgen})^{U_Q\cap M_P(\mBA)}
    \arrow[r,"\subset"] \arrow[d]
    \arrow[dr, phantom, "\lrcorner", very near start]
    &\DMod(\Bun_G^{Q\hgen}) \arrow[d,"(\ref{eqn-proof-cancellation-1})"]  \\
      \DMod(\Bun_{M_Q})\arrow[r,"\subset","*\hpull"']
    & \DMod(  \Bun_{Q/U_P}  ).
  \end{tikzcd}
\]
\end{defn}

\begin{rem} \label{rem-image-BunM-in-BunP-limit}
The essential image of the functor $ \DMod(\Bun_{M_Q})\to \DMod(  \Bun_{Q/U_P}  )$ is stable under limits. To see this, we can replace the $*$-pullback functor by the $!$-pullback functor and show the latter has a left adjoint. By the contraction principle (see e.g. \cite[Appendix C]{drinfeld2015compact}), the left adjoint, i.e., the $!$-pushforward functor along $\Bun_{Q/U_P}\to \Bun_{M_Q}$ is indeed well-defined and given by the $!$-pullback functor along $\Bun_{M_Q} \to \Bun_{Q/U_P} $.
\end{rem}

Using the base-change isomorphisms, it is easy to see that $ \mCI(G,Q)\subset \DMod(\Bun_G^{Q\hgen})^{U_Q\cap M_P(\mBA)}$.

Then to finish the proof of the cancellation lemma (Porposition \ref{prop-cancel-dual-weyl}), it is enough to prove the following two lemmas:

\begin{lem} \label{lem-proof-cancellation-2} For any arrow $Q\subset P$ in $\Par^\st_G$, the natural transformation
\[
\xymatrix{
	\DMod(\Bun_G^{Q\hgen})  \ar[r]^-{\Av_*^{U_{Q}(\mBA)}} \ar[d]_-{(p_{Q\to P}^\enh)_!} &
	\mCI(G,Q) \ar[d]^-{\Eis_{Q\to P}^\enh} \ar@{=>}[ld] \\
	\DMod(\Bun_G^{P\hgen}) \ar[r]_-{\Av_*^{U_{P}(\mBA)}} &
	\mCI(G,P)
}
\]
is invertible when restricted to $\DMod(\Bun_G^{Q\hgen})^{U_Q\cap M_P(\mBA)}$.
\end{lem}

\begin{lem} \label{lem-proof-cancellation-1} Let $\overline{w}\in W_G/W_R$ and $Q\subset P$ be an arrow in  $\Par_G^{\overline{w}}$. If $U_Q\cap M_P\subset U_R^w$, then the functor
 \[ 
  \begin{aligned}
  \epsilon_{Q}^{= \overline{w}}: \DMod(\Bun_L) \simeq \DMod(\Bun_{L^w}) \xrightarrow{*\hpull} \DMod(\Bun_{R^w}) \xrightarrow{!\hpush} \DMod( \widetilde{\Bun}_{R^w} ) \to\\ \xrightarrow{!\hpull} \DMod( \widetilde{\Bun}_{R^w}^{Q\cap R^w\hgen}) \xrightarrow{*\hpush} \DMod(\Bun_G^{Q\hgen}). 
  \end{aligned}
  \]
factors through $\DMod(\Bun_G^{Q\hgen})^{U_Q\cap M_P(\mBA)}$.
\end{lem}

\proof[Proof of Lemma \ref{lem-proof-cancellation-2}.] 
It suffices to prove the assertion after composing the natural transformation in question with the forgetful functor 
\[\iota_{M_P}^!: \mCI(G,P)\to \DMod(\Bun_{M_P}).\] 
Let $\mCF\in \DMod(\Bun_G^{Q\hgen})^{U_Q\cap M_P(\mBA)}$. We have a canonical arrow $\Av_*^{U_{Q}(\mBA)}(\mCF)\to \mCF $, where, abusing notation, we view both $\mCF$ and $\Av_*^{U_{Q}(\mBA)}(\mCF)$ as objects in $\DMod(\Bun_G^{Q\hgen})$.

Now note that the functor $\iota_{M_P}^!\circ \Av_*^{U_{P}(\mBA)}$ is equivalent to
\[ \DMod(\Bun_G^{P\hgen}) \xrightarrow{!\hpull} \DMod(\Bun_P) \xrightarrow{*\hpush}  \DMod(\Bun_{M_P}) \]
because they have the same left adjoint (Lemma \ref{lem-generator-IGP}). Then, unwinding the definitions, the assertion of the lemma boils down to proving that the functor
\[ \DMod(\Bun_G^{Q\hgen}) \xrightarrow{!\hpush} \DMod(\Bun_G^{P\hgen}) \xrightarrow{!\hpull} \DMod(\Bun_P) \xrightarrow{*\hpush}  \DMod(\Bun_{M_P})  \]
sends $\Av_*^{U_{Q}(\mBA)}(\mCF)\to \mCF $ to an invertible morphism.

To this end, observe that
\[ \Bun_G^{Q\hgen}\mt_{\Bun_G^{P\hgen}} \Bun_P \simeq \Bun_P^{Q\hgen};\]
thus, using proper base-change and Remark \ref{rem-!-push-gen-exist}, the above composition is equivalent to
\[ \DMod(\Bun_G^{Q\hgen}) \xrightarrow{!\hpull} \DMod( \Bun_P^{Q\hgen}) \xrightarrow{!\hpush} \DMod(\Bun_P) \xrightarrow{*\hpush}  \DMod(\Bun_{M_P}) .\]
Since $\Bun_P^{Q\hgen}\to \Bun_P$ and $\Bun_{M_P}^{Q/U_P\hgen}\to \Bun_{M_P}$ are pseudo-proper, this is equivalent to
\[ \DMod(\Bun_G^{Q\hgen}) \xrightarrow{!\hpull} \DMod( \Bun_P^{Q\hgen}) \xrightarrow{*\hpush} \DMod(\Bun_{M_P}^{Q/U_P\hgen}) \xrightarrow{!\hpush}  \DMod(\Bun_{M_P}) .\]
Then it suffices to prove that
\[\theta: \DMod(\Bun_G^{Q\hgen}) \xrightarrow{!\hpull} \DMod( \Bun_P^{Q\hgen}) \xrightarrow{*\hpush} \DMod(\Bun_{M_P}^{Q/U_P\hgen}) \]
sends $\Av_*^{U_{Q}(\mBA)}(\mCF)\to \mCF $ to an invertible morphism.

Using base-change, it is easy to check that $\theta$ sends $\DMod(\Bun_G^{Q\hgen})^{U_Q\cap M_P(\mBA)}$ into $\mCI(P,Q)\subset  \DMod(\Bun_{M_P}^{Q/U_P\hgen})$. Hence we just need to show that, for any $\mCK\in \mCI(P,Q)$, the morphism
\[ \Maps( \mCK, \theta( \Av_*^{U_{Q}(\mBA)}(\mCF)) ) \to \Maps(\mCK,\theta(\mCF) ) \]
is invertible. As shown in the proof of Proposition-Construction \ref{propconstr-IGP-functorial}, the left adjoint $\theta^L$ of $\theta$ sends $\mCI(P,Q)$ to $\mCI(G,Q)$, and, as such, it is just the functor $\Eis:\mCI(P,Q)\to \mCI(G,Q)$. Hence we just need to check that the natural arrow
\[ \Maps( \Eis(\mCK), \Av_*^{U_{Q}(\mBA)}(\mCF) ) \to \Maps(\Eis(\mCK),\mCF ) \]
is an isomorphism: this is obvious because $\Eis(\mCK)$ is contained in $\mCI(G,Q)$.

\qed[Lemma \ref{lem-proof-cancellation-2}] 

To prove Lemma \ref{lem-proof-cancellation-1}, we first need to perform some preliminary simplifications. By definition, we need to show that
\[
\begin{aligned}
 \DMod(\Bun_{L^w}) \xrightarrow{*\hpull} \DMod(\Bun_{R^w}) \xrightarrow{!\hpush} \DMod( \widetilde{\Bun}_{R^w} ) \xrightarrow{!\hpull} \DMod( \widetilde{\Bun}_{R^w}^{Q\cap R^w\hgen}) \to \\ \xrightarrow{*\hpush} \DMod(\Bun_G^{Q\hgen}) \xrightarrow{!\hpull} \DMod(\Bun_Q) \xrightarrow{*\hpush} \DMod(\Bun_{Q/U_P})
 \end{aligned}
\]
factors through the fully faithful pullback functor $\DMod(\Bun_{M_Q})\to \DMod(\Bun_{Q/U_P})$. We write
\[ \widetilde{\Bun}_{R^w}^{Q\cap R^w\hgen}\mt_{ \Bun_G^{Q\hgen} }\Bun_Q =: ( \widetilde{\Bun}_{R^w}\mt_{\Bun_G} \Bun_Q)^{Q\cap R^w\hgen},  \]
and similarly when $\widetilde{\Bun}_{R^w}$ is replaced by $\Bun_{R^w}$. By the base-change isomorphism, the above displayed functor is equivalent to
\begin{equation} \label{eqn-proof-lem-proof-cancellation-1-1}
\begin{aligned}
 \DMod(\Bun_{L^w}) \xrightarrow{*\hpull} \DMod(\Bun_{R^w}) \xrightarrow{!\hpush} \DMod( \widetilde{\Bun}_{R^w} )\to \\ \xrightarrow{!\hpull} \DMod( ( \widetilde{\Bun}_{R^w}\mt_{\Bun_G} \Bun_Q)^{Q\cap R^w\hgen})   \xrightarrow{*\hpush} \DMod(\Bun_{Q/U_P}).
 \end{aligned}
\end{equation}
Now we use the defect stratification on $\widetilde{\Bun}_{R^w}$ to reduce to a similar statement about $\Bun_{R^w}$. In other words, we prove the following result.

\begin{lem} \label{lem-proof-cancellation-3} For $\overline{w}\in W_G/W_R$, let $Q\subset P$ be an arrow in $\Par_G^{\overline{w}}$ such that $U_Q\cap M_P\subset U_R^w$. Assume that the functor
 \begin{equation}\label{eqn-proof-lem-proof-cancellation-3-0}
  \begin{aligned}
   \DMod(\Bun_{L^w}) \xrightarrow{!\hpull}\DMod((\Bun_{R^w}\mt_{\Bun_G} \Bun_Q)^{Q\cap R^w\hgen}) \xrightarrow{*\hpush} \DMod(\Bun_{Q/U_P})
  \end{aligned}
  \end{equation}
factors through the pullback functor $\DMod(\Bun_{M_Q})\to \DMod(\Bun_{Q/U_P})$. Then so does the functor (\ref{eqn-proof-lem-proof-cancellation-1-1}).
\end{lem}

\proof Recall that $\widetilde{\Bun}_{R^w}$ has a defect stratification such that the disjoint union of the strata is of the form $\Bun_{R^w}\mt_{\Bun_{L^w}} H$, where $H$ is certain Hecke ind-stack for $L^w$-torsors. Let 
\[\Bun_{R^w}\mt_{\Bun_{L^w}} H\to \widetilde{\Bun}_{R^w} \]
be the locally closed embedding for this stratification. Then the composition
\[ \Bun_{R^w}\mt_{\Bun_{L^w}} H\to \widetilde{\Bun}_{R^w} \to \Bun_G \]
is given by
\[ \Bun_{R^w}\mt_{\Bun_{L^w}} H\xrightarrow{\mathrm{pr}_1 } \Bun_{R^w}  \to \Bun_G. \]

By the same method as in \cite[Theorem 6.2.10]{braverman2002geometric}, the image of the composition
\[\DMod(\Bun_{L^w}) \xrightarrow{*\hpull} \DMod(\Bun_{R^w}) \xrightarrow{!\hpush} \DMod( \widetilde{\Bun}_{R^w} ) \xrightarrow{!\hpull} \DMod(\Bun_{R^w}\mt_{\Bun_{L^w}} H)   \]
is contained in the image of the pullback functor $\DMod(H)\to \DMod(\Bun_{R^w}\mt_{\Bun_{L^w}} H)$. Let $\mCF$ be any object in $ \DMod( \widetilde{\Bun}_{R^w} )$ such that its $!$-restriction to $\Bun_{R^w}\mt_{\Bun_{L^w}} H$ is an object pulled back from $H$. Then we only need to show the functor
\[
 \phi: \DMod( \widetilde{\Bun}_{R^w} ) \xrightarrow{!\hpull} \DMod( ( \widetilde{\Bun}_{R^w}\mt_{\Bun_G} \Bun_Q)^{Q\cap R^w\hgen})   \xrightarrow{*\hpush} \DMod(\Bun_{Q/U_P}) 
\]
sends $\mCF$ to an object pulled back from $\Bun_{M_Q}$. 

We can assume $\mCF$ is supported on a single connected component $\widetilde{\Bun}_{R^w,\lambda}$ of $\widetilde{\Bun}_{R^w}$. Write this connected component as a union $\cup U_\alpha$ of open substacks such that each $U_\alpha$ is the union of \emph{finitely} many defect strata. Then $\mCF \simeq \lim_\alpha \mCF_\alpha$, where $\mCF_\alpha$ is the $*$-extension of $\mCF|_{U_\alpha}$. 

We claim each $\phi(\mCF_\alpha)$ is an object pulled back from $\Bun_{M_Q}$. To prove the claim, we can replace $\mCF_\alpha$ by the $*$-extension of $\mCF|^!_{Y}$, where $Y$ is some defect stratum of $\widetilde{\Bun}_{R^w}$. By definition, such $*$-extension is contained in the image of the functor
\[
\varphi:  \DMod(H) \xrightarrow{!\hpull} \DMod( \Bun_{R^w}\mt_{\Bun_{L^w}} H) \xrightarrow{*\hpush} \DMod(\widetilde{\Bun}_{R^w}) .
 \]
By the base-change isomorphism, the composition $\phi\circ \varphi$ is equivalent to
\[  \DMod(H) \xrightarrow{*\hpush} \DMod(\Bun_{L^w}) \xrightarrow{(\ref{eqn-proof-lem-proof-cancellation-3-0})} \DMod(\Bun_{Q/U_P}) .  \]
Hence the assumption on the functor (\ref{eqn-proof-lem-proof-cancellation-3-0}) implies the claim.

By the above claim and Remark \ref{rem-image-BunM-in-BunP-limit}, the object $\lim_\alpha \phi(\mCF_\alpha)$ is pulled back from $\Bun_{M_Q}$.
It remains to show $\phi(\lim_\alpha \mCF_\alpha)\to \lim_\alpha \phi(\mCF_\alpha)$ is an isomorphism. Recall that each connected component of $\Bun_Q$ is quasi-compact over $\Bun_G$, hence each connected component of 
\[ ( \widetilde{\Bun}_{R^w}\mt_{\Bun_G} \Bun_Q)^{Q\cap R^w\hgen} \]
is quasi-compact over $\widetilde{\Bun}_{R^w}$. Thus, for a fixed such connected component, the inverse images of $U_\alpha\subset \widetilde{\Bun}_{R^w}$ in it stabilize for large enough $U_\alpha$. This makes the desired claim manifest.

\qed[Lemma \ref{lem-proof-cancellation-3}]

Next, let us prove that the assumption of Lemma \ref{lem-proof-cancellation-3} holds true. This will conclude the proof of Lemma \ref{lem-proof-cancellation-1}. By definition, we have
\[ (\Bun_{R^w}\mt_{\Bun_G} \Bun_Q)^{Q\cap R^w\hgen} \simeq \mathbf{Maps}_{\mathrm{gen}}(X, Q\backslash G/R^w \gets \mBB (Q\cap R^w) ), \]
which, up to non-reduced structures, it is isomorphic to
\[ \mathbf{Maps}_{\mathrm{gen}}(X, Q\backslash \overline{Q R^w}/R^w \gets \mBB (Q\cap R^w) ). \]
For any finite collection $\underline{x}=\big\{x_1,\dots,x_n \big\}$ of closed points of $X$, let
\[ (\Bun_{R^w}\mt_{\Bun_G} \Bun_Q)^{Q\cap R^w\textnormal{ near } \underline{x}} \subset \mathbf{Maps}_{\mathrm{gen}}(X, Q\backslash \overline{Q R^w}/R^w \gets \mBB (Q\cap R^w) )  \]
be the open substack classifying those maps $X\to Q\backslash \overline{Q R^w}/R^w  $ that send the points $x_i$, and therefore a Zariski neighborhood of them, into $ \mBB (Q\cap R^w)\subset Q\backslash \overline{Q R^w}/R^w$. By definition, when $\underline{x}$ varies, these open substacks form a cover.

The following result says we can replace Lemma \ref{lem-proof-cancellation-1} to a similar statement concerning these open substacks.

\begin{lem} \label{lem-proof-cancellation-4} Let $\overline{w}\in W_G/W_R$ and $Q\subset P$ be an arrow in $\Par_G^{\overline{w}}$ such that $U_Q\cap M_P\subset U_R^w$. For any finite collection $\underline{x}=\big\{x_1,\dots,x_n \big\}$ of closed points of $X$, suppose that the functor
 \begin{equation}\label{eqn-proof-lem-proof-cancellation-4-0}
  \begin{aligned}
   \DMod(\Bun_{L^w}) \xrightarrow{!\hpull}\DMod((\Bun_{R^w}\mt_{\Bun_G} \Bun_Q)^{Q\cap R^w\textnormal{ near } \underline{x}}) \xrightarrow{*\hpush} \DMod(\Bun_{Q/U_P}). 
  \end{aligned}
  \end{equation}
factors through the pullback functor $\DMod(\Bun_{M_Q})\to \DMod(\Bun_{Q/U_P})$. Then so does the functor (\ref{eqn-proof-lem-proof-cancellation-3-0}).

\end{lem}

\proof The lemma does \emph{not} follow immediately from Zariski descent and Remark \ref{rem-image-BunM-in-BunP-limit} because 
\[(\Bun_{R^w}\mt_{\Bun_G} \Bun_Q)^{Q\cap R^w\hgen}\] can not be covered by finitely many
\[(\Bun_{R^w}\mt_{\Bun_G} \Bun_Q)^{Q\cap R^w\textnormal{ near } \underline{x}}.\]
However, this problem is not essential because 
\begin{equation} \label{eqn-proof-lem-proof-cancellation-4-1}
 (\Bun_{R^w}\mt_{\Bun_G} \Bun_Q)^{Q\cap R^w\hgen}\to \Bun_{Q/U_P} 
 \end{equation}
is quasi-compact when restricted to each connected component of the source. 

In more detail, recall $\Bun_{R^w} \to \Bun_G$ is quasi-compact when restricted to each connected component of the source, and $\Bun_{Q} \to \Bun_{Q/U_P}$ is quasi-compact. Hence the last claim of the above paragraph is true. To prove the lemma, we can replace $\Bun_{M_Q}$ by any fixed quasi-compact open substack $V \subset \Bun_{M_Q}$ and replace the map (\ref{eqn-proof-lem-proof-cancellation-4-1}) by
\[  (\Bun_{R^w}\mt_{\Bun_G} \Bun_Q)^{Q\cap R^w\hgen}\mt_{ \Bun_{M_Q} } V\to \Bun_{Q/U_P} \mt_{\Bun_{M_Q}} V. \]
Recall that $ \Bun_{Q/U_P}$ is quasi-compact over $\Bun_{M_Q}$. Hence the target of the above map is quasi-compact and so is each connected component of the source. Now each such connected component can be covered by finitely many 
\[ (\Bun_{R^w}\mt_{\Bun_G} \Bun_Q)^{Q\cap R^w\textnormal{ near } \underline{x}}\mt_{ \Bun_{M_Q} } V,\]
and we are done by Zariski descent and Remark \ref{rem-image-BunM-in-BunP-limit}.

\qed[Lemma \ref{lem-proof-cancellation-4}]

\proof[Proof of Lemma \ref{lem-proof-cancellation-1}.] We just need to show that the assumption of Lemma \ref{lem-proof-cancellation-4} is valid. Write 
\[  1\to K_1\to G_1\to H_1 \to 1,\;\; 1\to K_2\to G_2\to H_2 \to 1 \]
for the short exact sequences
\[   1\to (Q\cap U_R^w)\to  Q\cap R^w \to (Q\cap R^w)/(Q\cap U_R^w) \to 1,\;\; 1\to U_Q/U_P \to Q/U_P \to Q/U_Q \to 1. \]
We claim both sequences have canonical splittings. For the second sequence, the splitting is provided by the isomorphism
\[ M_Q \to Q \to Q/U_P \to Q/U_Q. \]
For the first one, we only need to show that
\[ Q\cap L^w \to Q\cap R^w  \to (Q\cap R^w)/(Q\cap U_R^w) \]
is an isomorphism. Note that both $Q$ and $L^w$ contain the Cartan subgroup $T$. Hence, the weight decomposition of $\Lie(G)$ implies the above map induces an isomorphism of Lie algebras. On the other hand, choose $\mBG_m\to T$ such that the adjoint $\mBG_m$-action on $R^w$ contracts it onto $L^w$. The restricted action on the closed subscheme $Q\cap R^w$ must contract it onto $Q\cap L^w$. In particular, $Q\cap L^w$ has non-empty intersection with any connected component of $Q\cap R^w$. This finishes the proof of the claim.

Write $K_0:= U_Q\cap M_P$. By assumption, $K_0\subset K_1$, and the composition $K_0\to G_1\to G_2$ induces an isomorphism $K_0 \to K_2$. Note that $K_0,K_1,K_2$ are unipotent.

Consider the fiber product
\[\Hecke_{G_2,\underline{x}} \mt_{ \Hecke_{H_2,\underline{x}} } \Bun_{H_2} \]
studied in the proof of Lemma \ref{lem-sphP-preserve-BunM}. Informally, it classifies Hecke modifications of $G_2$-torsors at $\underline{x}$ that induce trivial Hecke modifications of the induced $H_2$-torsors. As in the proof of Lemma \ref{lem-sphP-preserve-BunM}, an object $\mCF\in \DMod(\Bun_{G_2})$ is pulled back from $\Bun_{H_2}$ iff the pullbacks of $\mCF$ along the two maps
\[  \Bun_{G_2} \xleftarrow{h^l} \Hecke_{G_2,\underline{x}} \mt_{ \Hecke_{H_2,\underline{x}} } \Bun_{H_2}\xrightarrow{h^r} \Bun_{G_2} \]
are isomorphic.

Now we reformulate the above condition using loop group actions. Consider the $G_2(\mCO_{\underline{x}})$-torsor $Z_2$ on $\Bun_{G_2}$ that classifies $G_2$-torsors equipped with a trivialization on $\Spec \mCO_{\underline{x}}$. The standard re-gluing construction entends the $G_2(\mCO_{\underline{x}})$-action on $Z_2$ to a $G_2(\mCK_{\underline{x}})$-action. Note that we have
\[\DMod(\Bun_{G_2}) \simeq \DMod(Z_2)^{G_2(\mCO_{\underline{x}}) }.\]
It follows the above condition is equivalent to $\mCF \in  \DMod(Z_2)^{ K_2(\mCK_{\underline{x}}) H_2(\mCO_{\underline{x}}) }$.

Similarly, we can define a $G_1(\mCO_{\underline{x}})$-torsor $Z_1$ on $(\Bun_{R^w}\mt_{\Bun_G} \Bun_Q)^{Q\cap R^w\textnormal{ near } \underline{x}}$ and extend the $G_1(\mCO_{\underline{x}})$-action on $Z_1$ to a $G_1(\mCK_{\underline{x}})$-action. Note that the $K_1(\mCK_{\underline{x}}) $-action on the source of
\[  Z_1\to  Z_1/ G_1(\mCK_{\underline{x}}) \simeq (\Bun_{R^w}\mt_{\Bun_G} \Bun_Q)^{Q\cap R^w\textnormal{ near } \underline{x}} \to \Bun_{L^w} \]
stabilizes the fibers (because $K_1$ is contained in the kernel of $Q\cap R^w \to L^w$). Hence the image of the $!$-pullback functor
\[\DMod(\Bun_{L^w}) \to \DMod((\Bun_{R^w}\mt_{\Bun_G} \Bun_Q)^{Q\cap R^w\textnormal{ near } \underline{x}}) \simeq   \DMod(Z_1)^{G_1(\mCO_{\underline{x}})}\]
is contained in the full subcategory
\[ \DMod(Z_1)^{ K_1(\mCK_{\underline{x}}) H_1(\mCO_{\underline{x}}) } \subset \DMod(Z_1)^{G_1(\mCO_{\underline{x}})} .\]

Also note that there is a canonical map $Z_1\to Z_2$ intertwining the $G_i(\mCK_{\underline{x}})$-actions and compatible with the map
\[ (\Bun_{R^w}\mt_{\Bun_G} \Bun_Q)^{Q\cap R^w\textnormal{ near } \underline{x}} \to \Bun_{M_Q}. \]

It remains to show the functor
\begin{equation} \label{eqn-proof-lem-proof-cancellation-1-5}
\DMod(Z_1)^{G_1(\mCO_{\underline{x}})} \simeq \DMod(Z_1/G_1(\mCO_{\underline{x}})) \xrightarrow{*\hpush} \DMod(Z_2/G_2(\mCO_{\underline{x}})) \simeq \DMod(Z_2)^{G_2(\mCO_{\underline{x}})} \end{equation}
sends $K_1(\mCK_{\underline{x}}) H_1(\mCO_{\underline{x}})$-equivariant objects to $K_2(\mCK_{\underline{x}}) H_2(\mCO_{\underline{x}})$-equivariant objects. This does \emph{not} follow immediately from $K_1\supset K_0 \simeq  K_2$ because $K_0$ is not a normal subgroup of $G_1$ hence we cannot define the notion of $K_0(\mCK_{\underline{x}}) H_1(\mCO_{\underline{x}})$-equivariant objects. 

Nevertheless, this can still be formally proved as follows.

We first replace the relevant ind-group schemes by group schemes. This step is standard. Namely, write $K_2(\mCK_{\underline{x}}) \simeq \bigcup_\alpha V_2^\alpha$ as union of group schemes stabilized by the adjoint action of $H_2(\mCO_{\underline{x}})$. For example, $V_2^\alpha$ can be chosen to be $h K_2(\mCO_{\underline{x}}) h^{-1}$ for certain central elements $h$ in $H_2(\mCO_{\underline{x}})$. We can assume that $V_2^\alpha$ contains $K_2(\mCO_{\underline{x}})$. For a fixed such $V_2:=V_2^\alpha$, let $V_0$ be the corresponding subgroup of $K_0(\mCK_x)$. Let $V_1\subset K_1(\mCK_{\underline{x}})$ be a sub-group scheme containing $V_0$ and stabilized by $H_1(\mCO_{\underline{x}})$. We only need to show that the functor (\ref{eqn-proof-lem-proof-cancellation-1-5}) sends $V_1 H_1(\mCO_{\underline{x}})$-equivariant objects to $V_2 H_2(\mCO_{\underline{x}})$-equivariant objects.

As in the proof of Lemma \ref{lem-proof-cancellation-4}, we replace $\Bun_{M_Q}$ by a fixed quasi-compact open substack and therefore assume $Z_i/G_i(\mCO_{\underline{x}})$ are quasi-compact algebraic stacks. In particular, we have Verdier self-duality on them. Passing to duality, we just need to show that the functor
\[ 
\begin{aligned}
\DMod(Z_2)^{G_2(\mCO_{\underline{x}})} \simeq  \DMod(Z_2/G_2(\mCO_{\underline{x}})) \xrightarrow{!\hpull}  \DMod(Z_1/G_1(\mCO_{\underline{x}})) \simeq \\
\simeq \DMod(Z_1)^{G_1(\mCO_{\underline{x}})} \xrightarrow{\Av_*} \DMod(Z_1)^{V_1 H_1(\mCO_{\underline{x}})}
\end{aligned}
\]
sends the kernel of
\[ \Av_*:  \DMod(Z_2)^{G_2(\mCO_{\underline{x}})} \to  \DMod(Z_2)^{V_2 H_2(\mCO_{\underline{x}})} \]
to $0$. Let $\mCC$ be this kernel. We have a commutative square
\[
\xymatrix{
    \DMod(Z_1)^{G_1(\mCO_{\underline{x}})} \ar[r]^-{\Av_*} \ar[d]^-{\oblv} & \DMod(Z_1)^{V_1 H_1(\mCO_{\underline{x}})} \ar[d]^-{\oblv} \\
    \DMod(Z_1)^{K_1(\mCO_{\underline{x}}) } \ar[r]^-{\Av_*} & \DMod(Z_1)^{V_1 },
}
\]
whose vertical arrows are conservative. Hence we just need to check that
\[ 
\begin{aligned}
\DMod(Z_2)^{G_2(\mCO_{\underline{x}})} \simeq  \DMod(Z_2/G_2(\mCO_{\underline{x}})) \xrightarrow{!\hpull}  \DMod(Z_1/K_1(\mCO_{\underline{x}})) \simeq \\
\simeq \DMod(Z_1)^{K_1(\mCO_{\underline{x}})} \xrightarrow{\Av_*} \DMod(Z_1)^{V_1}
\end{aligned}
\]
sends $\mCC$ to $0$. Since $ \DMod(Z_1)^{K_1(\mCO_{\underline{x}})} \subset  \DMod(Z_1)^{K_0(\mCO_{\underline{x}})}$, the above functor is equivalent to 
\[ 
\begin{aligned}
\DMod(Z_2)^{G_2(\mCO_{\underline{x}})} \simeq  \DMod(Z_2/G_2(\mCO_{\underline{x}})) \xrightarrow{!\hpull}  \DMod(Z_1/K_0(\mCO_{\underline{x}})) \simeq \\
\simeq \DMod(Z_1)^{K_0(\mCO_{\underline{x}})}  \xrightarrow{\Av_*} \DMod(Z_1)^{V_1},
\end{aligned}
\]
which in turn is equivalent to
\[ 
\begin{aligned}
\DMod(Z_2)^{G_2(\mCO_{\underline{x}})} \xrightarrow{\oblv} \DMod(Z_2)^{K_2(\mCO_{\underline{x}})} \xrightarrow{!\hpull}  \DMod(Z_1)^{K_0(\mCO_{\underline{x}})}  \xrightarrow{\Av_*} \DMod(Z_1)^{V_0} \xrightarrow{\Av_*} \DMod(Z_1)^{V_1}.
\end{aligned}
\]
Thus, it suffices to prove that
\[\DMod(Z_2)^{G_2(\mCO_{\underline{x}})} \xrightarrow{\oblv} \DMod(Z_2)^{K_2(\mCO_{\underline{x}})} \xrightarrow{!\hpull}  \DMod(Z_1)^{K_0(\mCO_{\underline{x}})}  \xrightarrow{\Av_*} \DMod(Z_1)^{V_0}\]
sends $\mCC$ to $0$. Note that the $!$-pullback functor $\DMod(Z_2)\to \DMod(Z_1)$ is compatible with the $K_2(\mCO_{\underline{x}})\simeq K_0(\mCO_{\underline{x}})$-actions. By \cite[Corollary 2.17.10]{raskin2016chiral}, the above composition is equivalent to
\[\DMod(Z_2)^{G_2(\mCO_{\underline{x}})} \xrightarrow{\oblv} \DMod(Z_2)^{K_2(\mCO_{\underline{x}})} \xrightarrow{\Av_*} \DMod(Z_2)^{V_2} \xrightarrow{!\hpull}  \DMod(Z_1)^{V_0}.\]
Then we are done because the composition of the first two functors is equivalent to
\[ \DMod(Z_2)^{G_2(\mCO_{\underline{x}})} \xrightarrow{\Av_*} \DMod(Z_2)^{V_2 H_2(\mCO_{\underline{x}})} \xrightarrow{\oblv} \DMod(Z_2)^{V_2}  .\]

\qed[Lemma \ref{lem-proof-cancellation-1}]

\qed[Proposition \ref{prop-cancel-dual-weyl}]

\subsection{\texorpdfstring{The equivalence $\Id\to \beta_G\circ \beta_G^L$}{Fully faithful of the G-anti-tempered de-gluing functor}: finishing the proof}
\label{ssec-surj-anti-temper-combo}

In this subsection, we use the cancellation lemma (Proposition \ref{prop-cancel-dual-weyl}) to finish the proof.

We first give the following convenient description of $\Par_G^{\overline{w}}$.

\begin{lem} \label{lem-nonzero-cond-gr}
For $P\in \Par_G^\st\setminus\{ G\}$, we have
\[P\in \Par_G^{\overline{w}}\;\;\;\text{iff}\;\;\;B\cap M_P \subset R^w,\] 
where $B$ is the fixed Borel subgroup, $M_P$ is the Levi subgroup of $P$, $w\in W_G$ is any representative of $\overline{w}$ and $R^w:=wRw^{-1}$.
\end{lem}

\proof By definition of the Bruhat order, $P\in \Par_G^{\overline{w}}$ iff $BwR$ is closed in $PwR$ iff $BR^w$ is closed in $PR^w$ iff $U_P\backslash BR^w$ is closed in $U_P\backslash PR^w$. 

Note that $U_P\backslash PR^w\simeq U_P\backslash (P\mt R^w)/(P\cap R^w) \simeq (M_P\mt R^w)/(P\cap R^w)$ and $U_P\backslash BR^w \simeq (B\cap M_P\mt R^w)/(B\cap R^w)$. Hence the above condition is equivalent to the map
\begin{equation} \label{eqn-lem-nonzero-cond-gr-1}
  (B\cap M_P\mt R^w)/(B\cap R^w)\to (M_P\mt R^w)/(P\cap R^w) 
 \end{equation}
being a closed embedding.

Note that $M_P R^w =(M_P\mt R^w)/(M_P\cap R^w)\to (M_P\mt R^w)/(P\cap R^w)$ is a $U_P\cap  R^w$-torsor. Hence (\ref{eqn-lem-nonzero-cond-gr-1}) is a closed embedding iff so is its base-change along this map. This base-change is $(B\cap M_P)R^w \to M_P R^w$. It is a closed embedding iff $(B\cap M_P)R^w/R^w \to M_P R^w/R^w$ is a closed embedding. This map can be rewritten as 
\begin{equation} \label{eqn-lem-nonzero-cond-gr-2}
  B\cap M_P/(B\cap R^w\cap M_P) \to M_P/(R^w\cap M_P).
\end{equation}
Note that $B\cap M_P$ is a Borel subgroup of $M_P$. We claim $R^w\cap M_P$ is a parabolic subgroup of $M_P$. To see this, first note that $R^w$ contains the Cartan subgroup $T$. We can choose a homomorphism $\mBG_m\to T$ such that the adjoint $\mBG_m$-action on $G$ has attractor equal to $R^w$. Then $R^w\cap M_P$ is the attractor of the restricted $\mBG_m$-action on $M_P$. This proves the claim.

Using the above claim, we see (\ref{eqn-lem-nonzero-cond-gr-2}) is a closed embedding iff $(B\cap M_P)(R^w\cap M_P)$ is closed in $M_P$ iff $B\cap M_P$ is contained in $R^w\cap M_P$ iff $B\cap M_P$ is contained in $R^w$. 

\qed[Lemma \ref{lem-nonzero-cond-gr}]

\begin{cor} \label{cor-condition-cancellation-lemma}
In the statement of Proposition \ref{prop-cancel-dual-weyl}, Condition (ii) follows from Condition (i). In other words, for $P,Q\in \Par_G^{\overline{w}}$ such that $P\cap R^w = Q\cap R^w$, we have $U_Q\cap M_P\subset U_R^w$. 
\end{cor}

\proof Note that $U_Q\cap M_P\subset B\cap M_P$. Hence by Lemma \ref{lem-nonzero-cond-gr}, $U_Q\cap M_P\subset R^w$. 

Since $U_Q\cap M_P$ and $U_R^w$ are connected, we only need to show that $\Lie(U_Q\cap M_P)\subset \Lie(U_R^w)$. Since $P,Q,R^w$ contain the Cartan subgroup $T$, by the weight decomposition of $\Lie(G)$, we just need to show that $\Lie(U_Q\cap M_P) \cap \Lie(L^w) = 0$. Applying the Cartan involution, this amounts to $\Lie(U_Q^-\cap M_P) \cap \Lie(L^w) = 0$. 

Note that $\Lie(U_Q^-\cap M_P) \subset \Lie(P)$ but $\Lie(U_Q^-\cap M_P)\cap \Lie(Q)=0$. Hence, by $P\cap R^w = Q\cap R^w$, we have $\Lie(U_Q^-\cap M_P)\cap \Lie(R^w) =0 $. In particular, $\Lie(U_Q^-\cap M_P) \cap \Lie(L^w) = 0$ as desired.

\qed[Corollary \ref{cor-condition-cancellation-lemma}]

Then we give a purely combinatorial sufficient condition for the assumption of the cancellation lemma. For convenience, we introduce the following definition:

\begin{defn} For any positive simple root $\check \alpha$ of $G$ and standard parabolics $Q\subset P \in  \Par_G^\st $, we write \emph{$P= Q\langle \check \alpha\rangle$} if the Dynkin diagram of $M_P$ is obtained from that of $M_Q$ by adding the vertex corresponding to $\check \alpha$. 
\end{defn}

\begin{lem} \label{lem-weyl-combinatorics} If $\overline{w} \neq \overline{w}_\circ \in W_G/W_R$, then there exists a negative simple root $-\check \alpha$ which is not a weight of $\Lie(R^w)$. Moreover, for $Q\subset P \in  \Par_G^\st\setminus \{ G \} $ such that $P=Q\langle \check \alpha\rangle$, we have
\begin{itemize}
  \item[(1)] $P\in  \Par_G^{\overline{w}}$ iff $Q\in  \Par_G^{\overline{w}}$;
  \item[(2)] If $P,Q\in \Par_G^{\overline{w}}$, then $P\cap R^w = Q\cap R^w$.
\end{itemize}
\end{lem}

\proof Suppose $\Lie(R^w)$ has all negative simple roots as weights. Then $N^-\subset R^w$ hence $R^w$ has transversal intersection with $B^-$, which contradicts $\overline{w} \neq \overline{w}_\circ$.

Let $-\check \alpha$ be a negative simple root which is not a weight of $\Lie(R^w)$. Then $\check \alpha$ must be a weight of $\Lie(R^w)$.

For (1), by Lemma \ref{lem-nonzero-cond-gr}, we only need to show that $B\cap M_P \subset R^w$ iff $B\cap M_Q \subset R^w$. The ``only if'' part is obvious. The ``if'' part follows from the above paragraph because $\Lie(B\cap M_P)$ is generated by $\Lie(B\cap M_Q)$ and the weight space of $\check \alpha$.

For (2), suppose $P\cap R^w \neq Q\cap R^w$. Then there exists a negative root $-\check\lambda$ which is a weight of $\Lie(P\cap R^w)$ but not of $\Lie(Q\cap R^w)$. Since $P=Q\langle \check \alpha\rangle$, we must have $-\check\lambda = -n\check \alpha-\check\mu$, where $n>0$ and $\check \mu$ is a non-negative linear sum of positive simple roots. We can assume there is no root strictly smaller than $\check\lambda$ satisfying the same assumption.

Note that $\check \mu \neq 0$, otherwise $-\check \alpha$ would be a weight of $\Lie(R^w)$ which contradicts our assumption. It follows that $\check \lambda$ is not simple: let us write $\check \lambda = \check \lambda_1+\check \lambda_2$ for two strictly positive roots $\check\lambda_1,\check\lambda_2$. Since $-\check \lambda$ is a weight of $\Lie(P)$, so are $\check\lambda_1$ and $-\check\lambda_1$.
In other words, $\pm \check\lambda_1$ are both weights of $\Lie(M_P)$. 
Hence, by the assumption $P\in \Par_G^{\overline{w}}$ and Lemma \ref{lem-nonzero-cond-gr}, $\check \lambda_1$ is a weight of $\Lie(R^w)$. Consequently, $-\check \lambda_2 = (-\check \lambda) + \check \lambda_1$ is also a weight of $\Lie(P\cap R^w)$. By the minimality assumption on $\check\lambda$, we see that $-\check \lambda_2$ must be a weight of $\Lie(Q\cap R^w)$. Write $\check \lambda_2$ as a sum of positive simple roots. Since $P=Q\langle \check \alpha\rangle$, this sum has no $\check \alpha$ term. By symmetry, the sum for $\check \lambda_1$ has no $\check \alpha$ term either. But $\check \lambda = \check \lambda_1+\check \lambda_2=n\check \alpha+\check \mu$ contains $\check \alpha$ terms. A contradiction!

\qed[Lemma \ref{lem-weyl-combinatorics}]

Now we are ready to finish the proof of Theorem \ref{thm-main-thm-atemp}.

\proof[Proof of Theorem \ref{thm-main-thm-atemp}.]

Recall our Goal \ref{goal-surj-anti-temper-colimit}: we need to show that
 \begin{equation} \label{eqn-goal-surj-anti-temper-colimit}
 \iota_L^! (\mCM_R)\to \colim_{P\in \Par_G^\st, P\neq G}  \alpha_P (\mCM_P)
\end{equation}
 is invertible, where $\alpha_P= \CT_{R,*}\circ \Eis_{P\to G}^\enh$, $\mCM_P\in \mCI(G,P)^{G\hatemp}$ and $\CT^\enh_{P\gets Q}(\mCM_Q) \simeq \mCM_P$ for $Q\subset P$.

Using the dual Weyl filtration on $\alpha_P$ (Construction \ref{constr-dual-weyl-filtration}), we obtain a filtration on the RHS by $W_G/W_R$, whose $\le \overline{w}$-term is
\[  \colim_{P\in \Par_G^\st, P\neq G}  \alpha_P'^{\le \overline{w}} (\mCM_P) .\]
To conclude, it suffices to prove the following two lemmas.

\begin{lem} \label{lem-colimit-cancel} We have:
\begin{itemize}
   \item[(1)] The morphism (\ref{eqn-goal-surj-anti-temper-colimit}) factors through an equivalence
  \[  \iota_L^! (\mCM_R)\to \colim_{P\in \Par_G^\st, P\neq G}  \alpha_P'^{= \overline{1}} (\mCM_P) ; \]
  \item[(2)] For $\overline{w}\in W_G/W_R$, if $\overline{w}\neq \overline{1} , \overline{w}_\circ$, then 
  \[\colim_{P\in \Par_G^\st, P\neq G}  \alpha_P'^{= \overline{w}} (\mCM_P)\simeq 0.\]

\end{itemize}
\end{lem}

\begin{lem} \label{lem-colimi-w0}
For $\overline{w}= \overline{w}_\circ$,
  \[\colim_{P\in \Par_G^\st, P\neq G}  \alpha_P'^{= \overline{w}_\circ} (\mCM_P)\simeq 0.\]
\end{lem}

We write the above results as two separate lemmas because they are proved using different ingredients. We will prove the former using the cancellation lemma, while prove the latter using the Deligne--Lusztig duality and Theorem \ref{thm-iota!*-temperedness}.

\proof[Proof of Lemma \ref{lem-colimit-cancel}] 
In this proof, we do not need the assumption that $\mCM_P$ is $G$-anti-tempered.

We assume $\overline{w}\neq \overline{w}_\circ$. Let $-\check \alpha$ be as in Lemma \ref{lem-weyl-combinatorics}, i.e., a negative simple root which is not a weight of $\Lie(R^w)$. Let $\Par^{-\check \alpha\notin}_{G}\subset \Par_{G}^\st\setminus\{G\}$ be the sub-poset containing those $Q$ for which $-\check \alpha$ is not a weight of $\Lie(Q)$. Then there is a retraction
\[\Par_{G}^\st\setminus\{G\}\to \Par^{-\check \alpha\notin}_{G} \]
sending both $Q\langle \alpha\rangle$ and $Q$ to $Q$. 

We claim that the functor
\[ \Par_{G}^\st\setminus\{G\}\to \DMod(\Bun_L) ,\; P\mapsto    \alpha_P'^{= \overline{w}} (\mCM_P) \]
factors through the above retraction. As we now explain, this follows easily from the cancellation lemma. We just need to show that the arrow $Q\to Q\langle \alpha\rangle=:P$ induces an isomorphism 
\[\alpha_Q'^{= \overline{w}} (\mCM_Q)\to  \alpha_P'^{= \overline{w}} (\mCM_P).\]
Suppose that $Q\notin \Par_G^{\underline{w}}$: then $P\notin \Par_G^{\underline{w}}$ by Lemma \ref{lem-weyl-combinatorics}(1), which implies that the source and the target are both $0$ (see Construction \ref{constr-description-gr-dual-weyl}). Now suppose that $Q\in \Par_G^{=\underline{w}}$: then Lemma \ref{lem-weyl-combinatorics}(1) implies that $P \in \Par_G^{=\underline{w}}$, too. By Corollary \ref{cor-condition-cancellation-lemma} and Lemma \ref{lem-weyl-combinatorics}(2), the assumptions of the cancellation lemma (Proposition \ref{prop-cancel-dual-weyl}) are satisfied. Hence the desired morphism is an isomorphism.

By the claim, we have
\[\colim_{Q\in \Par^{-\check \alpha\notin}_{G} }  \alpha_Q'^{= \overline{w}} (\mCM_Q)\simeq \colim_{P\in \Par_G^\st, P\neq G}  \alpha_P'^{= \overline{w}} (\mCM_P).\]
Note that $ \Par^{-\check \alpha\notin}_{G}$ has a final object: the maximal standard parabolic subgroup $S$ such that the Dynkin diagram of $M_{S}$ is obtained from that of $G$ by removing the vertex corresponding to $\check \alpha$. Consequently,
\[ \alpha_S'^{= \overline{w}} (\mCM_S) \simeq \colim_{P\in \Par_G^\st, P\neq G}  \alpha_P'^{= \overline{w}} (\mCM_P). \]

When $\underline{w}\neq \underline{1}$, we claim $S\notin \Par_G^{=\overline{w}}$. 
Of course, this implies (2) because in that case $\alpha_S'^{= \overline{w}}\simeq 0$ (see Construction \ref{constr-description-gr-dual-weyl}). To prove the claim, by Lemma \ref{lem-nonzero-cond-gr}, we only need to show that $B\cap M_S$ is not contained in $R^w$. Suppose the contrary, then $\Lie(R^w)$ contains all the positive simple roots of $M_S$. On the other hand, by definition, $\Lie(R^w)$ contains $\check \alpha$. Hence $R^w$ is a standard parabolic subgroup. Then $B\cap M_S \subset R^w$ implies $R^w = S$ because both $R^w$ and $S$ are maximal proper standard parabolic subgroups, and neither of them contains $-\check \alpha$ as a weight. But this forces $\underline{w}=\underline{1}$, because $S$ is standard. This finishes the proof of (2).

When $\underline{w}= \underline{1}$, we have $S=R$. Hence to prove (1), we need to show the natural transformation $\iota_L^! \to \alpha_R$ factors through an equivalence $\iota_{L}^!\simeq \alpha_R'^{= \overline{1}}$.

Recall (see Construction \ref{constr-description-gr-dual-weyl}) that $(\alpha_R^\vee)^{=\overline{1}}$ is equivalent to 
\[ 
  \begin{aligned}
  \DMod(\Bun_L)  \xrightarrow{*\hpull} \DMod(\Bun_{R^-}) \xrightarrow{!\hpush} \DMod( \widetilde{\Bun}_{R^-} )  \xrightarrow{*\hpush} \DMod(\Bun_G^{R^-\hgen}) \xrightarrow{\Av_*} \mCI(G,R^-),
  \end{aligned}
  \]
which by definition is equivalent to $\iota_{L,!}^-: \DMod(\Bun_L) \to \mCI(G,R^-)$. Passing to duality, we see $\alpha_R'^{= \overline{1}} \simeq \iota_{L}^!: \mCI(G,R) \to \DMod(\Bun_L)$ (see Corollary \ref{cor-inv-inv-duality-functors-global}). It follows from the construction that $\iota_L^! \to \alpha_R$ indeed factors through this equivalence. This finishes the proof of (1).

\qed[Lemma \ref{lem-colimit-cancel}]

It remains to prove Lemma \ref{lem-colimi-w0}. We first need some preparations. Let $S:= (R^-)^{w_\circ}:= w_\circ R^- w_\circ^{-1}$. Note that $S$ is a maximal standard (proper) parabolic subgroup of $G$. Let $-\check \beta$ be the unique negative simple root which is not a weight of $\Lie(S)$.

Recall $\Par_S^\st\subset \Par_G^\st$ is the full subcategory of standard parabolic subgroups $P$ contained in $S$. We claim 
\[\Par_S^\st= \Par_G^{=\overline{w}_\circ}.\]
By Lemma \ref{lem-nonzero-cond-gr}, we only need to show that $B\cap M_P \subset S^-$ iff $P\subset S$, but this is obvious.

Let $[1]$ be the 1-simplex viewed as a category, i.e., it contains two objects $0$ and $1$, and a unique non-invertible arrow $0\to 1$. Then we have an equivalence
\begin{equation} \label{eqn-glu-ParG-from-ParS}
(\Par_{S}^\st\setminus\{S\} \mt [1]) \bigsqcup_{\Par_{S}^\st\setminus\{S\} \mt \{ 0\} } \Par_{S}^\st   \to  \Par_{G}^\st\setminus\{G\}  
\end{equation}
which sends
\begin{itemize}
  \item $(Q,0)\in (\Par_{S}^\st\setminus\{S\} \mt [1])$ to $Q$;
  \item $(Q,1)\in (\Par_{S}^\st\setminus\{S\} \mt [1])$ to $Q\langle \check \beta\rangle$;
  \item $P\in \Par_{S}^\st$ to $P$.
\end{itemize}

By Construction \ref{constr-description-gr-dual-weyl}, if $P\notin \Par_G^{=\overline{w}_\circ}=\Par_S^\st$, then $ \alpha_P'^{= \overline{w}_\circ}\simeq 0$. Hence the equivalence (\ref{eqn-glu-ParG-from-ParS}) implies
\[ \colim_{P\in \Par_G^\st, P\neq G}  \alpha_P'^{= \overline{w}_\circ} (\mCM_P)  \simeq \coFib( \colim_{ Q \in \Par_S^\st, Q\neq S } \alpha_Q'^{= \overline{w}_\circ} (\mCM_Q)  \to \alpha_S'^{= \overline{w}_\circ} (\mCM_S)    )
 \]
By assumption, we have $\mCM_Q\simeq \CT^\enh_{S\gets Q}(\mCM_S)$. Hence it remains to prove that
\begin{equation} \label{eqn-proof-w0-1}
 \coFib( \colim_{ Q \in \Par_S^\st, Q\neq S } \alpha_Q'^{= \overline{w}_\circ} \circ  \CT^\enh_{S\gets Q} \to \alpha_S'^{= \overline{w}_\circ})
 \end{equation}
sends $\mCI(G,S)^{G\hatemp}$ to $0$.

By Construction \ref{constr-dual-weyl-filtration}, the functor (\ref{eqn-proof-w0-1}) is the dual functor of
\begin{equation} \label{eqn-proof-w0-2}
 \coFib( \colim_{ Q \in \Par_S^\st, Q\neq S }   \Eis^\enh_{Q^-\to S^-} \circ (\alpha_Q^\vee)^{= \overline{w}_\circ} \to (\alpha_S^\vee)^{= \overline{w}_\circ}),
 \end{equation}
 where each $(\alpha_Q^\vee)^{= \overline{w}_\circ}$ is equivalent to (see Construction \ref{constr-description-gr-dual-weyl})
 \begin{equation} \label{eqn-alpha-w0}
 \begin{aligned}
\DMod(\Bun_L) \simeq \DMod(\Bun_{M_S}) \xrightarrow{*\hpull} \DMod(\Bun_S) \xrightarrow{!\hpush} \DMod( \widetilde{\Bun}_S ) \to\\ \xrightarrow{!\hpull} \DMod( \widetilde{\Bun}_S^{Q^-\cap S\hgen}) \xrightarrow{*\hpush} \DMod(\Bun_G^{Q^-\hgen}) \xrightarrow{\Av_*} \mCI(G,Q^-).
\end{aligned}
 \end{equation}

We claim the following diagram is Cartesian:
\[
\begin{tikzcd}
   \widetilde{\Bun}_S^{Q^-\cap S\hgen}
    \arrow[r] \arrow[d]
    \arrow[dr, phantom, "\lrcorner", very near start]
    & \Bun_G^{Q^-\hgen} \arrow[d]  \\
      \widetilde{\Bun}_S^{S^-\cap S\hgen}\arrow[r]
    & \Bun_G^{S^-\hgen}.
  \end{tikzcd}
\]
The claim follows from
\[\mBB (Q^-\cap S) \simeq \mBB (S^-\cap S) \mt_{ \mBB S^- } \mBB Q^-,\]
which in turn is a consequence of $(S^-\cap S)Q^- = S^-$.

Applying the base-change isomorphism to this diagram, we see that $(\alpha_Q^\vee)^{= \overline{w}_\circ}$ is equivalent to
 \[ 
 \begin{aligned}
\DMod(\Bun_L) \simeq \DMod(\Bun_{M_S}) \xrightarrow{*\hpull} \DMod(\Bun_S) \xrightarrow{!\hpush} \DMod( \widetilde{\Bun}_S ) \to \\ \xrightarrow{!\hpull} \DMod( \widetilde{\Bun}_S^{M_S\hgen})  \xrightarrow{*\hpush} \DMod(\Bun_G^{S^-\hgen}) 
 \xrightarrow{!\hpull} \DMod(\Bun_G^{Q^-\hgen})\xrightarrow{\Av_*} \mCI(G,Q^-).
\end{aligned}
 \]
By Definition \ref{defn-Eis-CT-enh}, we have a commutative diagram
\[
\xymatrix{
  \DMod(\Bun_G^{S^-\hgen}) \ar[r]^-{!\hpull} \ar[d]^-{\Av_*} &
   \DMod(\Bun_G^{Q^-\hgen})  \ar[d]^-{\Av_*} \\
   \mCI(G,S^-) \ar[r]^-{\CT^\enh_{S^-\gets Q^-}} &
   \mCI(G,Q^-).
}
\]
We obtain
\[ (\alpha_Q^\vee)^{= \overline{w}_\circ}\simeq  \CT^\enh_{S^-\gets Q^-}  \circ (\alpha_S^\vee)^{= \overline{w}_\circ}\]
and can thus rewrite (\ref{eqn-proof-w0-2}) as
\[ \coFib(  \colim_{ Q \in \Par_S^\st, Q\neq S }   \Eis^\enh_{Q^-\to S^-} \circ \CT^\enh_{S^-\gets Q^-} \to \Id ) \circ (\alpha_S^\vee)^{= \overline{w}_\circ} . \]
Now consider the endo-functor of $\mCI(G,S^-)$ defined by
\begin{equation} \label{eqn-dl-s}
\mathrm{DL}_{S^-}:= \coFib(  \colim_{ Q \in \Par_S^\st, Q\neq S }   \Eis^\enh_{Q^-\to S^-} \circ \CT^\enh_{S^-\gets Q^-} \to \Id ).
\end{equation}
We have reduced Lemma \ref{lem-colimi-w0} (and therefore the entire proof) to the following goal.


\begin{goal} \label{goal-lem-colimi-w0}
For any $\mCF \in \DMod(\Bun_L)$ and $\mCM\in \mCI(G,S)^{G\hatemp}$, we have
\[\langle \mathrm{DL}_{S^-}\circ (\alpha_S^\vee)^{= \overline{w}_\circ}(\mCF),\mCM \rangle \simeq 0,\]
where $\langle -,-\rangle$ is the pairing functor between $\mCI(G,S^-)$ and $\mCI(G,S)$. 
\end{goal}

We first deduce this from the following lemmas.

\begin{lem} \label{lem-alpha-is-*-gen}
Let $\mCI(G,S^-)^{\iota_*\hgen}\subset \mCI(G,S^-)$ be the full subcategory generated under colimits by the image of $\iota_{M_S,*}^-: \DMod(\Bun_{M_S}) \to \mCI(G,S^-)$ (see Definition \ref{defn-iota-*-gen}). Then the image of $(\alpha_S^\vee)^{= \overline{w}_\circ}$ is contained in $\mCI(G,S^-)^{\iota_*\hgen}$.
\end{lem}

\begin{lem} \label{lem-*-gen-and-DL}
The image of the functor 
\[ \mCI(G,S) \xrightarrow{\mathrm{DL}_S} \mCI(G,S) \xrightarrow{\iota_{M_S}^*} \DMod(\Bun_{M_S})  \]
is contained in $\DMod(\Bun_{M_S})^{*\hgen}$.
\end{lem}

\proof[Proof of Lemma \ref{lem-colimi-w0}] By Goal \ref{goal-lem-colimi-w0} and Lemma \ref{lem-alpha-is-*-gen}, it suffices to prove that
\[\langle \mathrm{DL}_{S^-}\circ \iota_{M_S,*}^-(\mCK),\mCM \rangle \simeq 0,\]
for any $\mCK \in \DMod(\Bun_{M_S})$ and $\mCM\in \mCI(G,S)^{G\hatemp}$. By Corollary \ref{cor-inv-inv-duality-functors-global} and Remark \ref{rem-inv-inv-duality-functors-global}, the functor $ \mathrm{DL}_{S^-}\circ \iota_{M_S,*}^-$ is dual to $  \iota_{M_S}^*\circ \mathrm{DL}_{S}$. Hence we need to prove that
\[   \iota_{M_S}^*\circ \mathrm{DL}_{S}(\mCM) \simeq 0.\]
By Lemma \ref{lem-*-gen-and-DL}, 
\[\iota_{M_S}^*\circ \mathrm{DL}_{S}(\mCM)\in \DMod(\Bun_{M_S})^{*\hgen} \]
and therefore is $M_S$-tempered (\cite[Theorem B]{beraldo2021geometric}). 

On the other hand, the functor $\mathrm{DL}_{S}$ is $\Sph_{G,x}$-linear, hence $ \mathrm{DL}_{S}(\mCM)$ is $G$-anti-tempered because so is $\mCM$. Theorem \ref{thm-iota!*-temperedness}(4) implies that $\iota_{M_S}^*\circ \mathrm{DL}_{S}(\mCM)$ is also $M_S$-anti-tempered.

Combining the above two paragraphs, we see $\iota_{M_S}^*\circ \mathrm{DL}_{S}(\mCM) \simeq 0$ as desired.

\qed[Lemma \ref{lem-colimi-w0}]

It remains to prove Lemma \ref{lem-alpha-is-*-gen} and Lemma \ref{lem-*-gen-and-DL}.

\proof[Proof of Lemma \ref{lem-alpha-is-*-gen}] By (\ref{eqn-alpha-w0}), the functor $(\alpha_S^\vee)^{= \overline{w}_\circ}$ factors through 
\[ \DMod(\widetilde{\Bun}_S^{M_S\hgen}) \xrightarrow{*\hpush} \DMod(\Bun_G^{S^-\hgen}) \xrightarrow{\Av_*} \mCI(G,S^-). \]
We claim that the image of above functor already lands in $\mCI(G,S^-)^{\iota_*\hgen}$. Let $\lambda\in \Lambda_{G,S}$ and consider the corresponding connected component $\widetilde{\Bun}_{S,\lambda}^{M_S\hgen}$. By Corollary \ref{cor-iota-*-gen-description}, we only need to show that
\[  \DMod(\widetilde{\Bun}_{S,\lambda}^{M_S\hgen}) \xrightarrow{*\hpush} \DMod(\Bun_G^{S^-\hgen}) \xrightarrow{\Av_*} \mCI(G,S^-) \xrightarrow{(\iota_{M_S}^{-,\mu})^!} \DMod(\Bun_{M_S,\mu})\]
is $0$ unless $\mu\ge \lambda$. Here we flip the sign of the inequality in the statement of Corollary \ref{cor-iota-*-gen-description} because we apply it to $S^-$.

By definition, the above composition is given by
\[  \DMod(\widetilde{\Bun}_{S,\lambda}^{M_S\hgen}) \xrightarrow{*\hpush} \DMod(\Bun_G^{S^-\hgen}) \xrightarrow{!\hpull} \DMod(\Bun_{S^-,\mu}) \xrightarrow{*\hpush} \DMod(\Bun_{M_S,\mu}).\]
By the base-change isomorphism, it suffices to check that
\[ \widetilde{\Bun}_{S,\lambda}^{M_S\hgen}\mt_{ \Bun_G^{S^-\hgen} } \Bun_{S^-,\mu} =\emptyset \]
unless $\mu\ge \lambda$. Then we are done: this fiber product is just
\[ (\widetilde{\Bun}_{S,\lambda}\mt_{ \Bun_G } \Bun_{S^-,\mu})^{M_S\hgen}, \]
i.e. the parabolic Zastava space defined in \cite[\S 2]{braverman2002intersection}, and the desired inequality was proved there.

\qed[Lemma \ref{lem-alpha-is-*-gen}]

\proof[Proof of Lemma \ref{lem-*-gen-and-DL}] We will show that the image of
\begin{equation} \label{eqn-proof-lem-*-gen-and-DL-1}
 \DMod(\Bun_{M_S})\xrightarrow{\iota_{M_S,!}} \mCI(G,S) \xrightarrow{\mathrm{DL}_S} \mCI(G,S) \xrightarrow{\iota_{M_S}^*} \DMod(\Bun_{M_S}) 
 \end{equation}
is contained in $\DMod(\Bun_{M_S})^{*\hgen}$. 

Consider the functor $\mathrm{DL}_{M_S}: \DMod(\Bun_{M_S}) \to\DMod(\Bun_{M_S}) $ in Theorem \ref{thm-DL} (applied to the reductive group $M_S$). Using notations in \S \ref{ssec-IGP-functorial}, it is a certain colimit of functors
\[ \mCI(S,S) \xrightarrow{\CT} \mCI(S,Q) \xrightarrow{\Eis} \mCI(S,S) . \]
On the other hand, by definition (see (\ref{eqn-dl-s}), the functor $\mathrm{DL}_S$ is a similar colimit of functors
\[ \mCI(G,S) \xrightarrow{\CT} \mCI(G,Q) \xrightarrow{\Eis} \mCI(G,S) . \]
By the base-change isomorphism for $\mCI(-,-)$ (Lemma \ref{lem-Eis-CT-adjointable}), it is easy to see that
\[ \iota_{M_S}^!\circ \mathrm{DL}_S \simeq \mathrm{DL}_{M_S}\circ \iota_{M_S}^!.  \]
Passing to dualities (Corollary \ref{cor-inv-inv-duality-functors-global} and Remark \ref{rem-inv-inv-duality-functors-global}), we obtain that
\[ \mathrm{DL}_S\circ \iota_{M_S,!} \simeq \iota_{M_S,!}\circ \mathrm{DL}_{M_S}.  \]
Hence the functor (\ref{eqn-proof-lem-*-gen-and-DL-1}) is equivalent to $\mathrm{DL}_{M_S}$. By Theorem \ref{thm-DL}, it suffices to show that the image of the functor $\PsId_{M_S,*}$ is contained in $\DMod(\Bun_{M_S})^{*\hgen}$. But this is known, see \cite[\S 4.4.3]{drinfeld2013some}.

\qed[Lemma \ref{lem-*-gen-and-DL}]

\qed[Theorem \ref{thm-main-thm-atemp}]

\appendix
\section{Lax \texorpdfstring{$\Sph_G$}{SphG}-linear functors are strict}
\label{sect-lax-linear}
\setcounter{subsection}{1}
In this appendix, we prove the following result:

\begin{prop} \label{prop-lax-linear-is-strict-SphG}
Any \emph{continuous} (left or right)-lax $\Sph_G$-linear is strictly $\Sph_G$-linear.
\end{prop}

\begin{defn} \label{defn-rigid}
Following \cite[Chapter 1, Definition 9.1.2]{gaitsgory2019study}, a stable presentable monoidal category $\mCA$ is \emph{rigid} if the following conditions hold:
\begin{itemize}
	\item[(i)] The right adjoint $m^R: \mCA\to \mCA\ot\mCA$ of the multiplication functor $m$ is continuous;
	\item[(ii)] The functor $m^R$, which a priori is right-lax $(\mCA,\mCA^\rev)$-linear, is strictly linear;
	\item[(iii)] The unit object $1_\mCA$ is compact.
\end{itemize}
\end{defn}

\begin{prop}[\!\!{\cite[Chaper 1, Lemma 9.3.6]{gaitsgory2019study}}] \label{prop-lax-linear-is-strict-rigid}
For a presentable monoidal category $\mCA$ satisfying Condition (i) and (ii) in Definition \ref{defn-rigid}, any continuous lax $\mCA$-linear functor is strictly linear\footnote{Warning: the proof in \emph{loc.cit.} contains \emph{lots} of typos.}.
\end{prop}

\qed[Proposition \ref{prop-lax-linear-is-strict-rigid}]
 
\begin{rem} \label{rem-quasi-rigid-is-enough}
In fact, \cite{gaitsgory2019study} stated the above results for rigid monoidal categories. However, the proofs there did not use Condition (iii), nor the stable assumption on $\mCA$.
\end{rem}

\begin{rem}
$\Sph_G$ is \emph{not} rigid because its unit object is not compact. Indeed, the unit object is the pushforward of $\omega$ along $ \mBB G(\mCO) \to G(\mCO)\backslash G(\mCK)/G(\mCO)$ and  $\omega\in \DMod(\mBB G(\mCO) ) \simeq \DMod(\mBB G )$ is not compact (see \cite[\S 7.1.4]{drinfeld2013some}). 
\end{rem}

Nevertheless, we have:

\begin{lem} \label{lem-sphG-quasi-rigid}
The monoidal category $\Sph_G$ satisfies Condition (i) and (ii) in Definition \ref{defn-rigid}.

\end{lem}

\proof Recall the multiplication functor $m$ is given by $!$-pull-$*$-push along the following correspondence
\[  G(\mCO)\backslash G(\mCK)/G(\mCO) \xleftarrow{u} G(\mCO)\backslash G(\mCK)\mt^{G(\mCO)} G(\mCK)/G(\mCO) \xrightarrow{v}  G(\mCO)\backslash G(\mCK)/G(\mCO) \mt  G(\mCO)\backslash G(\mCK)/G(\mCO),\]
where $u$ is ind-proper and $v$ is pro-smooth and in particular placid (see \cite[\S 4.10]{raskin2015d}). It follows that $m\simeq u_*\circ v^!\simeq u_!\circ v^!$ has a continuous right adjoint $m^R\simeq v_{*,\ren}\circ u^!$, where $ v_{*,\ren}$ is the renormalized pushforward functor defined in \emph{loc.cit.}. Hence we verified Condition (i). Also, Condition (ii) follows from the base-change isomorphisms between such renormalized functors and usual $!$-pullback functors (see \cite[Proposition. 4.11.1]{raskin2015d}).

\qed[Lemma \ref{lem-sphG-quasi-rigid}]

\noindent\emph{Proof of Proposition \ref{prop-lax-linear-is-strict-SphG}.}
Now Proposition \ref{prop-lax-linear-is-strict-SphG} follows from Lemma \ref{lem-sphG-quasi-rigid} and Proposition \ref{prop-lax-linear-is-strict-rigid}.

\qed[Proposition \ref{prop-lax-linear-is-strict-SphG}]

We also need the following technical lemma.

\begin{lem}[\!\!{\cite[Chapter 1, Lemma 9.3.2]{gaitsgory2019study}.}] \label{lem-action-functor-quasi-rigid}
For a presentable monoidal category $\mCA$ satisfying Condition (i) and (ii) in Definition \ref{defn-rigid} and any left $\mCA$-module $\mCM$, the action functor $a:\mCA\ot \mCM\to \mCM$ has a continuous right adjoint $a^R$, which is canonically isomorphic to
\[ \mCM \xrightarrow{1_\mCA} \mCA\ot \mCM \xrightarrow{m^R\ot \Id} \mCA\ot\mCA\ot\mCM \xrightarrow{\Id \ot a} \mCA\ot \mCM. \]
\end{lem}

\qed[Lemma \ref{lem-action-functor-quasi-rigid}]

\bibliography{mybiblio.bib}{}
\bibliographystyle{alpha}

\end{document}